\documentclass[11pt,a4paper]{article}
 \usepackage{indentfirst, latexsym,bm}
 \usepackage{amsmath}
 \usepackage{mdframed}
 \usepackage{pifont}
 \usepackage{amsfonts}
 \usepackage{mathrsfs}
 \usepackage{array}
 \usepackage{multirow} 
 \usepackage{graphicx}
 \usepackage{subfigure}
 \usepackage{picinpar}
 \usepackage{setspace}
 \usepackage[textsize=footnotesize]{todonotes}
 \RequirePackage[colorlinks,citecolor=blue,urlcolor=blue,linkcolor=blue]{hyperref}
 \usepackage[top=2cm,bottom=2cm, outer=2cm, inner=2cm]{geometry}
 
 \usepackage{caption}
 \usepackage{subcaption}
 \usepackage[labelfont={bf,small},textfont={small}]{caption}
 
 \usepackage{tikz}
 \usetikzlibrary{arrows.meta}
 
 \RequirePackage[numbers]{natbib}
 \usepackage{enumitem} 
 
 \usepackage{booktabs}
 \usepackage{authblk}
 
 \usepackage{threeparttable}
 \usepackage{mathrsfs,amsfonts,amsmath}
 
 \usepackage{color}
 
 \makeatletter
 \def\namedlabel#1#2{\begingroup
 	#2%
 	\def\@currentlabel{#2}%
 	\phantomsection\label{#1}
 	\endgroup}
 \makeatother

 \numberwithin{figure}{section}

 \newcommand\email[1]{\href{mailto:#1}{ \nolinkurl{#1}}}

 \newtheorem{theorem}{Theorem}[section]
 \newtheorem{definition}[theorem]{Definition}
 \newtheorem{lemma}[theorem]{Lemma}
 \newtheorem{corollary}[theorem]{Corollary}
 \newtheorem{proposition}[theorem]{Proposition}
 \newtheorem{remark}[theorem]{Remark}
 \newtheorem{condition}[theorem]{Condition}
 \newtheorem{example}{Example}[section]

 \def\blemma{\begin{lemma}}\def\elemma{\end{lemma}}
 \def\bproposition{\begin{proposition}}\def\eproposition{\end{proposition}}
 \def\ttheorem{\begin{theorem}}\def\etheorem{\end{theorem}}
 \def\bcorollary{\begin{corollary}}\def\ecorollary{\end{corollary}}
 \def\bremark{\begin{remark}}\def\eremark{\end{remark}}
 \def\bcondition{\begin{condition}}\def\econdition{\end{condition}}

 \def\benumerate{\begin{enumerate}}\def\eenumerate{\end{enumerate}}
 \def\bitemize{\begin{itemize}}\def\eitemize{\end{itemize}}

 \def\beqlb{\begin{eqnarray}}\def\eeqlb{\end{eqnarray}}
 \def\beqnn{\begin{eqnarray*}}\def\eeqnn{\end{eqnarray*}}
 \def\ar{\!\!\!&}

 \def\proof{\noindent{\it Proof.~~}}\def\qed{\hfill$\Box$\medskip}

\begin{document} 
  \title{\bf  Local Limit Theorems for Heavy-Tailed L\'evy Processes with Negative Drift}
 
 \author{Hao Wu\footnote{School of Mathematical Sciences, Nankai University, China; email: wuhao.math@outlook.com} 
 	\quad\ and\quad 
 	Wei Xu\footnote{School of Mathematics and Statistics, Beijing Institute of Technology,  China; email: xuwei.math@gmail.com. Xu gratefully acknowledges financial support from the Fundamental Research Funds for the Central Universities, the National Natural Science Foundation of China (No. 11531001) and the National Key R\&D Program of China (No. 2023YFA1010103).} 
 }  
 \maketitle

 \begin{abstract} 
 
  This paper establishes uniform local asymptotics for the first-passage time below zero and the stopped maximum of a L\'evy process with negative drift and a locally regularly varying L\'evy measure. 
  We show that their probabilities of lying in distant intervals are asymptotic to the expected first-passage time multiplied by the corresponding L\'evy masses, uniformly over interval widths bounded away from zero, including unbounded intervals. 
  Although the two local events have asymptotically equivalent probabilities, their conditional overlap vanishes for fixed-width intervals. Nevertheless, we establish conditional concentration of the first-passage time and stopped maximum around their deterministic linear relation on scales dominating the fluctuations of the centered process. 
  We further prove conditional functional limit theorems yielding a single-step process with random jump time on the original time scale and a deterministic linear descent on the accelerated scale. 
  Finally, we identify the joint limiting distribution of the unscaled pre-passage position and terminal overshoot, distinguishing continuous passage from passage by a downward jump.
 
 	\bigskip
 	
 	\noindent {\it MSC 2020 subject classifications:} Primary 60G51, 60F17;  secondary 60F10, 60G50. 
 	
 	\smallskip
 	
 	\noindent  {\it Keywords and phrases:} L\'evy process,  local limit theorem, regular variation, one-big-jump principle.

 \end{abstract} 
 
 \section{Introduction and main results} 
  \label{Sec.Introduction}
  \setcounter{equation}{0}
  
 The asymptotic behaviour of a stochastic process conditioned to remain in a half-line for a long time depends markedly on the mechanism by which it avoids the boundary. 
 In the stable-domain-of-attraction regime, the path typically exhibits fluctuations on the natural scaling of the underlying process, and the associated conditional limits are closely related to stable meanders and excursions. 
 For a L\'evy process with negative drift and heavy positive jumps, a different picture emerges: survival over a long time interval is typically sustained by a single exceptionally large jump that compensates for the accumulated negative drift. 
 The distinction becomes particularly delicate when a tail event is replaced by a local event, since the latter detects fluctuations that are invisible in the leading-order tail asymptotics.
  
  The purpose of this paper is to establish a uniform local theory for the first-passage time and the maximum before first passage of a heavy-tailed L\'evy process, and to identify the corresponding conditional path and terminal-passage limits. 
  Our results complement the excursion-level asymptotics recently developed in \cite{CuiWuXu2025,CuiWuXu2026b,WuXu2026b}. 
  In particular, they show that the first-passage time and the pre-passage maximum have the same leading-order local masses, even though their fixed-width conditioning events become asymptotically disjoint in a conditional sense.

  \subsection{Background and related results}
  
  There is an extensive literature on first-passage problems and conditional limit theorems for random walks and L\'evy processes. 
  For random walks in the attraction domain of a stable law, Vatutin and Wachtel~\cite{VatutinWachtel2009} and Doney~\cite{Doney2012} obtained the precise local first-passage asymptotics in different regimes. 
  The analogous continuous-time counterparts were investigated by Doney and Rivero~\cite{DoneyRivero2013,DoneyRivero2016Erratum}, who derived local estimates for first-passage times and excursion lifetimes in the stable-domain-of-attraction setting, distinguishing continuous from discontinuous passage.  
  
  A complementary line of research concerns the heights of excursions of L\'evy processes reflected at their past infimum under exponentially decaying tail assumptions. 
  Doney and Maller~\cite{DoneyMaller2005} obtained a Cram\'er-type asymptotic for the excursion-height tail and a Poisson limit theorem for the number of high excursions. 
  Later, Doney and Griffin~\cite{DoneyGriffin2018} revisited these results, corrected an issue in the original proof of the Cram\'er estimate, and established further relationships between excursion-height tails, upward first-passage probabilities and convolution-equivalent L\'evy measures beyond the Cram\'er case. 
  
  The heavy-tailed negative-drift regime requires different techniques. Denisov and Shneer~\cite{DenisovShneer2013} established first-passage tail asymptotics for L\'evy processes and random walks in terms of one-time large-deviation probabilities. For random walks with negative drift, Denisov et al.~\cite{DenisovVatutinWachtel2014} studied local probabilities under the constraint of staying nonnegative and obtained local first-passage asymptotics under suitable regularity and moment hypotheses. 
  Related local estimates for the maximum attained during a cycle were proved in \cite{DenisovShneer2007}. 
   
  At the level of an unrestricted process, the big-jump domain for random walks was systematically investigated by Denisov et al.~\cite{DenisovDiekerShneer2008}. 
  Local subexponentiality is particularly relevant when a probability is evaluated in a prescribed spatial interval rather than in a half-line; see e.g. \cite{AsmussenFossKorshunov2003}. 
  More recently, a continuous-time uniform local large-deviation principle was established in \cite{WuXu2026a} for L\'evy processes with subexponential jumps.  
  
  The conditional path structure was first studied for heavy-tailed random walks by Durrett~\cite{Durrett1980} and later developed for random walks and L\'evy processes by Doney and Jones~\cite{DoneyJones2012}. 
  In particular, the occurrence of a single dominant positive jump produces a discontinuous limiting path with an essentially deterministic descent. 
  In the L\'evy setting, Xu~\cite{Xu2021a} obtained conditional functional limits for negative-drift processes under long survival. 
  The excursion formulation in \cite{CuiWuXu2025} further identifies the asymptotic equivalence of the lifetime and height tail events and the random arrival time of the distinguished large jump.
 
  The subsequent development of the local theory provides a finer description of these phenomena. 
 The uniform local large-deviation asymptotics were established in \cite{CuiWuXu2026b} for the entrance laws of reflected L\'evy excursions and the associated killed semigroup. 
 Meanwhile, Wu and Xu~\cite{WuXu2026b} obtained the uniform local asymptotics for excursion lifetimes and heights, showed that the corresponding fixed-width events have asymptotically negligible conditional overlap, and derived conditional path and terminal-passage limits. 
 Extending these results to a L\'evy process started from a fixed positive point and killed upon its first entrance into the negative half-line is not merely a formal consequence of excursion theory. 
 Indeed, the stopped maximum involves the entire path before killing, while the local first-passage constraint requires precise control of both the initial segment and the terminal crossing. 
 These difficulties are addressed in this paper by combining excursion decomposition, compensation formulas, and fluctuation-theoretic identities.

  \subsection{Main results}

  Consider a real-valued L\'evy process $X=\{ X_t:t\geq 0 \}$ with negative drift $ \beta:=-\mathbf{E}\big[X_1\big]\in (0,\infty)$, diffusion coefficient $\sigma\geq 0$ and L\'evy measure $\nu$ satisfying
  \beqnn
  \int_{\mathbb{R}} \big( |y|\wedge |y|^2\big)\,\nu(dy)<\infty. 
  \eeqnn
  For $x\in \mathbb{R}$, we write $\tau_x^-$ and $\tau_x^+$ for the first passage times of $X$ into $(-\infty, x)$ and $(x,\infty)$, respectively. 
   Denote by $\underline{X}$ and $\overline{X}$ the \textsl{running infimum} and \textsl{supremum processes} defined  by
  \beqnn
  \underline{X}_t:=\inf\big\{ X_s:0\leq s\leq t \big\} 
  \quad \mbox{and}\quad 
  \overline{X}_t:=\sup\big\{ X_s:0\leq s\leq t \big\},\quad t\geq 0.
  \eeqnn
  We write $\mathbf P_x$ and $\mathbf E_x$ for the law and expectation of $X$ started from $x$, and put $\mathbf P=\mathbf P_0$ and $\mathbf E=\mathbf E_0$. For a fixed initial position $x>0$, the main objects of interest are the first-passage time $\tau_0^-$, the stopped maximum $\overline X_{\tau_0^-}$, and the positions immediately before and after passage, $(X_{\tau_0^--},X_{\tau_0^-})$. Throughout, we assume that $0$ is regular for $(0,\infty)$ for $X$.

 Given a constant $\kappa\in\mathbb R$, a positive measurable function $f$ on $\mathbb{R}$ is said to be \textsl{regularly varying} at infinity with index $\kappa$, written $f \in \mathrm{RV}^\infty_\kappa$, if $f(\lambda t)/f(t) \to \lambda^\kappa$ as $t\to\infty$ for any $\lambda>0$.   
 We write $f\sim g$ whenever $f(x)/g(x)\to1$ in the asymptotic regime under consideration. 
 For $\delta\in(0,\infty]$ and $a>0$, let 
  \beqnn
 \varDelta_\delta:=(0,\delta],\qquad
 x+\varDelta_\delta:=(x,x+\delta],\qquad
 a(x+\varDelta_\delta):=(ax,ax+a\delta], 
 \eeqnn
 where $x+\varDelta_\infty=(x,\infty)$.
   The following local regular-variation condition is imposed throughout. It strengthens ordinary regular variation of the positive tail by controlling the L\'evy mass in every fixed interval.

  \begin{condition}\label{Main.Condition}
  	For some $\alpha>1$, assume that $x\mapsto\nu(x+\varDelta_\delta)$ belongs to $\mathrm{RV}^\infty_{-\alpha-1}$ for every $\delta>0$.
  \end{condition}
 
  Set $\overline\nu(x):=\nu\big((x,\infty)\big)$ for $x>0$. 	
  By Proposition~3.3 of \cite{WuXu2026a}, Condition~\ref{Main.Condition} is equivalent to the requirement that $\overline\nu\in\mathrm{RV}_{-\alpha}^\infty$ and, for every fixed $\delta>0$,
  	\beqlb\label{eqn.30}
  	\nu(x+\varDelta_\delta)\sim \alpha\cdot \frac{\overline\nu(x)}{x}\cdot \delta,\quad \mbox{as } x\to\infty.
  	\eeqlb
  	In particular, the assumption allows the L\'evy measure to have a singular component; no density of $\nu$ is required. A fixed-window approximation such as \eqref{eqn.30} cannot, in general, be used uniformly over arbitrarily large window widths. This explains why the full interval mass is retained in the principal results below.

  Our first theorem concerns the local distribution of the first-passage time. It is the negative-mean, heavy-tailed counterpart of local first-passage estimates in the stable-domain-of-attraction regime. The appearance of $\mathbf E_x[\tau_0^-]$ reflects the total time during which the path can experience the exceptional positive jump.
  
  \begin{theorem}\label{MainThm01}
  	For any fixed $x,\delta_0>0$, we have as $t\to\infty$, uniformly in $\delta\in[\delta_0,\infty]$,
  	\beqlb\label{eqn.MainThm01}
  	\mathbf P_x\big(\tau_0^-\in t+\varDelta_\delta\big)
  	\sim \mathbf E_x[\tau_0^-]\cdot \nu\big(\beta(t+\varDelta_\delta)\big).
  	\eeqlb
  \end{theorem}
  
  Our second result gives the corresponding local estimate for the maximum attained before the first passage below zero. The same coefficient appears, although the events in the two theorems constrain different functionals of the stopped path.
  
  \begin{theorem}\label{MainThm02}
  	For any fixed $x,\delta_0>0$,  we have as $h\to\infty$, uniformly in $\delta\in[\delta_0,\infty]$,
  	\beqnn
  	\mathbf P_x\big(\overline X_{\tau_0^-}\in h+\varDelta_\delta\big)
  	\sim \mathbf E_x[\tau_0^-]\cdot \nu\big(h+\varDelta_\delta\big).
  	\eeqnn
  \end{theorem}
  
  Taking $\delta=\infty$ in Theorems~\ref{MainThm01} and \ref{MainThm02} recovers the corresponding tail relations. For fixed $\delta>0$, the local regular variation assumption gives the asymptotic scale $\alpha\beta\delta\,\mathbf E_x[\tau_0^-]\overline\nu(\beta t)/(\beta t)$. However, asymptotic equality of local probabilities does not imply that the events themselves coincide with high conditional probability. Our next theorem makes this distinction precise.
  
  \begin{theorem}\label{MainThm03}
  	For any fixed $x,\delta_1,\delta_2>0$, we have as $t\to\infty$,
  	\beqnn
  	\mathbf P_x\Big(\overline X_{\tau_0^-}\in\beta(t+\varDelta_{\delta_1})\,\Big|\,
  	\tau_0^-\in t+\varDelta_{\delta_2}\Big)
  	\sim\frac{\delta_1}{\delta_2}
  	\mathbf P_x\Big(\tau_0^-\in t+\varDelta_{\delta_2}\,\Big|\,
  	\overline X_{\tau_0^-}\in\beta(t+\varDelta_{\delta_1})\Big)\to 0.
  	\eeqnn
  \end{theorem}
  
  The failure of fixed-window equivalence is a small-scale phenomenon. The first-passage time and the stopped maximum remain strongly coupled once their deviations are measured on a scale dominating the typical fluctuations of $X_t+\beta t$.
  
 \begin{theorem}\label{MainThm04}
	For any  $x, \delta>0$ and any positive function $\theta$ on $\mathbb{R}_+$ such that
	$|X_t+\beta t|/\theta_t \overset{\rm p}\to0$, we have as $t\to\infty$,
	\beqlb\label{eqn.MainThm04}
	\mathbf P_x\Big(
	\big|\overline{X}_{\tau_0^{-}}-\beta t\big|\leq\theta_t
	\,\Big|\,
	\tau_0^-\in t+\varDelta_\delta
	\Big)
	\to1
	\quad\mbox{and}\quad
	\mathbf P_x\Big(
	|\tau_0^--t|\leq\theta_t
	\,\Big|\,
	\overline{X}_{\tau_0^{-}}\in\beta(t+\varDelta_\delta)
	\Big)
	\to1.
	\eeqlb
\end{theorem}
  
  The next result describes the entire conditioned path. Let $\mathcal T_x$ be a nonnegative random variable with distribution function
  \beqlb\label{eqnn1.9}
  \mathbf P(\mathcal T_x\leq s)
  =\frac{\mathbf E_x[\tau_0^-\wedge s]}{\mathbf E_x[\tau_0^-]},\qquad s\geq0.
  \eeqlb
 The random time $\mathcal T_x$ describes the asymptotic arrival time of the distinguished large jump. Under local conditioning, its rescaled size becomes deterministic, in contrast to the Pareto-distributed jump size under tail conditioning. 
 For a locally compact interval $I$, the space $D(I;\mathbb R)$ is equipped with the Skorokhod $J_1$ topology, with convergence on compact subintervals when $I$ is unbounded.
 We also introduce the class 
  \beqnn
  \boldsymbol{\varTheta}:=\big\{\theta:\mathbb R_+\to(0,\infty):\theta\text{ is nondecreasing},\ \theta_t\to\infty,\ \theta_t/t\to0\big\}.
  \eeqnn
  
  \begin{theorem}\label{MainThm05}
  	For any $x>0$ and $\theta\in\boldsymbol{\varTheta}$, as $t\to\infty$, uniformly in $\delta\in[\theta_0,\theta_t]$, the following assertions hold under either of the conditional laws $\mathbf P_x\big(\,\cdot\mid\tau_0^-\in t+\varDelta_\delta\big)$ or
  	$\mathbf P_x\big(\,\cdot\mid\overline X_{\tau_0^-}\in\beta(t+\varDelta_\delta)\big)$.
  	\begin{enumerate}
  		\item[\rm(1)]  $\{X_s/t:s\geq0\}\overset{\rm d}\to \{\beta\cdot \mathbf 1_{\{s\geq\mathcal T_x\}}:s\geq0\}$ in $D\big([0,\infty);\mathbb R\big)$. 
  		\item[\rm(2)]  $	\{X_{ts}/t:s>0\}\overset{\rm d}\to \{\beta(1-s):s>0\}$ in $D\big((0,\infty);\mathbb R\big)$. 
  	\end{enumerate}
  \end{theorem}
  
  The two scales capture complementary aspects of the same path. On the original time scale, only the location of the exceptional jump remains random; on the accelerated scale, this jump occurs at time zero in the limit and the negative mean produces a deterministic linear descent. Since the original process is not killed after $\tau_0^-$ in this theorem, the second limiting path continues below zero for $s>1$.
  
  Our final theorem resolves the terminal crossing on its unscaled spatial scale. Let $V$ be the renewal function of the ascending ladder-height process and let $\widehat H$ be the descending ladder-height process, with the paired normalization specified in Section~\ref{Sec.Preliminaries}. Write $\underline{\tt a}$ for the drift coefficient of $\widehat H$.  Theorem~7.8 in \cite[p.~213]{Kyprianou2014} gives
  \beqlb\label{eqn1.9}
  \mathbf E[\widehat H_1]
  =\underline{\tt a}+\int_0^\infty V(z)\nu(( -\infty,-z))\,dz. 
  \eeqlb
  
  \begin{theorem}\label{MainThm06}
  	For any fixed $x,\delta_0>0$ and any bounded continuous function $f$ on $\mathbb R_+^2$, as $t\to\infty$, uniformly in $\delta\in[\delta_0,\infty]$,
  	\beqlb\label{eqn.MainThm06}
  	\mathbf E_x\Big[f\big(X_{\tau_0^--},-X_{\tau_0^-}\big)\,\Big|\,
  	\tau_0^-\in t+\varDelta_\delta\Big]
  	\to \frac{1}{\mathbf E[\widehat H_1]}
  	\Big(\underline{\tt a}\,f(0,0)
  	+\int_0^\infty du\int_0^\infty f(u,v)V(u)\,\nu(-u-dv)\Big).
  	\eeqlb
  \end{theorem}
  
  The atom at $(0,0)$ in the limiting measure corresponds to continuous first passage, while the second term describes first passage by a negative jump. This limit is particularly informative because the preceding functional limits do not retain the $O(1)$ positions at the terminal crossing. It also gives a common limiting terminal mechanism under large local first-passage conditioning and, at the excursion level, under large local lifetime conditioning; compare \cite{WuXu2026b}.
  
  \smallskip
  {\it\textbf{Organization of this paper.}}
  Section~\ref{Sec.Preliminaries} recalls the fluctuation and excursion identities used throughout.
 Sections~\ref{Sec.ProofThm01} and~\ref{Sec.Thm02} establish the local first-passage and maximum estimates. Sections~\ref{Sec.Thm03} and~\ref{Sec.Thm04} study the relation between the two local conditioning events; Sections~\ref{Sec.ProofThm05} and~\ref{Sec.ProofThm06} prove the conditional functional limits and the terminal-passage limit, respectively.

   \medskip
  
  \textit{\textbf{Notation.}} 
  Let $\mathbb{N}$ be the space of all natural numbers including $0$. 
  For any $x\in\mathbb{R}$, let $x^+:=x\vee 0$, $x^-=-(x\wedge 0)$ and $[x]$ be the integer part of $x$. 
  Let $\overset{\rm f.d.d.}\to $,  $\overset{\rm a.s.}\to$, $\overset{\rm d}\to$  and  $\overset{\rm p}\to$ be the convergence in the sense of finite dimensional distributions, almost sure convergence,  convergence in distribution and convergence in probability, respectively.  
  We make the conventions that for $x, y\in \mathbb{R}$ with $y\geq x$,
  \beqnn
  \int_x^y=-\int_y^x= \int_{(x,y]},\quad \int_{x-}^{y-}=  \int_{[x,y)}
  \quad\mbox{and}\quad
  \int_x^\infty = \int_{(x,\infty)}.
  \eeqnn
  Integrals with respect to an excursion measure are written in expectation notation, for example $\underline n(\zeta)$ and $\underline n(\zeta\wedge t)$. 
  Without confusion, we also write $A\triangle B$ for the symmetric difference of the two sets $A$ and $B$. 
  We use $C$ to denote a positive constant whose value might change from line to line.

   \section{Preliminaries} 
  \label{Sec.Preliminaries}
  \setcounter{equation}{0} 
  
  We collect the fluctuation-theoretic and excursion-theoretic identities of L\'evy processes needed in the proofs. 
  Particular attention is paid to the normalization of ladder processes, the two compensation formulas, and the representation of the killed semigroup. 
  We refer to  \cite{Bertoin1996,Doney2007,Kyprianou2014} for more details.

  \subsection{Fluctuation theory}
 
 Assume that $X$ is defined on a complete filtered probability space $(\Omega,\mathscr F,(\mathscr F_t)_{t\geq0},\mathbf P)$ satisfying the usual hypotheses. 
  Let $\Delta X_t := X_t-X_{t-}$ denote the jump size of $X$ at time $t$. 
  The point measure of all jumps 
  \beqnn
  \sum_{t>0} \mathbf{1}_{\{ \Delta X_t\neq 0 \}}\cdot \delta_{(t,\Delta X_t)}(ds,dy),
  \eeqnn
  is a Poisson random measure on $(0,\infty)\times \mathbb{R}$ with intensity $ds\, \nu(dy)$.  
  By the compensation formula for Poisson random measures (see \cite[p.7]{Bertoin1996}), we have that for any predictable process $\mathbf{F}=\{ \mathbf{F}_t:t\geq 0 \}$ taking values in the space of non-negative measurable functions on $\mathbb{R}$, 
  \beqlb\label{eqn.CompenstationX}
  \mathbf{E}\bigg[ \sum_{t>0} \mathbf{F}_t(\Delta X_t) \cdot \mathbf{1}_{\{ \Delta X_t\neq 0 \}}  \bigg]= \mathbf{E}\bigg[\int_0^\infty dt \int_{\mathbb{R}\setminus\{0\}}\mathbf             {F}_t(y)\, \nu(dy)\bigg].
  \eeqlb

 Recall the running infimum and supremum processes $\underline{X}$ and $\overline{X}$. 
 It is well known that the reflected processes $X-\underline X$ and $\overline X-X$ are nonnegative strong Markov processes.
 Let $\widehat L=(\widehat L_t)_{t\geq0}$ be the local time at zero of the reflected process $X-\underline X$. When zero is regular for $(-\infty,0)$, we use the continuous local time; otherwise, we use the right-continuous version constructed by means of independent exponential variables, as in Sections~IV.3 and IV.5 of \cite{Bertoin1996}. 
 Without loss of generality, we normalize it so that 
 \beqnn
 \mathbf{E}\bigg[ \int_0^\infty e^{-t} d \widehat{L}_t\bigg] =1. 
 \eeqnn
 Its right-inverse $\widehat{L}^{-1}$ is referred to as the \textit{descending ladder time process} and $\widehat{H}= \{ \widehat{H}_t :t\geq 0 \}$ with $\widehat{H}_t := -X_{\widehat{L}^{-1}_t}$ is called the \textit{descending ladder height process}, i.e., 
 \beqnn
 \widehat L^{-1}_u:=\inf\{t\geq0:\widehat L_t>u\} 
 \quad \mbox{and}\quad 
 \widehat H_u:=-X_{\widehat L^{-1}_u},\qquad u\geq0,
 \eeqnn
 with $\widehat L^{-1}_{0-}=0$ and $\widehat H_{0-}=0$ for convention.
  Then $(\widehat L^{-1},\widehat H)$ is a possibly killed bivariate subordinator with Laplace exponent $\widehat\kappa(\lambda,q)
  :=-\log\mathbf E\big[e^{-\lambda\widehat L^{-1}_1-q\widehat H_1}\big]$. 
 We write $\underline{\tt d}$ and $\underline{\tt a}$ for the drift coefficients of $\widehat L^{-1}$ and $\widehat H$, respectively. The regularity of zero for $(0,\infty)$ implies $\underline{\tt d}=0$.
  Denote the bivariate potential measure of $(\widehat L^{-1},\widehat H)$ by
 \beqlb\label{eqn.U}
 \widehat U(dt,dy)
 :=\int_0^\infty\mathbf P\big(\widehat L^{-1}_s\in dt,\widehat H_s\in dy\big)\,ds.
 \eeqlb
 The descending ladder-height renewal function is
 \beqlb\label{eqn.DefV}
 \widehat V(x)
 :=\widehat U(\mathbb R_+,[0,x])
 =\int_0^\infty\mathbf P(\widehat H_s\leq x)\,ds
 =\mathbf E[\widehat L_{\tau_{-x}^-}],\qquad x\geq0.
 \eeqlb
 Quantities associated with $\overline{X}$ and $\overline{X}-X$ are denoted without hats, e.g., $\kappa(\lambda,u)$, $U(dt,dy)$ and $V(x)$,  and satisfy analogous properties. 
 The local time and its inverse are normalized jointly so that $\kappa(\lambda,0)\widehat\kappa(\lambda,0)=\lambda$ for every $\lambda\geq0$. 
 Since $X$ drifts to $-\infty$, the ascending ladder process is
 killed, with $\kappa(0,0)>0$, whereas the descending ladder process
 is not killed.

 \subsection{Excursion theory} 
 
 Let $\mathcal E$ denote the space of excursion paths, and let $\Upsilon$ be an isolated cemetery mark.  
 Each excursion $\epsilon \in \mathcal{E}$ is a non-negative c\`adl\`ag function on $\mathbb{R}_+$ with \textsl{lifetime} $\zeta$ and \textsl{height} $\overline\epsilon$ defined by 
 \beqnn
 \zeta:= \inf\big\{t>0: \epsilon_s= 0 \mbox{ for all } s\geq t \big\}
 \quad\mbox{and}\quad 
 \overline\epsilon:=\sup_{0\leq s<\zeta}\epsilon_s, 
 \eeqnn
 with $\inf \emptyset =\infty$ for convention.
 Here and below, $\epsilon_\zeta$ denotes the terminal value in the excursion representation, so that $-\epsilon_\zeta$ is the nonnegative terminal overshoot when the excursion terminates by a downward jump. 
 The excursion process $e=\{e_u: u\geq0\}$ of $X-\underline X$ away from $0$ is the $\mathcal{E}\cup\{\Upsilon\}$-valued process defined by
 \beqnn
 e_u :=
 \begin{cases}
 	\big\{(X-\underline X)_{s+\widehat L^{-1}_{u-}}: 0\leq s<  \widehat L^{-1}_{u}-\widehat L^{-1}_{u-} \big\}, & \text{if } \widehat L^{-1}_{u-}<\widehat L^{-1}_{u}, \\[1.2ex]
 	\Upsilon, & \text{if } \widehat L^{-1}_{u-}=\widehat L^{-1}_{u}.
 \end{cases}
 \eeqnn
 Here $\delta_\Upsilon$ denotes the zero measure on $\mathcal{E}$.  
 It is a Poisson point process on $\mathcal E$ governed by the characteristic measures $\underline n(d\epsilon)$. 
 Similarly,  excursions of  $\overline{X}-  X$ away from $0$ also form a Poisson point process on $\mathcal E$ with characteristic measure $\overline n(d\epsilon)$.


 By the strong Markov property of $X$ and the memoryless property of the independent exponential variables used to construct $\widehat L$; see Section~IV.5 in \cite[pp.122-123]{Bertoin1996}, the joint process $ \boldsymbol{Y}=\{  \boldsymbol{Y}_t:t\geq 0\}$ with 
 \beqnn
 \boldsymbol{Y}_t:=\Big(
 \widehat L^{-1}_t,\widehat H_t,\,
 \sum_{0<u\leq t}\delta_{e_u}
 \Big)  
 \eeqnn
 has stationary independent increments and is strong Markov with respect to its natural filtration $\{\mathscr{G}_t:t\geq 0\}$\footnote{When $0$ is regular for $(-\infty,0)$, we have almost surely  $\widehat L^{-1}_{t-}\leq\widehat L^{-1}_{t}$ and then $\mathscr{G}_t:=\mathscr{F}_{\widehat L^{-1}_t}$ for all $t\geq 0$.} that is a filtration generated by itself and satisfying the usual hypotheses.  
 More precisely, for any finite $(\mathscr{G}_t)$-stopping time $\tau$, Lemma~2(i) in \cite[p.~157]{Bertoin1996}, applied to $-X$, together with the strong Markov property of $X$ at $\widehat L^{-1}_\tau$, yields that $\{  \boldsymbol{Y}_{\tau+t}- \boldsymbol{Y}_\tau:t\geq 0 \}$ under $\mathbf{P}(\cdot\,|\, \mathscr{G}_\tau)$ is equal in law to $ \boldsymbol{Y}$. 
 In particular, the process $e$ is a $(\mathscr{G}_t)$-Poisson point process with characteristic measure $\underline{n}(d\epsilon)$ that is a $\sigma$-finite measure on $\mathcal{E}$; see Theorem~10 in \cite[p.~118]{Bertoin1996}.

 Let $\mathbf{F}=\{\mathbf{F}_t: t\geq0\}$ be a $(\mathscr{G}_t)$-predictable process taking values in the space of nonnegative measurable functions on $\mathcal{E}\cup\{\Upsilon\}$ with $\mathbf{F}_t(\Upsilon)=0$. 
 The compensation formula, see Section~0.5 in \cite[p.~7]{Bertoin1996}, gives
 \beqlb\label{eqnn5.11}
 \mathbf E\bigg[\sum_{t\geq0} \mathbf{F}_t(e_t) \bigg] = \mathbf E\Bigg[ \int_0^\infty dt \int_{\mathcal E} \mathbf{F}_t(\epsilon)\,\underline n(d\epsilon) \Bigg].
 \eeqlb
 In particular, recall the first passage time $\tau_x^-$ of $X$ into $(-\infty,x)$,   similarly as  in the proof of Theorem~19 in \cite[p.~175]{Bertoin1996} we have 
 \beqnn
 \widehat L_{\tau_{-x}^- -}
 = \inf\big\{ t\geq0:\widehat H_t>x \big\}
 \quad\mbox{and}\quad \big\{ t\geq0:\widehat H_{t-}\leq x \big\} = \big[ 0,\widehat L_{\tau_{-x}^- -} \big].
 \eeqnn
 Since $\widehat H$ is $(\mathscr{G}_t)$-adapted and c\`adl\`ag, the first identity shows that $\widehat L_{\tau_{-x}^- -}$ is a $(\mathscr{G}_t)$-stopping time.
 Hence, the mapping  $ t\mapsto \mathbf 1_{\{t\leq\widehat L_{\tau_{-x}^- -}\}} $
 is left-continuous and $(\mathscr{G}_t)$-adapted, and hence $(\mathscr{G}_t)$- predictable.
 Applying \eqref{eqnn5.11} to the $(\mathscr{G}_t)$-predictable process $ t\mapsto \mathbf 1_{\{t\leq\widehat L_{\tau_{-x}^- -}\}}\cdot \mathbf{F}_t $ yields
 \beqlb\label{eqnn5.12}
 \mathbf E\bigg[ \sum_{t\geq0} \mathbf 1_{\{t\leq\widehat L_{\tau_{-x}^- -}\}}\cdot \mathbf{F}_t(e_t) \bigg] 
 \ar=\ar \mathbf E\bigg[ \int_0^{\widehat L_{\tau_{-x}^- -}} dt \int_{\mathcal E} \mathbf{F}_t(\epsilon)\,\underline n(d\epsilon)\bigg].
 \eeqlb

  Under $\underline n$, the coordinate process $\epsilon$ has the transition semigroup of $X$ killed at its first entrance into $(-\infty,0]$, after the initial entrance from zero. In particular, for nonnegative measurable functionals $F$ and $G$ on the indicated path intervals,
 \beqnn
 \lefteqn{\underline n\Big(F(\epsilon_r:0\leq r\leq t)\,
 	G(\epsilon_{t+r}:0\leq r\leq s),\ \zeta>t+s\Big)}\ar\ar\cr\cr
 \ar=\ar\underline n\Big(F(\epsilon_r:0\leq r\leq t)\,
 \mathbf E_{\epsilon_t}\big[G(X_r:0\leq r\leq s),\ \underline X_s>0\big],\ \zeta>t\Big).
 \eeqnn
 Combining the Markov property with the L\'evy system yields an excursion-level compensation formula. With $\mathscr G_s:=\sigma(\epsilon_r:0\leq r\leq s)$, for any nonnegative $(\mathscr G_s)$-predictable family $G_s(y)$,
 \beqlb\label{eqn.Compensation}
 \underline n\bigg[\sum_{0<s\leq\zeta}G_s(\Delta\epsilon_s)
 \mathbf 1_{\{\Delta\epsilon_s\ne0\}}\bigg]
 =\underline n\bigg[\int_0^\zeta ds
 \int_{\mathbb R\setminus\{0\}}G_s(y)\,\nu(dy)\bigg].
 \eeqlb 
 Here, the terminal jump	$\Delta\epsilon_\zeta:=\epsilon_\zeta-\epsilon_{\zeta-}$	is computed using the unreflected terminal position retained	as an additional mark.
 The two compensation formulas \eqref{eqnn5.12} and \eqref{eqn.Compensation} apply at different levels: the former sums over excursions indexed by ladder local time, while the latter sums over jumps within a single excursion.

 The excursion representation of the bivariate exponent  of $(\widehat L^{-1},\widehat H)$ takes the form
  \beqlb\label{eqn.kappa}
  \widehat\kappa(\lambda,q)
  =\underline{\tt a}\,q
  +\int_0^\infty\int_{[0,\infty)}(1-e^{-\lambda t-qz})\,
  \underline n(\zeta\in dt,-\epsilon_\zeta\in dz),\quad \lambda,q\geq0.
  \eeqlb 
  The negative-drift implies $\kappa(0,0)>0$. Setting $q=0$ in \eqref{eqn.kappa} and integrating by parts gives 
  \beqlb\label{eqn.501}
  \underline n(\zeta)
  :=\int_0^\infty\underline n(\zeta>t)\,dt
  =\lim_{\lambda\downarrow0}\frac{\widehat\kappa(\lambda,0)}{\lambda}
  =\frac1{\kappa(0,0)}<\infty.
  \eeqlb
  Moreover, the Wiener--Hopf factorization and the normalization above yield
  \beqlb\label{eqn.EH1}
  \mathbf E[\widehat H_1]
  =\left.\frac{\partial\widehat\kappa(0,q)}{\partial q}\right|_{q=0+}
  =\frac{\beta}{\kappa(0,0)}
  =\beta \cdot \underline n(\zeta).
  \eeqlb
  These identities specify the constants appearing in Theorems~\ref{MainThm01} and \ref{MainThm06}.
  The renewal functions also admit representations in terms of excursion entrance laws. Writing $\overline{\tt d}$ for the drift of $L^{-1}$, we have
  \beqlb\label{eqn.447}
  \widehat V(x)
  =\overline{\tt d}
  +\int_0^\infty\overline n(\epsilon_s\leq x,\zeta>s)\,ds,
  \qquad x\geq0,
  \eeqlb
  see also (16) in \cite{DoneyRivero2013}.
  Similarly, we also have for any $x\geq0$,
  \beqnn
  V(x)
  =
  \int_0^\infty \mathbf P(H_s\leq x)\,ds
  =
  \int_0^\infty
  \underline n(\epsilon_s\leq x,\zeta>s)\,ds
  \quad\mbox{and}\quad
  V(\infty)=\underline n(\zeta)<\infty.
  \eeqnn
 The preceding identities connect the expected first-passage
 time with the descending ladder renewal function. In fact,
 for each $x>0$,  Proposition~17~(ii) in
 \cite[p.172]{Bertoin1996} and Theorem~2.2 in
 \cite{DenisovShneer2013} give
 \beqlb\label{Asytau}
 \mathbf E_x[\tau_0^-]
 =
 \underline n(\zeta)\cdot\widehat V(x)
 \quad\mbox{and}\quad
 \mathbf P_x(\tau_0^->t)
 \sim
 \underline n(\zeta)\cdot\widehat V(x)\cdot
 \overline\nu(\beta t)
 \in\mathrm{RV}_{-\alpha}^{\infty},
 \quad\mbox{as }t\to\infty.
 \eeqlb
   
 By Theorem~5 in \cite{Chaumont2013}, the transition semigroup of $X$
 killed upon entering $(-\infty,0)$ admits the following representation:
 for any $t, x,y>0$,
 \beqnn
 \mathbf P_x\big(X_t\in dy,\tau_0^->t\big)
 \ar=\ar
 \int_0^t ds
 \int_{(x-y)^+}^{x}
 \overline n\big(\epsilon_s\in dz,\zeta>s\big)
 \underline n\big(\epsilon_{t-s}\in dy+z-x,\zeta>t-s\big)
 \cr\cr
 \ar\ar
 +\overline{\tt d}\cdot 
 \underline n\big(\epsilon_t\in dy-x,\zeta>t\big)\cdot 
 \mathbf 1_{\{y\geq x\}}.
 \eeqnn
%
%
%
%
 The measures on the right-side hand are understood to vanish whenever
 their spatial arguments are outside $(0,\infty)$.
 Integrating over $y>0$ and applying Tonelli's theorem gives
 \beqlb\label{eqn.209}
 \mathbf P_x(\tau_0^->t)
 =\overline{\tt d}\cdot\underline n(\zeta>t)
 +\int_0^t\overline n(\epsilon_s\in(0,x],\zeta>s)
 \,\underline n(\zeta>t-s)\,ds.
 \eeqlb
 
%
%
%
%
%
%

 \section{Proofs of Theorem \ref{MainThm01}}
 \label{Sec.ProofThm01}
 \setcounter{equation}{0}
 
 This section is devoted to the proof of the local limit theorem
 for the first-passage time $\tau_0^-$ stated in
 Theorem~\ref{MainThm01}.
 For any $t,\delta>0$, we have
 $\mathbf P_x\big(\tau_0^-\in t+\varDelta_\delta\big)
 =\mathbf P_x(\tau_0^->t)-\mathbf P_x(\tau_0^->t+\delta)$.
 In view of \eqref{eqn.209}, this probability difference admits the following decomposition
 \beqlb\label{eqn.4011}
 \mathbf P_x\big(\tau_0^-\in t+\varDelta_\delta\big)
 \ar=\ar
 \overline{\tt d}\cdot 
 \underline n\big(\zeta\in t+\varDelta_\delta\big)
 +
 \int_0^t 
 \overline n\big(\epsilon_s\leq x,\zeta>s\big)\, 
 \underline n\big(\zeta\in t-s+\varDelta_\delta\big)\, ds
 \cr\cr
 \ar-\ar
 \int_t^{t+\delta} 
 \overline n\big(\epsilon_s\leq x,\zeta>s\big)\, 
 \underline n\big(\zeta>t+\delta-s\big)\, ds.
 \eeqlb
 Consequently, the long-time asymptotics of
 $\mathbf P_x\big(\tau_0^-\in t+\varDelta_\delta\big)$ can be established by analyzing the three terms on the right-hand side.  
 Firstly,  Theorem~1.3 in \cite{WuXu2026b} tells that as $t\to\infty$,
 \beqlb \label{eqn.301}
 \underline n(\zeta\in t+\varDelta_\delta)
 \sim \underline n(\zeta)\cdot
 \nu\big(\beta(t+\varDelta_\delta)\big).
 \eeqlb
 It remains to analyze the last two terms in \eqref{eqn.4011}.  
 For any $\varepsilon\in(0,1/2)$, we first decompose them into three parts:
 \beqlb
 J_1(t,\varepsilon,\delta)
 \ar:=\ar
 \int_0^{\varepsilon t}
 \overline n\big(\epsilon_s\leq x,\zeta>s\big)\,
 \underline n\big(\zeta\in t-s+\varDelta_\delta\big)\,ds, \label{eqn.3031} \\
 J_2(t,\varepsilon,\delta)
 \ar:=\ar
 \int_{\varepsilon t}^{t-\varepsilon t}
 \overline n\big(\epsilon_s\leq x,\zeta>s\big)\,
 \underline n\big(\zeta\in t-s+\varDelta_\delta\big)\,ds, \label{eqn.3032}\\
 J_3(t,\varepsilon,\delta)
 \ar:=\ar
 \int_{t-\varepsilon t}^{t}
 \overline n\big(\epsilon_s\leq x,\zeta>s\big)\,
 \underline n\big(\zeta\in t-s+\varDelta_\delta\big)\,ds
 \cr\cr
 \ar\ar-
 \int_t^{t+\delta}
 \overline n\big(\epsilon_s\leq x,\zeta>s\big)\,
 \underline n\big(\zeta>t+\delta-s\big)\,ds, \label{eqn.3033}
 \eeqlb
 and then investigate their long-term behavior one-by-one in the following three lemmas.
 
 \begin{lemma}\label{Lemma.301}
 	For any $\delta>0$, we have 
 	\beqnn
 	\lim_{\varepsilon\to 0+}\limsup_{t\to\infty}
 	\bigg| \frac{ J_1(t,\varepsilon,\delta) }{ \nu\big(\beta(t+\varDelta_\delta)\big)}- \underline n(\zeta)\cdot  \big(\widehat V(x)-\overline{\tt d}\big)  \bigg|
 	=0.
 	\eeqnn
 \end{lemma} 
 \proof In view of \eqref{eqn.3031} and \eqref{eqn.447}, we first have 
 \beqnn
 \lefteqn{\frac{ J_1(t,\varepsilon,\delta) }{ \nu\big(\beta(t+\varDelta_\delta)\big)} - \underline n(\zeta)\cdot  \big(\widehat V(x)-\overline{\tt d}\big) }\cr\cr 
 \ar=\ar  \int_0^{\varepsilon t} \overline n\big(\epsilon_s\leq x,\zeta>s\big)\,
 \frac{\underline n\big(\zeta\in t-s+\varDelta_\delta\big)}{\nu\big(\beta(t+\varDelta_\delta)\big)}\,ds - \underline n(\zeta) \int_0^\infty
 \overline n\big(\epsilon_s\leq x,\zeta>s\big)\,ds\cr
 \ar=\ar \int_0^{\varepsilon t} \overline n\big(\epsilon_s\leq x,\zeta>s\big)\,\bigg(
 \frac{\underline n\big(\zeta\in t-s+\varDelta_\delta\big)}{\nu\big(\beta(t+\varDelta_\delta)\big)}-\underline n(\zeta)\bigg) \,ds - \underline n(\zeta) \int_{\varepsilon t}^\infty
 \overline n\big(\epsilon_s\leq x,\zeta>s\big)\,ds.
 \eeqnn
 This along with  the triangle inequality yields that 
 \beqnn
 \lefteqn{\bigg| \frac{ J_1(t,\varepsilon,\delta) }{ \nu\big(\beta(t+\varDelta_\delta)\big)}- \underline n(\zeta)\cdot  \big(\widehat V(x)-\overline{\tt d}\big)  \bigg|}\cr\cr
 \ar\leq\ar \int_0^{\varepsilon t} \overline n\big(\epsilon_s\leq x,\zeta>s\big)\,\bigg|
 \frac{\underline n\big(\zeta\in t-s+\varDelta_\delta\big)}{\nu\big(\beta(t+\varDelta_\delta)\big)}-\underline n(\zeta)\bigg| \,ds + \underline n(\zeta) \int_{\varepsilon t}^\infty
 \overline n\big(\epsilon_s\leq x,\zeta>s\big)\,ds.
 \eeqnn
 By \eqref{eqn.447}, the last integral on the right-hand side
 	tends to $0$ as $t\to\infty$ for every fixed $\varepsilon>0$,
 	and the first one is bounded above by
 	\beqlb\label{eqn.303}
 	\lefteqn{
 		\big(\widehat V(x)-\overline{\tt d}\big)\cdot
 		\sup_{s\in[0,\varepsilon t]}
 		\Big|
 		\frac{
 			\underline n\big(\zeta\in t-s+\varDelta_\delta\big)
 		}{
 			\nu\big(\beta(t+\varDelta_\delta)\big)
 		}
 		-\underline n(\zeta)
 		\Big|
 	}
 	\cr\cr
 	\ar\leq\ar
 	\underline n(\zeta)\cdot
 	\big(\widehat V(x)-\overline{\tt d}\big)\cdot
 	\sup_{s\in[0,\varepsilon t]}
 	\Big|
 	\frac{
 		\nu\big(\beta(t-s+\varDelta_\delta)\big)
 	}{
 		\nu\big(\beta(t+\varDelta_\delta)\big)
 	}
 	-1
 	\Big|
 	\cr\cr
 	\ar\ar+
 	\big(\widehat V(x)-\overline{\tt d}\big)\cdot
 	\sup_{s\in[0,\varepsilon t]}
 	\frac{
 		\nu\big(\beta(t-s+\varDelta_\delta)\big)
 	}{
 		\nu\big(\beta(t+\varDelta_\delta)\big)
 	}
 	\cdot
 	\sup_{s\in[0,\varepsilon t]}
 	\Big|
 	\frac{
 		\underline n\big(\zeta\in t-s+\varDelta_\delta\big)
 	}{
 		\nu\big(\beta(t-s+\varDelta_\delta)\big)
 	}
 	-\underline n(\zeta)
 	\Big|.
 	\eeqlb
 By Condition~\ref{Main.Condition} and the uniform convergence theorem for regularly varying functions; see Theorem~1.5.2 in \cite{BinghamGoldieTeugels1987},  we have as $t\to\infty$, 
 \beqlb\label{eqn.313}
 \sup_{s\in[0,\varepsilon t]}
 \frac{\nu\big(\beta(t-s+\varDelta_\delta)\big)}{\nu\big(\beta(t+\varDelta_\delta)\big)} \to (1-\varepsilon)^{-\alpha-1}.
 \eeqlb
 This induces that  the first term on the right-hand side of \eqref{eqn.303} vanishes as $t\to \infty$ and then $\varepsilon \to 0+$, and also that the second term can be bounded uniformly for all sufficiently large $t$ by
 \beqnn
 C\cdot \sup_{s\in[0,\varepsilon t]} \bigg| \frac{ \underline n\big(\zeta\in t-s+\varDelta_\delta\big)  }{\nu\big(\beta(t-s+\varDelta_\delta)\big)} -\underline n(\zeta) \bigg|,
 \eeqnn
 which also goes to $0$ as $t\to\infty$ by Theorem~1.3 in \cite{WuXu2026b}.
 Combining all estimates together, we can get the desired limit immediately. 
 \qed
 
 \begin{lemma}\label{Lemma.302}
 	For any $\delta>0$ and $\varepsilon \in(0,1/2)$, we have as $t\to\infty$,
 	\beqnn
 	\frac{J_2(t,\varepsilon,\delta)}{\nu\big(\beta(t+\varDelta_\delta)\big)}\to 0. 
 	\eeqnn
 \end{lemma}
 \proof In view of \eqref{eqn.3032}, we first write the objective fraction as   
 \beqnn
 \frac{J_2(t,\varepsilon,\delta)}{\nu\big(\beta(t+\varDelta_\delta)\big)}= \int_{\varepsilon t}^{t-\varepsilon t}
 \frac{ \overline n\big(\epsilon_s\leq x,\zeta>s\big) }{ \overline\nu(\beta s)/s
 }
 \cdot
 \frac{
 	\underline n\big(\zeta\in t-s+\varDelta_\delta\big)
 }{
 	\underline n(\zeta)
 	\nu\big(\beta(t-s+\varDelta_\delta)\big)
 }\cdot
 \frac{
 	\nu\big(\beta(t-s+\varDelta_\delta)\big)
 }{
 	\nu\big(\beta(t+\varDelta_\delta)\big)
 }
 \frac{\underline n(\zeta)\overline\nu(\beta s)}{s}\,ds,
 \eeqnn
 which obviously can be bounded by
 \beqnn 
 \sup_{r\in[\varepsilon t,t-\varepsilon t]}
 \frac{
 	\overline n\big(\epsilon_r\leq x,\zeta>r\big)
 }{
 	\overline\nu(\beta r)/r
 }
 \cdot
 \sup_{r\in[\varepsilon t,t-\varepsilon t]}
 \frac{
 	\underline n\big(\zeta\in t-r+\varDelta_\delta\big)
 }{ 
 	\nu\big(\beta(t-r+\varDelta_\delta)\big)
 }	\cdot  \sup_{r\in[\varepsilon t,t-\varepsilon t]}\frac{
 	\nu\big(\beta(t-r+\varDelta_\delta)\big)
 }{
 	\nu\big(\beta(t+\varDelta_\delta)\big)
 }
 \cdot 
 \int_{\varepsilon t}^{t-\varepsilon t}
 \frac{ \overline\nu(\beta s)}{s}\,ds.
 \eeqnn
 By Lemma~2.8 in \cite{CuiWuXu2026b} and Theorem~1.3 in \cite{WuXu2026b},  the first two suprema can be bounded uniformly for all sufficiently large $t$. 
 Similarly as in \eqref{eqn.313}, the third supremum can be bounded by $C\cdot \varepsilon^{-\alpha-1}$ uniformly for all sufficiently large $t$.  
 Finally, since $\overline{\nu} $ is integrable over $[1,\infty)$,  the last integral vanishes as $t\to\infty$. This completes the proof. 
 \qed 
 
 \begin{lemma} \label{Lemma.303}
 	For any $\delta>0$, we have 
 	\beqnn
 	\lim_{\varepsilon\to 0+}
 	\limsup_{t\to\infty} \frac{\big|J_3(t,\varepsilon,\delta) \big|}{\nu\big(\beta(t+\varDelta_\delta)\big) } =0.
 	\eeqnn
 \end{lemma}
 \proof  
 Using the change of variables to the two integrals in \eqref{eqn.3033} allows us to write $J_3(t,\varepsilon,\delta)$ as 
 \beqlb\label{eqnn4.7}
 J_3(t,\varepsilon,\delta)
 \ar=\ar
 \int_0^{\varepsilon t}
 \overline n\big(\epsilon_{t-s}\leq x,\zeta>t-s\big)\,
 \underline n(\zeta>s)\,ds
 \cr\cr
 \ar\ar-
 \int_0^{\varepsilon t+ \delta}
 \overline n\big(\epsilon_{t+\delta-s}\leq x,\zeta>t+\delta-s\big)\,
 \underline n(\zeta>s)\,ds
 \cr\cr
 \ar=\ar
 \int_0^{\varepsilon t}
 \Big[
 \overline n\big(\epsilon_{t-s}\leq x,\zeta>t-s\big)
 -
 \overline n\big(\epsilon_{t+\delta-s}\leq x,\zeta>t+\delta-s\big)
 \Big]\,
 \underline n(\zeta>s)\,ds
 \cr\cr
 \ar\ar-
 \int_{\varepsilon t}^{\varepsilon t + \delta}
 \overline n\big(\epsilon_{t+\delta-s}\leq x,\zeta>t+\delta-s\big)\,
 \underline n(\zeta>s)\,ds.
 \eeqlb
 By Lemma~2.8 in \cite{CuiWuXu2026b} again, we have as $t\to\infty$,
 \beqnn
 \eta_t:=
 \sup_{r\geq(1-\varepsilon)t}
 \bigg|
 \frac{
 	\overline n\big(\epsilon_r\leq x,\zeta>r\big)
 }{ \alpha\int_0^x\widehat V(z)\,dz
 	\cdot\overline\nu(\beta r)/(\beta r) }
 -1 \bigg| \to0.
 \eeqnn
 This firstly induces that uniformly for all
 sufficiently large $t$, 
 \beqnn
 \sup_{s\in [\varepsilon t, \varepsilon t + \delta]}\overline n\big(
 \epsilon_{t+\delta-s}\leq x,\zeta>t+\delta-s
 \big)
 \ar\leq\ar
 C\cdot   \sup_{s\in [\varepsilon t, \varepsilon t + \delta]}\frac{
 	\overline\nu\big(\beta(t+\delta-s)\big) }{ \beta(t+\delta-s) }\leq 
 C\cdot\frac{ \overline\nu\big(\beta(1-\varepsilon)t\big) }{ \beta(1-\varepsilon)t},
 \eeqnn
 Additionally, this along with the triangle inequality and the monotonicity of $r\mapsto\overline\nu(\beta r)/r$ further yields that  for some constant $C>0$ independent of $t$,
 \beqnn
 \lefteqn{ \sup_{s\in[0,\varepsilon t]}\Big| \overline n\big(\epsilon_{t-s}\leq x,\zeta>t-s\big) - \overline n\big(\epsilon_{t+\delta-s}\leq x,\zeta>t+\delta-s\big) \Big| }
 \cr\cr
 \ar\leq\ar \alpha\int_0^x\widehat V(z)\,dz \cdot \Bigg( \sup_{s\in[0,\varepsilon t]}
 \bigg[ \frac{\overline\nu\big(\beta(t-s)\big)}{\beta(t-s)} - \frac{\overline\nu\big(\beta(t+\delta-s)\big)}{\beta(t+\delta-s)} \bigg]+
 2 \cdot
 \frac{\overline\nu\big(\beta(1-\varepsilon)t\big)}
 {\beta(1-\varepsilon)t}
 \cdot\eta_t \Bigg)
 \cr\cr
 \ar\leq\ar
 C\cdot \Big[
 \frac{\overline\nu\big(\beta(1-\varepsilon)t\big)}
 {\beta(1-\varepsilon)t}
 -
 \frac{\overline\nu\big(\beta(t+\delta)\big)}
 {\beta(t+\delta)}
 \Big]
 +
 C\cdot \frac{\overline\nu\big(\beta(1-\varepsilon)t\big)}
 {\beta(1-\varepsilon)t}\cdot\eta_t.
 \eeqnn
 Here the last inequality follows from the local integrability of $\widehat{V}$. 
 Substituting these two upper bound estimates back into the right-hand side of the second equality in \eqref{eqnn4.7}, we obtain, for all sufficiently large $t$,
 \beqnn
 \frac{
 	|J_3(t,\varepsilon,\delta)|
 }{
 	\nu\big(\beta(t+\varDelta_\delta)\big)
 }
 \ar\leq\ar
 \frac{C}{
 	\nu\big(\beta(t+\varDelta_\delta)\big)
 }
 \Big[
 \frac{\overline\nu\big(\beta(1-\varepsilon)t\big)}
 {\beta(1-\varepsilon)t}
 -
 \frac{\overline\nu\big(\beta(t+\delta)\big)}
 {\beta(t+\delta)}
 \Big]
 \cr\cr
 \ar+\ar 
 C\cdot
 \frac{
 	\overline\nu\big(\beta(1-\varepsilon)t\big)
 }{
 	\beta(1-\varepsilon)t
 	\cdot\nu\big(\beta(t+\varDelta_\delta)\big)
 }
 \Big[
 \eta_t+
 \int_{\varepsilon t}^\infty
 \underline n(\zeta>s)\,ds
 \Big],
 \eeqnn
  which, by \eqref{eqn.30} and \eqref{eqn.501}, converges to $C\cdot\big((1-\varepsilon)^{-\alpha-1}-1\big)$
 	as $t\to\infty$. Letting $\varepsilon\to0+$ yields the desired limit.
 \qed

 \textbf{\textit{Proof of Theorem~\ref{MainThm01}. }}
 Combining the three limits in Lemmas~\ref{Lemma.301}-\ref{Lemma.303} together, then plugging them along with \eqref{eqn.301} back into the right-hand side of \eqref{eqn.4011}, we can  verify the asymptotic equivalence \eqref{eqn.MainThm01} for any fixed $\delta>0$. 
 We now extend it to hold uniformly for $\delta\in[\delta_0,\infty]$.
 For any large integer $n\geq1/\delta_0$ and
 $\delta\in[\delta_0,\infty]$, we have
 \beqlb\label{eqn.302}
 \sum_{k=0}^{[n\delta]-1} \mathbf P_x\big(\tau_0^-\in t+k/n+\varDelta_{1/n}\big) \ar\leq\ar \mathbf P_x\big(\tau_0^-\in t+\varDelta_\delta\big) \leq \sum_{k=0}^{[n\delta]} \mathbf P_x\big(\tau_0^-\in t+k/n+\varDelta_{1/n}\big).
 \eeqlb
 Applying the asymptotic equivalence \eqref{eqn.MainThm01} with $\delta=1/n$, we obtain a positive function $\rho_n$ on $\mathbb{R}_+$ that decreases to $0$ at infinity and satisfies, for all sufficiently large $t$ and all $k\geq0$, the probability $\mathbf P_x\big(\tau_0^-\in t+k/n+\varDelta_{1/n}\big)$ can be bounded from below and above respectively by 
 \beqnn
 \mathbf E_x[\tau_0^-] \cdot \big(1-\rho_n(t)\big)\cdot \nu\big(\beta(t+k/n+\varDelta_{1/n})\big)
 \quad \mbox{and}\quad 
 \mathbf E_x[\tau_0^-]\cdot \big(1+\rho_n(t)\big)\cdot
 \nu\big(\beta(t+k/n+\varDelta_{1/n})\big). 
 \eeqnn
 Substituting the upper bound into the second sum in \eqref{eqn.302} yields 
 \beqnn
 \frac{ \mathbf P_x\big(\tau_0^-\in t+\varDelta_\delta\big) }{ \mathbf E_x[\tau_0^-]\cdot \nu\big(\beta(t+\varDelta_\delta)\big) }
 \ar\leq\ar \frac{1+\rho_n(t)}{\nu\big(\beta(t+\varDelta_\delta)\big)} \sum_{k=0}^{[n\delta]} \nu\big(\beta(t+k/n+\varDelta_{1/n})\big) \cr
 \ar\leq\ar  \big(1+\rho_n(t)\big)\cdot \frac{\nu\big(\beta(t+\varDelta_{\delta+1/n})\big)}{\nu\big(\beta(t+\varDelta_\delta)\big) }.
 \eeqnn
 Notice that $\nu\big(\beta(t+\varDelta_{\delta+1/n})\big) =\nu\big(\beta(t+\varDelta_{\delta})\big) +\nu\big(\beta(t+\delta)+\varDelta_{\beta/n}\big) $. These along with \eqref{eqn.30} yield
 \beqnn
 \limsup_{t\to\infty} \sup_{\delta\in[\delta_0,\infty]}\frac{ \mathbf P_x\big(\tau_0^-\in t+\varDelta_\delta\big) }{ \mathbf E_x[\tau_0^-]\cdot \nu\big(\beta(t+\varDelta_\delta)\big) } 
 \leq 1+ \limsup_{t\to\infty}\sup_{\delta\in[\delta_0,\infty]}  \frac{ \nu\big(\beta(t+\delta)+\varDelta_{\beta/n}\big)}{\nu\big(\beta(t+\varDelta_\delta)\big) } \leq
 1+\frac{1}{n\delta_0}.
 \eeqnn 
 On the other hand, an analogous argument gives
 \beqnn
 \liminf_{t\to\infty}
 \inf_{\delta\in[\delta_0,\infty]}
 \frac{
 	\mathbf P_x\big(\tau_0^-\in t+\varDelta_\delta\big)
 }{
 	\mathbf E_x[\tau_0^-]\cdot
 	\nu\big(\beta(t+\varDelta_\delta)\big)
 }\geq
 \liminf_{t\to\infty}
 \inf_{\delta\in[\delta_0,\infty]}
 \frac{
 	\nu\big(\beta(t+\varDelta_{\delta-1/n})\big)
 }{
 	\nu\big(\beta(t+\varDelta_\delta)\big)
 }
 \geq
 1-\frac{1}{n\delta_0}.
 \eeqnn
 Combining these two bounds together and then
 letting $n\to\infty$ yields the desired uniform
 asymptotic equivalence  \eqref{eqn.MainThm01}.
 \qed

  \section{Proofs of Theorem~\ref{MainThm02}}
 \label{Sec.Thm02}
 \setcounter{equation}{0}
 
 This section is devoted to proving the local asymptotics for the maximum of $X$ before its first entrance into $(-\infty,0)$. By the spatial homogeneity of $X$, it is equivalent to establish that
 \beqlb\label{eqn.MainThm0201}
 \mathbf P\Big(
 \overline{X}_{\tau_{-x}^-}\in h+\varDelta_\delta
 \Big)
 \sim
 \mathbf E_x[\tau_0^-]\cdot\nu\big(h+\varDelta_\delta\big),\quad \mbox{as }h\to\infty. 
 \eeqlb
 
 Our proof relies on the It\^o reconstruction of $X$ from the excursions of $X-\underline X$, together with the key observation that, for $X$ to attain a large maximum before $\tau_{-x}^-$, the dominant scenario is that a single excursion of $X-\underline X$ produces this maximum; see Figure~\ref{Fig.01}. 
 More precisely, for any excursion $e_u\neq\Upsilon$ satisfying $\widehat H_{u-}\leq x$, since $X_{\widehat{L}_{u-}^{-1}}=\underline{X}_{\widehat{L}_{u-}^{-1}}=-\widehat H_{u-}$, we have  
 \beqnn
 e_u(s)= X_{\widehat{L}_{u-}^{-1}+s}- X_{\widehat{L}_{u-}^{-1}}= X_{\widehat{L}_{u-}^{-1}+s}+ \widehat H_{u-},\quad 0\leq s<
 \widehat L_u^{-1}-\widehat L_{u-}^{-1}. 
 \eeqnn
 Thus, in the original L\'evy process, the excursion starts from the level $-\widehat H_{u-}\in[-x,0]$, and its maximal position is $-\widehat H_{u-}+\overline e_u$, where $\overline e_u$ is the height of $e_u$.
 Consequently,
 \beqlb\label{eqnn5.14}
 \sup_{\widehat L_{u-}^{-1}\leq s<\widehat L_u^{-1}}X_s
 =
 -\widehat H_{u-}+\overline e_u,
 \quad\mbox{and}\quad
 \overline X_{\tau_{-x}^-}
 =
 0\vee
 \sup_{\substack{u\geq0:\,e_u\neq\Upsilon, \widehat H_{u-}\leq x}}
 \big(
 -\widehat H_{u-}+\overline e_u
 \big).
 \eeqlb
 Hence, after translating each excursion by the level of the running infimum at which it starts, its contribution to the running maximum can be analyzed through the corresponding excursion measure. 
 
 \begin{figure}[htbp]
 	\centering
 	
 	\begin{tikzpicture}[
 		x=1cm,
 		y=1cm,
 		font=\normalsize,
 		line cap=round,
 		line join=round,
 		axisline/.style={
 			draw=black!80,
 			line width=0.6pt,
 			-{Stealth[length=2mm,width=1.35mm]}
 		},
 		guide/.style={
 			draw=black!50,
 			line width=0.4pt,
 			dash pattern=on 0.5pt off 2pt
 		},
 		measurearrow/.style={
 			draw=black!85,
 			line width=0.55pt,
 			-{Stealth[length=1.8mm,width=1.2mm]}
 		}
 		]
 		\path[use as bounding box]
 		(-3.5,-2) rectangle (12.5,4.5);
 		
 		\draw[guide]
 		(-2.5,-0.8) -- (11.5,-0.8);
 		\draw[guide]
 		(-2.5,-1.2) -- (11.5,-1.2);
 		
 		\draw[guide]
 		(0.30,0) -- (0.30,-0.8);
 		\draw[guide]
 		(5.0,0) -- (5.0,-0.8);
 		\draw[guide]
 		(10.35,0) -- (10.35,-1.2);
 		
 		\draw[axisline]
 		(-2.5,-1.5) -- (-2.5,4);
 		\draw[axisline]
 		(-2.5,0) -- (11.5,0);
 		
 		\fill (-2.5,0) circle (1.2pt);
 		\fill (0.30,-0.8) circle (1.1pt);
 		\fill (5.0,-0.8) circle (1.1pt);
 		\fill (10.35,-1.2) circle (1.1pt);
 		
 		\node[above left=2pt]
 		at (-2.5,3.6) {\small $X_t$};
 		\node[right=4pt]
 		at (11.5,0) {\small $t$};
 		\node[left=1pt]
 		at (-2.5,0) {\small $0$};
 		\node[left=1pt]
 		at (-2.5,-0.8) {\small $-\widehat H_{\ell-}$};
 		\node[left=1pt]
 		at (-2.5,-1.2) {\small $-x$};
 		
 		\node[above=1pt]
 		at (0.30,0) {\small $\widehat L_{\ell-}^{-1}$};
 		\node[above=1pt]
 		at (5.0,0) {\small $\widehat L_{\ell}^{-1}$};
 		\node[above=1pt]
 		at (10.35,0) {\small $\tau_{-x}^-$};
 		
 		\node[fill=white,inner sep=2pt]
 		(heightlabel) at (1.5,1.2)
 		{\small $\overline e_\ell$};
 		
 		\draw[measurearrow,shorten >=2pt]
 		(1.5,-0.8) -- (heightlabel.south);
 		\draw[measurearrow,shorten >=2pt]
 		(1.5,3.2) -- (heightlabel.north);
 		
 		\draw[black,line width=0.8pt]
 		(1.425,-0.8) -- (1.575,-0.8);
 		\draw[black,line width=0.8pt]
 		(1.425,3.2) -- (1.575,3.2);
 		
 		\node[anchor=west,inner sep=0]
 		(peaklabel) at (2.82,3.65)
 		{\small $-\widehat H_{\ell-}+\overline e_\ell$};
 		\draw[measurearrow,shorten <=3pt,shorten >=4pt]
 		(peaklabel.south west) -- (1.5,3.2);
 		
 		\node[anchor=south west,inner sep=0,outer sep=0]
 		at (-2.5,-1.5)
 		{\includegraphics[width=14cm,height=5.5cm]{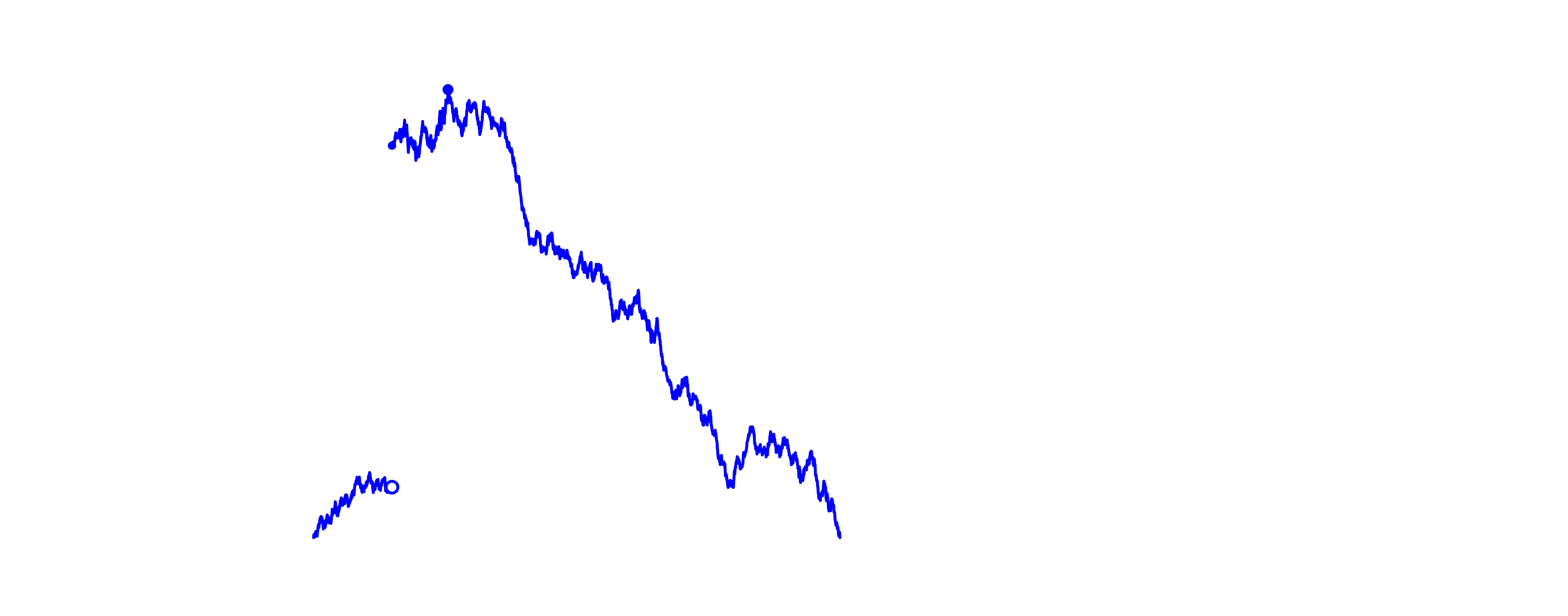}};
 		
 		\node[anchor=south west,inner sep=0,outer sep=0]
 		at (-2.5,-1.5)
 		{\includegraphics[width=14cm,height=5.5cm]{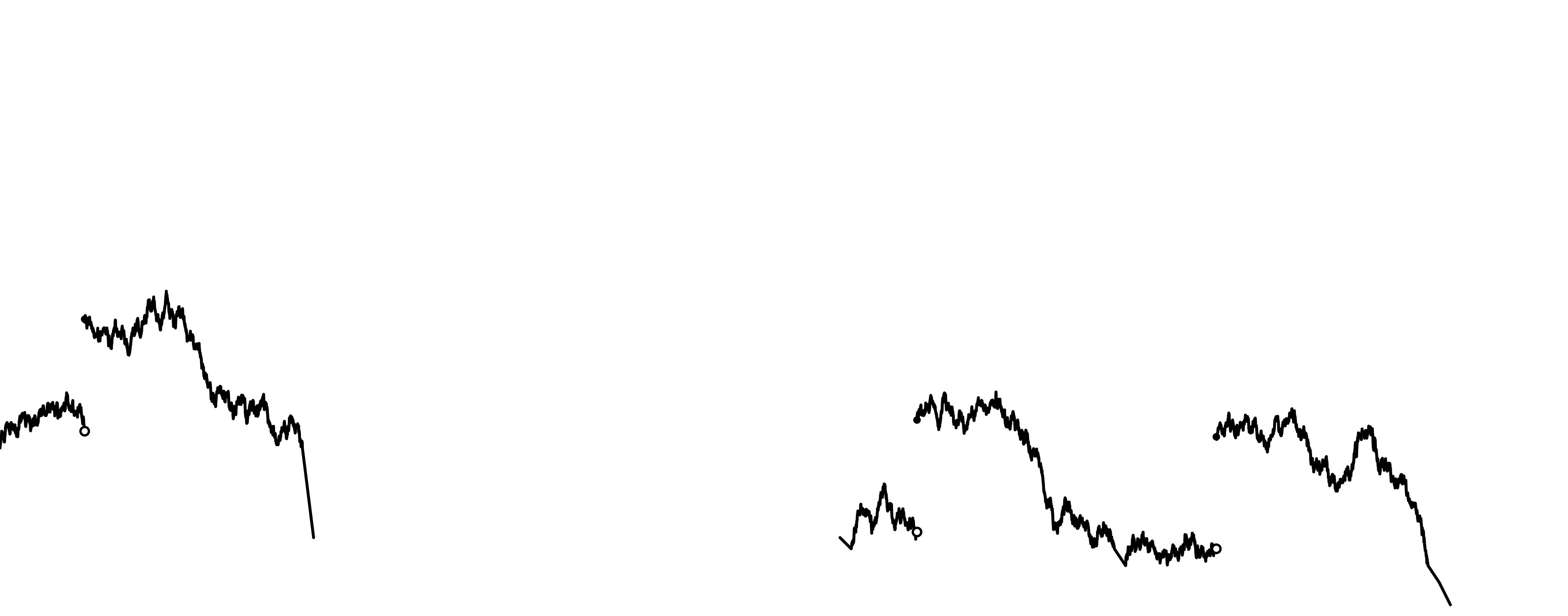}};
 	\end{tikzpicture}
 	\caption{This is a typical It\^o reconstruction for the sample path of $X$ from the excursions of $X-\underline{X}$. 
 	The maximum of $X$ before passing into $(-\infty,-x)$ is determined the heights of all excursions.}
 	\label{Fig.02}
 \end{figure}
 
 For any  interval $I\subset(0,\infty)$, we need the following set in our proof:  
 \beqlb\label{eqn5.11}
 \ar\ar
 \mathcal A_x(I)  :=  \big\{
 0\leq\ell\leq\widehat L_{\tau_{-x}^- -}: e_\ell\neq\Upsilon,\, -\widehat H_{\ell-}+\overline e_\ell\in I  \big\},
 \eeqlb
 In the proof below, we shall apply this definition with $I=h+\delta + \varDelta_\infty$ and $I=h+\varDelta_\delta$.  
 The second identity in \eqref{eqnn5.14} tells that 
 $\overline X_{\tau_{-x}^-}\in h+\varDelta_\delta$ if and only if $\#\mathcal A_x(h+\varDelta_\delta)  \geq 1$ and $\# \mathcal A_x( h+\delta + \varDelta_\infty ) =0$, which follows that
 \beqnn
 \mathbf P\Big(
 \overline{X}_{\tau_{-x}^-}\in h+\varDelta_\delta
 \Big) = \mathbf P\big(  \#\mathcal A_x(h+\varDelta_\delta)  \geq 1, \,\# \mathcal A_x\big( h+\delta + \varDelta_\infty \big) =0\big) . 
 \eeqnn
 Here and throughout, $\#A$ denotes the cardinality of the set $A$.  
 A simple calculation gives 
	\beqlb\label{eqn.405}
	\mathbf P\Big(
	\overline{X}_{\tau_{-x}^-}\in h+\varDelta_\delta
	\Big)
	\ar=\ar
	\mathbf E\big[
	\#\mathcal A_x(h+\varDelta_\delta);
	\,\#\mathcal A_x(h+\delta+\varDelta_\infty)=0
	\big]
	\cr\cr
	\ar\ar-
	\mathbf E\big[
	\#\mathcal A_x(h+\varDelta_\delta)-1;
	\,\#\mathcal A_x(h+\delta+\varDelta_\infty)=0,\,
	\#\mathcal A_x(h+\varDelta_\delta)\geq1
	\big]
	\cr\cr
	\ar=\ar
	\mathbf E\big[
	\#\mathcal A_x(h+\varDelta_\delta)
	\big]
	-
	\mathbf E\big[
	\#\mathcal A_x(h+\varDelta_\delta);
	\,\#\mathcal A_x(h+\delta+\varDelta_\infty)\geq1
	\big]
	\cr\cr
	\ar\ar-
	\mathbf E\big[
	\#\mathcal A_x(h+\varDelta_\delta)-1;
	\,\#\mathcal A_x(h+\delta+\varDelta_\infty)=0,\,
	\#\mathcal A_x(h+\varDelta_\delta)\geq1
	\big].
	\eeqlb
 Therefore, to establish the desired asymptotic equivalence~\eqref{eqn.MainThm0201}, it suffices to determine the asymptotic behavior of the three expectations appearing on the right-hand side. This is carried out in the next three lemmas. 
 
 \begin{lemma}\label{Lemma.401}
 	For any $\delta>0$, we have as $h\to\infty$, 
 	\beqnn
 	\frac{\mathbf E\big[ \#\mathcal A_x\big(h+\varDelta_\delta \big) \big] }{\nu\big(h+\varDelta_\delta\big)}
 	\to  \mathbf E_x\big[\tau_0^-\big] .
 	\eeqnn
 	
 \end{lemma}
 \proof By spatial homogeneity, \eqref{Asytau} and \eqref{eqn.DefV}, we first have 
 \beqnn
 \mathbf E_x[\tau_{0}^-] = \underline n(\zeta)\widehat V(x) = \underline n(\zeta) \mathbf E\big[ \widehat L_{\tau_{-x}^- -} \big].
 \eeqnn  
 Moreover, an application of the compensation formula \eqref{eqnn5.12} gives
 \beqlb\label{eqn.403}
 \mathbf E\big[ \#\mathcal A_x\big(h+\varDelta_\delta\big)  \big]
 = \mathbf E\bigg[ \int_0^{\widehat L_{\tau_{-x}^- -}} \underline n\Big( \overline\epsilon\in h+\widehat H_{u-} +\varDelta_\delta \Big)\,du
 \bigg].
 \eeqlb
 These together yield
 \beqnn
 \frac{\mathbf E\big[ \#\mathcal A_x\big(h+\varDelta_\delta \big) \big] }{ \nu\big(h+\varDelta_\delta\big)} -  \mathbf E_x\big[\tau_0^-\big]
 \ar=\ar  \mathbf E\bigg[ \int_0^{\widehat L_{\tau_{-x}^- -}} \bigg( \frac{\underline n\big( \overline\epsilon\in h+\widehat H_{u-} +\varDelta_\delta \big)}{\nu\big(h+\varDelta_\delta\big)} -\underline n(\zeta) \bigg)\,du \bigg].
 \eeqnn
 Since $0\leq\widehat H_{\ell-}\leq x$ for $0\leq\ell\leq\widehat L_{\tau_{-x}^- -}$,  we have 
 \beqnn
 \bigg|	\frac{\mathbf E\big[ \#\mathcal A_x\big(h+\varDelta_\delta \big) \big] }{  \nu\big(h+\varDelta_\delta\big)} - \mathbf E_x\big[\tau_0^-\big]  \bigg| 
 \leq \widehat V(x) \cdot \sup_{0\leq y\leq x} \bigg| \frac{\underline n\big( \overline\epsilon\in h+y +\varDelta_\delta \big)}{  \nu\big(h+\varDelta_\delta\big)} -\underline n(\zeta)\bigg|,
 \eeqnn
 which goes to $0$ as $h\to\infty$ by Theorem~1.4 in
 \cite{WuXu2026b} and \eqref{eqn.30}. This completes the proof.
 \qed 
 
 \begin{lemma}\label{Lemma.402}
 	For any $\delta>0$, we have as $h\to\infty$, 
 	\beqnn
 	\frac{1 }{\nu\big(h+\varDelta_\delta\big)} \cdot \mathbf E\big[ \#\mathcal A_x\big(h+\varDelta_\delta \big) ; \,\# \mathcal A_x\big( h+\delta + \varDelta_\infty \big) \geq 1 \big]
 	\to 0. 
 	\eeqnn
 \end{lemma}
 \proof  It is obvious that the objective expectation can be bounded by $ \mathbf E\big[ \#\mathcal A_x\big(h+\varDelta_\delta \big) \,\# \mathcal A_x\big( h+\delta + \varDelta_\infty \big)  \big]$. 
 Since $ \mathcal A_x\big(h+\delta + \varDelta_\infty\big) \cap \mathcal A_x\big(h+\varDelta_\delta\big) = \emptyset $,   we have
 \beqlb\label{eqnn5.17} 
 \mathbf E\big[ \#\mathcal A_x\big(h+\varDelta_\delta \big) \,\# \mathcal A_x\big( h+\delta + \varDelta_\infty \big)  \big] 
 \ar=\ar
 \mathbf E\bigg[
 \sum_{0\leq u<t\leq\widehat L_{\tau_{-x}^- -}}
 \mathbf 1_{\{ u\in\mathcal A_x(h+\delta + \varDelta_\infty) \}}
 \mathbf 1_{\{ t\in\mathcal A_x(h+\varDelta_\delta) \}} \bigg]
 \cr\cr
 \ar\ar+
 \mathbf E\bigg[ \sum_{0\leq u<t\leq\widehat L_{\tau_{-x}^- -}}
 \mathbf 1_{\{ u\in\mathcal A_x(h+\varDelta_\delta) \}}
 \mathbf 1_{\{ t\in\mathcal A_x(h+\delta + \varDelta_\infty) \}} \bigg].
 \eeqlb
 Since the right-end point of every successive
 excursion in $\mathcal A_x((h,\infty))$ is an $(\mathscr{F}_t)$-stopping time,
 by the strong Markov property of $X$ and the
 spatial homogeneity we have 
 \beqlb\label{eqn.404}
 \lefteqn{\mathbf E\bigg[ \sum_{0\leq u<t\leq\widehat L_{\tau_{-x}^- -}} \mathbf 1_{\{ u\in\mathcal A_x(h+\delta + \varDelta_\infty) \}} \mathbf 1_{\{ t\in\mathcal A_x(h+\varDelta_\delta) \}} \bigg]}\cr\cr
 \ar=\ar \mathbf E\bigg[ \sum_{0\leq u<\widehat L_{\tau_{-x}^- -}} \mathbf 1_{\{ u\in\mathcal A_x (h+\delta + \varDelta_\infty)\}} \cdot \mathbf E_{-\widehat H_u}\Big[ \#\mathcal A_x\big( h+\varDelta_\delta \big) \Big]
 \bigg]\cr\cr
 \ar\leq\ar \sup_{0\leq z\leq x} \mathbf E\Big[ \#\mathcal A_{x-z}\big(
 h+z+\varDelta_\delta \big) \Big]\cdot \mathbf E\big[
 \#\mathcal A_x\big(h+\delta+ \varDelta_\infty\big) \big] .
 \eeqlb
 Here the last inequality follows from the fact that $0\leq\widehat H_u\leq x$ when  $u<\widehat L_{\tau_{-x}^- -}$. 
 By Lemma~\ref{Lemma.401} and \eqref{eqn.30}, we have uniformly for all sufficiently large $h$,
 \beqnn
 \mathbf E\bigg[ \sum_{0\leq u<t\leq\widehat L_{\tau_{-x}^- -}} \mathbf 1_{\{ u\in\mathcal A_x(h+\delta + \varDelta_\infty) \}} \mathbf 1_{\{ t\in\mathcal A_x(h+\varDelta_\delta) \}} \bigg]  
 \ar\leq\ar  C\cdot \widehat V(x)\cdot \nu\big(h+\varDelta_\delta\big)\cdot \mathbf E\big[  \#\mathcal A_x\big(h+\delta + \varDelta_\infty\big) \big]. 
 \eeqnn 
 We now consider the second expectation on the right-hand side of \eqref{eqnn5.17}.  
 By \eqref{eqn.403}, we have 
 \beqnn
 \mathbf E\big[ \#\mathcal A_x\big(h+\delta+\varDelta_\infty\big) \big] = \mathbf E\bigg[ \int_0^{\widehat L_{\tau_{-x}^- -}} \underline n\Big( \overline\epsilon>h+\widehat H_{u-}+\delta \Big)\,du \bigg]
 \eeqnn
 Similarly to our preceding arguments, this along with \eqref{eqn.DefV} and the
 monotonicity of $\widehat V$ yields that uniformly for all sufficiently large $h$, 
 \beqnn
 \mathbf E\bigg[ \sum_{0\leq u<t\leq\widehat L_{\tau_{-x}^- -}} \mathbf 1_{\{ u\in\mathcal A_x(h+\varDelta_\delta) \}} \mathbf 1_{\{ t\in\mathcal A_x(h+\delta + \varDelta_\infty) \}} \bigg]
 \leq C\cdot  
 \widehat V(x) \cdot
 \underline n\big(\overline\epsilon>h\big)\cdot 
 \mathbf E\big[
 \#\mathcal A_x\big(h+\varDelta_\delta \big)
 \big].
 \eeqnn 
 Substituting these two upper bound estimates back into  \eqref{eqnn5.17}, we obtain that uniformly for all sufficiently large $h$, 
 \beqnn
 \lefteqn{  \frac{1 }{\nu\big(h+\varDelta_\delta\big)} \cdot \mathbf E\big[ \#\mathcal A_x\big(h+\varDelta_\delta \big) ; \,\# \mathcal A_x\big( h+\delta + \varDelta_\infty \big) \geq 1 \big] }\cr
 \ar\leq\ar 
 \frac{1 }{\nu\big(h+\varDelta_\delta\big)} \cdot \mathbf E\big[ \#\mathcal A_x\big(h+\varDelta_\delta \big) \,\# \mathcal A_x\big( h+\delta + \varDelta_\infty \big)  \big] \cr
 \ar \leq\ar  C\cdot  
 \widehat V(x) \cdot \bigg(\mathbf E\big[  \#\mathcal A_x\big(h+\delta + \varDelta_\infty\big) \big] + \frac{\mathbf E\big[ \#\mathcal A_x\big(h+\varDelta_\delta \big) \big]}{\nu\big(h+\varDelta_\delta\big)} \cdot  \underline n\big(\overline\epsilon>h\big) \bigg),
 \eeqnn
 which goes to $0$ as $h\to\infty$ by Lemma~\ref{Lemma.401}.  This completes the proof. 
 \qed 
 
 \begin{lemma}\label{Lemma.403}
 	For any $\delta>0$, we have as $h\to\infty$, 
 	\beqnn
 	\frac{1 }{\nu\big(h+\varDelta_\delta\big)} \cdot \mathbf E\big[ \#\mathcal A_x\big(h+\varDelta_\delta \big) -1, \,\# \mathcal A_x\big( h+\delta + \varDelta_\infty \big) =0,\, \#\mathcal A_x\big(h+\varDelta_\delta \big) \geq 1\big]
 	\to 0. 
 	\eeqnn
 \end{lemma}
 \proof Obviously, the objective expectation can be bounded by   
 \beqlb\label{eqn.406}
 \mathbf E\big[ \#\mathcal A_x\big(h+\varDelta_\delta \big) \,\big( \#\mathcal A_x\big(h+\varDelta_\delta \big) -1 \big) \big]
 \ar=\ar 2\cdot \mathbf E\bigg[
 \sum_{0\leq u<t\leq\widehat L_{\tau_{-x}^- -}}
 \mathbf 1_{\{ u\in\mathcal A_x(h+ \varDelta_\delta) \}}
 \mathbf 1_{\{ t\in\mathcal A_x(h+\varDelta_\delta) \}} \bigg].
 \eeqlb
 By the same strong Markov argument as in \eqref{eqn.404}, Lemma~\ref{Lemma.401} and  \eqref{eqn.30}, the last expression can be bounded uniformly for all sufficiently large $h$ by
 \beqlb\label{eqn.407}
 C\cdot \nu \big(h+\varDelta_\delta \big) \cdot \widehat V(x)\cdot 
 \mathbf E\big[
 \#\mathcal A_x\big(h+\varDelta_\delta \big) \big] \leq
 C\cdot \big|\widehat V(x)\cdot \nu\big(h+\varDelta_\delta \big)\big|^2.
 \eeqlb
 The desired limit follows directly and the proof ends. 
 \qed 
 
 \medskip
 \textbf{\textit{Proof of Theorem~\ref{MainThm02}. }} 
 Substituting the three limits established in Lemmas~\ref{Lemma.401}-\ref{Lemma.403} into \eqref{eqn.405} immediately yields the asymptotic equivalence \eqref{eqn.MainThm0201}  for any fixed $\delta>0$. 
 The extension to uniformity in $\delta\in[\delta_0,\infty]$
 follows by the same discretization argument as in the proof
 of Theorem~\ref{MainThm01}.
 \qed

 \section{Proofs of Theorem~\ref{MainThm03}}
 \label{Sec.Thm03}
 \setcounter{equation}{0}
 
 In this section, we prove Theorem~\ref{MainThm03}, which states that, as $t\to\infty$, the simultaneous occurrence of the two events $ \big\{\overline X_{\tau_0^-} \in\beta(t+\varDelta_{\delta_1})\big\}$ and $\big\{\tau_0^-\in t+\varDelta_{\delta_2} \big\}$ becomes asymptotically negligible under $\mathbf P_x$.
 By the spatial homogeneity, it is equivalent to prove that as $t\to\infty$, 
 \beqlb \label{eqn.5011}
 \lefteqn{ \mathbf P\Big(x+ \overline{X}_{\tau_{-x}^{-}}\in\beta(t+\varDelta_{\delta_1})
 	\,\Big|\,
 	\tau_{-x}^-\in t+\varDelta_{\delta_2}
 	\Big)}\cr\cr
 \ar\sim\ar
 \frac{\delta_1}{\delta_2}\cdot
 \mathbf P\Big(
 \tau_{-x}^-\in t+\varDelta_{\delta_2}
 \,\Big|\, x+\overline{X}_{\tau_{-x}^{-}}\in\beta(t+\varDelta_{\delta_1})
 \Big)
 \to0. 
 \eeqlb 
 The proof follows a strategy similar to that used for
 Theorem~\ref{MainThm02}, but requires an additional decomposition
 according to whether the combined length of the following two
 time intervals exceeds $R$: the time from $0$ to the start of
 the excursion with a large maximum and the time from its end
 to $\tau_{-x}^-$.
 
 Recall the set $\mathcal A_x(I)$ defined in \eqref{eqn5.11}.
 The main idea is as follows.
 For all sufficiently large $t$, on the event $\{x+\overline X_{\tau_{-x}^-}\in\beta(t+\varDelta_{\delta_1})\}$,
 we have $\mathcal A_x\big(\beta(t+\varDelta_{\delta_1})-x\big)\neq\emptyset$.
 For fixed $R>x/\beta$, we first consider excursions indexed by this set for which the combined length of the two time intervals is at most $R$.
 Their maxima are close to $\beta t$.
 On the second event, their lifetimes are also close to $t$.
 We use the joint local estimate under $\underline n$ to show that their contribution is asymptotically negligible as $t\to\infty$. For the remaining excursions indexed by $\mathcal A_x\big(\beta(t+\varDelta_{\delta_1})-x\big)$, at least one of the two time intervals has length greater than $R/2$.
 We control their contribution using the compensation formula and the strong Markov property.
 To this end, for any constants $x,R>0$, define
 \beqnn 
 \mathcal A_{x}^{\leq R}(I)
 \ar:=\ar \mathcal A_x(I) \cap \big\{ 0\leq s\leq\widehat L_{\tau_{-x}^- -}: \widehat L^{-1}_{s-} +\tau_{-x}^- -\widehat L^{-1}_s \leq R \big\}
 \quad \mbox{and}\quad
 \mathcal A_{x}^{> R}(I) := \mathcal A_x(I)-\mathcal A_x^{\leq R}(I).
 \eeqnn

 \medskip
 \textbf{\textit{Proof of Theorem~\ref{MainThm03}. }} 
 By the multiplication rule for conditional probabilities,
 \beqnn
 \lefteqn{ \mathbf P\Big( x+
 	\overline{X}_{\tau_{-x}^{-}}\in\beta(t+\varDelta_{\delta_1})
 	\,\Big|\,
 	\tau_{-x}^-\in t+\varDelta_{\delta_2}
 	\Big)}\cr\c\r
 \ar=\ar \frac{
 	\mathbf P\big(x+
 	\overline X_{\tau_{-x}^-}
 	\in\beta(t+\varDelta_{\delta_1})
 	\big)
 }{
 	\mathbf P\big(
 	\tau_{-x}^-\in t+\varDelta_{\delta_2}
 	\big)
 }\cdot
 \mathbf P\Big(
 \tau_{-x}^-\in t+\varDelta_{\delta_2}
 \,\Big|\,
 x+\overline{X}_{\tau_{-x}^{-}}\in\beta(t+\varDelta_{\delta_1})
 \Big).
 \eeqnn 
 By Theorems~\ref{MainThm01} and~\ref{MainThm02}, the ratio on the right-hand side converges to $\delta_1/\delta_2>0$ and hence the asymptotic equivalence in \eqref{eqn.5011} holds. 
 It remains to prove the limit holds.  
 For all sufficiently large $t>x/\beta $,  if $ x+\overline X_{\tau_{-x}^-} \in\beta(t+\varDelta_{\delta_1})$,  we have $  \mathcal A_x\big( \beta t-x+\varDelta_{\beta\delta_1} \big) \neq \emptyset$, which follows that 
 \beqnn
 \mathbf P\big(x+\overline X_{\tau_{-x}^-} \in\beta(t+\varDelta_{\delta_1}),\, \tau_{-x}^-\in t+\varDelta_{\delta_2} \big) 
 \ar\leq\ar \mathbf P\big( \# \mathcal A_x\big( \beta t-x+\varDelta_{\beta\delta_1} \big)\geq 1,\, \tau_{-x}^-\in t+\varDelta_{\delta_2} \big) \cr\cr
 \ar\leq\ar\mathbf{E}\big[ \# \mathcal A_x\big( \beta t-x+\varDelta_{\beta\delta_1} \big),\, \tau_{-x}^-\in t+\varDelta_{\delta_2} \big],
 \eeqnn
 which, for any $R>0$, can be decomposed into the following two terms
 \beqlb\label{eqn.I2}
 I_1(t,R) \ar:=\ar \mathbf{E}\big[ \# \mathcal A_x^{\leq R}\big( \beta t-x+\varDelta_{\beta\delta_1} \big),\, \tau_{-x}^-\in t+\varDelta_{\delta_2} \big],\cr\cr
 I_2(t,R) \ar:=\ar \mathbf{E}\big[ \# \mathcal A_x^{>R}\big( \beta t-x+\varDelta_{\beta\delta_1} \big),\, \tau_{-x}^-\in t+\varDelta_{\delta_2} \big].
 \eeqlb
 Consequently, it suffices to prove that as $t\to\infty$ and then $R\to\infty$, 
 \beqlb\label{eqn.5012}
 \frac{I_i(t,R)}{\mathbf P\big( \tau_{-x}^-\in t+\varDelta_{\delta_2} \big)} \to 0,\quad i=1,2. 
 \eeqlb 
 
 \textbf{Case $i=1$.} For any fixed $R>x/\beta$, on the event $\{\tau_{-x}^-\in t+\varDelta_{\delta_2}\}$, one can easily identify that for any $u\in\mathcal A_{x}^{\leq R}\big(
 \beta t-x+\varDelta_{\beta\delta_1}\big)$, 
 \beqnn
 \overline e_u\in(\beta t-x,\beta t+\beta\delta_1] \subset\beta(t-R+\varDelta_{R+\delta_1})
 \quad\mbox{and}\quad
 \zeta(e_u) \in(t-R,t+\delta_2] = t-R+\varDelta_{R+\delta_2},
 \eeqnn
 Consequently, we have 
 \beqnn
 I_1(t,R)  \leq \mathbf{E}\bigg[ \sum_{u\in[0, \widehat L_{\tau_{-x}^- -]}} \mathbf{1}_{ \{\overline e_u\in \beta(t-R+\varDelta_{R+\delta_1})\} }\cdot \mathbf{1}_{\{ \zeta(e_u) \in t-R+\varDelta_{R+\delta_2}\}}\bigg]
 \eeqnn
 Applying the compensation formula \eqref{eqnn5.12} and then \eqref{eqn.DefV} to the right-hand side  gives
 \beqnn
 \frac{I_1(t,R)}{\mathbf P\big( \tau_{-x}^-\in t+\varDelta_{\delta_2} \big)} 
 \ar\leq\ar  \frac{\mathbf E\big[ \widehat L_{\tau_{-x}^- -} \big]}{\mathbf P\big( \tau_{-x}^-\in t+\varDelta_{\delta_2} \big)}  \cdot \underline n\big( \bar\epsilon \in\beta(t-R+\varDelta_{R+\delta_1}),\, \zeta \in t-R+\varDelta_{R+\delta_2} \big)\cr\cr
 \ar\leq\ar \widehat{V}(x)  \cdot \frac{\underline{n}\big( \zeta \in t-R+\varDelta_{R+\delta_2} \big)}{\mathbf P\big( \tau_{-x}^-\in t+\varDelta_{\delta_2} \big)}  \cdot \underline n\big( \bar\epsilon \in\beta(t-R+\varDelta_{R+\delta_1})\,\big|\, \zeta \in t-R+\varDelta_{R+\delta_2} \big),
 \eeqnn
 By Theorem~1.3 in \cite{WuXu2026b}, Theorem~\ref{MainThm01} and \eqref{eqn.30}, the last ratio is bounded for all sufficiently large $t$. 
 This along with Theorem~1.5 in \cite{WuXu2026b} proves that the limit \eqref{eqn.5012} with $i=1$ holds as $t\to\infty$ for any $R>x/\beta$. 
 
 \textbf{Case $i=2$.} It is obvious that $I_2(t,R) \leq \mathbf{E}\big[ \# \mathcal A_x^{>R}\big( \beta t-x+\varDelta_{\beta\delta_1} \big)\big]$.  
 For any $u\in  \mathcal A_x^{>R}\big( \beta t-x+\varDelta_{\beta\delta_1} \big)$, we have the following properties: (1) $\widehat L^{-1}_{u-}>R/2$ or $ \tau_{-x}^--\widehat L^{-1}_u>R/2$; (2) $\widehat{H}_{u-} \in[0,x]$; (3) $\overline{e}_u \in \widehat{H}_{u-}+ \beta t-x+\varDelta_{\beta\delta_1} \subset \beta t-x+\varDelta_{x+\beta\delta_1} $. 
 This yields 
 \beqlb\label{eqn.MainThm03.I2split}
 \mathbf{E}\big[ \# \mathcal A_x^{>R}\big( \beta t-x+\varDelta_{\beta\delta_1} \big)\big]\ar\leq\ar   \mathbf{E}\bigg[ \sum_{0\leq u\leq\widehat L_{\tau_{-x}^- -}}
 \mathbf 1_{\{ u\in\mathcal A_x( \beta t-x+\varDelta_{\beta\delta_1}) \}} \cdot
 \mathbf 1_{\{\widehat L^{-1}_{u-}>R/2\}} \bigg]\cr
 \ar\ar + \mathbf{E}\bigg[ \sum_{0\leq u\leq\widehat L_{\tau_{-x}^- -}}
 \mathbf 1_{\{ u\in\mathcal A_x( \beta t-x+\varDelta_{\beta\delta_1}) \}} \cdot
 \mathbf 1_{\{\tau_{-x}^--\widehat L^{-1}_u>R/2\}} \bigg]. 
 \eeqlb 
 Applying the compensation formula \eqref{eqnn5.12} to the first expectation on the right-hand side and then using property (3) gives
 \beqnn
 \lefteqn{\mathbf{E}\bigg[ \sum_{0\leq u\leq\widehat L_{\tau_{-x}^- -}}
 	\mathbf 1_{\{ u\in\mathcal A_x( \beta t-x+\varDelta_{\beta\delta_1}) \}} \cdot
 	\mathbf 1_{\{\widehat L^{-1}_{u-}>R/2\}} \bigg]}\cr\cr
 \ar\leq\ar \underline n\big( \bar\epsilon\in \beta t-x+\varDelta_{x+\beta\delta_1}\big)\cdot \mathbf{E}\bigg[ \int_0^{\widehat L_{\tau_{-x}^- -}} \mathbf 1_{\{\widehat L^{-1}_{s-}>R/2\}}\,ds
 \bigg].
 \eeqnn 
 For the second expectation, by using the strong Markov property of $X$ at the  right-end time $d_u=\widehat L^{-1}_u$, then the compensation formula \eqref{eqnn5.12} and finally the property (2), 
 \beqnn
 \lefteqn{\mathbf{E}\bigg[ \sum_{0\leq u\leq\widehat L_{\tau_{-x}^- -}}
 	\mathbf 1_{\{ u\in\mathcal A_x( \beta t-x+\varDelta_{\beta\delta_1}) \}} \cdot
 	\mathbf 1_{\{\tau_{-x}^--\widehat L^{-1}_u>R/2\}} \bigg] }\cr\cr
 \ar=\ar
 \mathbf E\Big[\sum_{0\leq u<\widehat L_{\tau_{-x}^- -}} \mathbf 1_{\{ u\in\mathcal A_x( \beta t-x+\varDelta_{\beta\delta_1}) \}} \cdot \mathbf P_{-\widehat H_u}\big( \tau_{-x}^->R/2 \big) \Big]\cr\cr
 \ar\leq \ar \widehat{V}(x)\cdot  \underline n\big( \bar\epsilon\in \beta t-x+\varDelta_{x+\beta\delta_1}\big)\cdot \mathbf P\big( \tau_{-x}^->R/2 \big).
 \eeqnn
 Combining these two estimates with \eqref{eqn.MainThm03.I2split} gives
 \beqnn
 \frac{ \mathbf{E}\big[ \# \mathcal A_x^{>R}\big( \beta t-x+\varDelta_{\beta\delta_1} \big)\big]}{\mathbf P\big( \tau_{-x}^-\in t+\varDelta_{\delta_2} \big)} 
 \ar\leq\ar
 C\cdot\frac{\underline n\big(\bar\epsilon\in \beta t-x+\varDelta_{x+\beta\delta_1}\big)}{\mathbf P\big( \tau_{-x}^-\in t+\varDelta_{\delta_2} \big)}\cdot \Bigg( \mathbf E\Big[\int_0^{\widehat L_{\tau_{-x}^- -}} \mathbf 1_{\{\widehat L^{-1}_{s-}>R/2\}}\,ds \Big] + \mathbf P\big( \tau_{-x}^->R/2 \big)\Bigg).
 \eeqnn
 By Theorem~1.4 in \cite{WuXu2026b},  Theorem~\ref{MainThm01} and \eqref{eqn.30}, the ratio is bounded for large $t$.
 The expectation and the probability both go to $0$ as $R\to\infty$
 by \eqref{eqn.DefV} and the finiteness of $\tau_{-x}^-$. 
 Hence, 
	\beqlb\label{eqn5.5}
	\lim_{R\to\infty}\limsup_{t\to\infty}
	\frac{
	I_{2}(t,R)
	}{
	\mathbf P\big(\tau_{-x}^-\in t+\varDelta_{\delta_2}\big)
	}\leq \lim_{R\to\infty}\limsup_{t\to\infty}
	\frac{
	\mathbf E\big[
	\#\mathcal A_x^{>R}
	\big(\beta t-x+\varDelta_{\beta\delta_1}\big)
	\big]
	}{
	\mathbf P\big(\tau_{-x}^-\in t+\varDelta_{\delta_2}\big)
	}
	=0.
	\eeqlb
	This proves the limit \eqref{eqn.5012} with $i=2$. The proof ends. 
 \qed

 \section{Proofs of Theorem \ref{MainThm04}}
 \label{Sec.Thm04}
 \setcounter{equation}{0}
 
 This section is devoted to the proof of Theorem~\ref{MainThm04}, which provides second-order localization estimates for the maximum $\overline X_{\tau_0^-}$ and the first-passage time $\tau_0^-$ under $\mathbf P_x$.
 Similarly to before, by the spatial homogeneity, we prove the corresponding result under $\mathbf P=\mathbf P_0$, with $\overline X_{\tau_0^-}$ and $\tau_0^-$ replaced by $x+\overline X_{\tau_{-x}^-}$ and $\tau_{-x}^-$, respectively.
 
 To prove the first limit in \eqref{eqn.MainThm04}, we recall from \cite{CuiWuXu2025} that $X$ typically needs an early large jump to remain in $[-x,\infty)$ for a long time.
 Conditional on $\tau_{-x}^-\in t+\varDelta_\delta$, we show that the excursion containing this jump typically occupies almost the entire interval $[0,\tau_{-x}^-]$.
 We then transfer height localization from $\underline n$ to $\mathbf P$ and show that multiple high excursions make a negligible contribution.
 For the second limit in \eqref{eqn.MainThm04}, we use the excursion decomposition from the proof of Theorem~\ref{MainThm02}.
 A large maximum is typically produced by a single excursion of $X-\underline X$, shifted by its starting position in $X$.
 We apply lifetime localization under $\underline n$ and use the estimates from the preceding section to control the time before and after this excursion.

 \begin{lemma}\label{Lemma.MainThm04.transfer}
 	For any $x,\delta>0$ and any positive function $\theta$ on $\mathbb R_+$ satisfying $|X_t+\beta t|/\theta_t\overset{\rm p}\to 0$, we have
 	\beqnn
 	\lim_{R\to\infty} \limsup_{t\to\infty} \mathbf P\Big( \#\mathcal A_{x}^{\leq R}\big(\beta t-x-\theta_t/2+\varDelta_{\theta_t} \big)=0 \,\Big|\, \tau_{-x}^-\in t+\varDelta_\delta \Big) = 0.
 	\eeqnn
 \end{lemma}

 \textbf{\textit{Proof of Theorem~\ref{MainThm04}.}} 
 Firstly, we prove the first limit in \eqref{eqn.MainThm04}. 
 The law of large numbers gives $|X_t+\beta t|/(\beta t )\overset{\rm p}\to0$, we may assume that $\theta_t\leq\beta t$.
 By spatial homogeneity,
 \beqnn
 \mathbf P_x\Big( \big|\overline X_{\tau_0^-}-\beta t\big|<\theta_t \,\Big|\, \tau_0^-\in t+\varDelta_\delta \Big)
 = \mathbf P\Big( \big|x+\overline X_{\tau_{-x}^-}-\beta t\big|<\theta_t \,\Big|\, \tau_{-x}^-\in t+\varDelta_\delta \Big).
 \eeqnn
 Since $\theta_t\leq\beta t$ for all sufficiently large $t$, we have $\beta t-x-\theta_t/2\geq\beta t/2-x>0$ and then 
 \beqnn
 \mathcal A_{x}^{\leq R}\big(\beta t-x-\theta_t/2+\varDelta_{\theta_t} \big)
 \subset \mathcal A_x\big( \beta t/2-x+\varDelta_\infty \big).
 \eeqnn
 Combining this together with \eqref{eqnn5.14} yields
 \beqnn
 \big\{ \#\mathcal A_{x}^{\leq R}\big(\beta t-x-\theta_t/2+\varDelta_{\theta_t} \big)\geq1,\, \#\mathcal A_x\big( \beta t/2-x+ \varDelta_\infty \big)=1 \big\} \subset \big\{ \big| x+\overline X_{\tau_{-x}^-}-\beta t \big|<\theta_t \big\}.
 \eeqnn
 Consequently,
 \beqnn
 \mathbf P_x\Big( \big| \overline X_{\tau_0^-}-\beta t \big|<\theta_t \,\Big|\, \tau_0^-\in t+\varDelta_\delta \Big)
 \ar\geq\ar 1- \mathbf P\Big( \#\mathcal A_{x}^{\leq R}\big(\beta t-x-\theta_t/2+\varDelta_{\theta_t} \big)=0 \,\Big|\, \tau_{-x}^-\in t+\varDelta_\delta \Big)\cr\cr
 \ar\ar -  \mathbf P\Big( \#\mathcal A_x\big( \beta t/2-x+ \varDelta_\infty\big)\geq2 \,\Big|\, \tau_{-x}^-\in t+\varDelta_\delta \Big).
 \eeqnn
 In view of Lemma~\ref{Lemma.MainThm04.transfer}, the first probability on the right-hand side goes to $0$ as $t\to\infty$ and then $R\to\infty$. 
 Hence it suffices to prove that the second probability also tends to $0$ as $t\to\infty$. 
 Indeed, an application of \eqref{eqn.406}-\eqref{eqn.407} with $h=\beta t/2-x$ and $\delta=\infty$ gives 
 \beqnn
 \mathbf E\Big[ \#\mathcal A_x\big( (\beta t/2-x,\infty) \big) \Big( \#\mathcal A_x\big( \beta t/2-x +\varDelta_\infty \big)-1 \Big) \Big]\leq C\cdot \widehat V(x)^2\cdot  \overline\nu(\beta t/2-x)^2.
 \eeqnn
 Therefore, by Theorem~\ref{MainThm01} and \eqref{eqn.30}, we have uniformly for all large $t$, 
 \beqnn
 \lefteqn{ \mathbf P\Big( \#\mathcal A_x\big( \beta t/2-x + \varDelta_\infty\big)\geq2 \,\Big|\, \tau_{-x}^-\in t+\varDelta_\delta \Big)} \cr\cr
 \ar\leq\ar \frac{\mathbf E\big[ \#\mathcal A_x\big(\beta t/2-x+ \varDelta_\infty\big)\,\big( \#\mathcal A_x\big(\beta t/2-x+ \varDelta_\infty\big)-1\big)\big]}{2\cdot \mathbf P\big( \tau_{-x}^-\in t+\varDelta_\delta\big)}
 \leq C\cdot t\cdot \overline\nu(\beta t), 
 \eeqnn
 which vanishes as $t\to\infty$.
 These prove the first limit in \eqref{eqn.MainThm04}. 
 
 We now start to prove the second limit in \eqref{eqn.MainThm04}.
 Fix any large $R>0$. By \eqref{eqnn5.14}, for all sufficiently large $t$,
 conditioned on the event $ \{x+\overline X_{\tau_{-x}^-}\in\beta(t+\varDelta_\delta)\}$ we have  $ \#\mathcal A_x\big(
 (\beta t-x,\beta t-x+\beta\delta]
 \big)\geq 1.
 $
 Hence, by the spatial homogeneity,
 \beqnn
 \lefteqn{\mathbf P_x\Big( |\tau_0^--t|\leq\theta_t \,\Big|\, \overline X_{\tau_0^-}\in \beta(t+\varDelta_\delta) \Big)}\cr\cr
 \ar=\ar 1-\frac{\mathbf P\big(|\tau_{-x}^--t|>\theta_t,\, x+\overline X_{\tau_{-x}^-}  \in \beta(t+\varDelta_\delta) \big)}{\mathbf P_x\big( \overline X_{\tau_0^-} \in \beta(t+\varDelta_\delta) \big) }  
 \geq 1-\frac{J_1(t,R)+ J_2(t,R)}{\mathbf P_x\big( \overline X_{\tau_0^-} \in \beta(t+\varDelta_\delta) \big)}
 \eeqnn
 with $J_1(t,R) := \mathbf E\big[ \#\mathcal A_{x}^{\leq R}\big( \beta t-x+\varDelta_{\beta\delta}\big); \,|\tau_{-x}^--t|>\theta_t \big] $ and $J_2(t,R) :=  \mathbf E\big[ \#\mathcal A_{x}^{> R}\big( \beta t-x+\varDelta_{\beta\delta}\big) \big] $. 
 Consequently, it suffices to prove that as $t\to\infty$ and then $R\to\infty$,  
 \beqlb\label{eqn.601}
 \frac{J_i(t,R)}{\mathbf P_x\big( \overline X_{\tau_0^-} \in \beta(t+\varDelta_\delta) \big)} \to 0,\quad i=1,2. 
 \eeqlb
 For all sufficiently large $t$, on the event $\{|\tau_{-x}^--t|>\theta_t\}$, we have for any $u\in\mathcal A_{x}^{\leq R}\big(\beta t-x+\varDelta_{\beta\delta}\big)$,  
 \beqnn
 \bar e_u \in \beta\big(t-x/\beta+\varDelta_{\delta+x/\beta} \big)
 \quad\mbox{and}\quad
 \big|\zeta(e_u)-(t-x/\beta) \big|\geq|\tau_{-x}^--t|-R-x/\beta >\theta_{t}/2.
 \eeqnn
 This along with the compensation formula \eqref{eqnn5.12} and \eqref{eqn.DefV} induces that 
 \beqnn
 J_1(t,R) 
 \ar\leq\ar \mathbf E\bigg[  \sum_{0\leq u\leq\widehat L_{\tau_{-x}^- -}}
 \mathbf 1_{\{ \bar e_u \in \beta( t-x/\beta +\varDelta_{\delta+x/\beta} ) \}} \cdot \mathbf 1_{\{|\zeta(e_u)-(t-x/\beta)| >\theta_t/2\big\}} \bigg] \cr
 \ar=\ar  \widehat V(x)\cdot \underline n\Big(
 \bar\epsilon \in \beta\big( t-x/\beta +\varDelta_{\delta+x/\beta} \big),\,
 \big| \zeta-(t-x/\beta) \big|
 > \theta_t/2 \Big). 
 \eeqnn
 Consequently, we have uniformly for all sufficiently large $t$, 
 \beqnn
 \lefteqn{\frac{J_1(t,R)}{\mathbf P_x\big( \overline X_{\tau_0^-} \in \beta(t+\varDelta_\delta) \big)}}\cr\cr 
 \ar\leq\ar   \widehat V(x)\cdot \frac{
 	\underline n\Big(
 	\bar\epsilon
 	\in
 	\beta\big(
 	t-x/\beta
 	+\varDelta_{\delta+x/\beta}
 	\big)
 	\Big)
 }{
 	\mathbf P_x\big(
 	\overline X_{\tau_0^-}
 	\in
 	\beta(t+\varDelta_\delta)
 	\big)
 }\cdot \underline n\Big(
 \big|
 \zeta-(t-x/\beta)
 \big|
 >
 \frac{\theta_t}{2}
 \,\Big|\,
 \bar\epsilon
 \in
 \beta\big(
 t-x/\beta
 +\varDelta_{\delta+x/\beta}
 \big)
 \Big).
 \eeqnn 
 In view of Theorem~1.4 in \cite{WuXu2026b}, Theorem~\ref{MainThm02} and \eqref{eqn.30}, the ratio on the right-hand side can be bounded uniformly for all large $t$. 
 This along with Theorem~1.6 in \cite{WuXu2026b} yields that  the limit \eqref{eqn.601} with $i=1$ holds as $t\to\infty$ for all large $R$. 
 
  For $J_2(t,R)$, we first write it as 
 	\beqnn
 	\frac{J_2(t,R)}
 	{\mathbf P_x\big(
 		\overline X_{\tau_0^-}\in\beta(t+\varDelta_\delta)
 		\big)}
 	=
 	\frac{\mathbf P_x\big(\tau_0^-\in t+\varDelta_\delta\big)}
 	{\mathbf P_x\big(
 		\overline X_{\tau_0^-}\in\beta(t+\varDelta_\delta)
 		\big)}
 	\cdot
 	\frac{	\mathbf E\big[
 		\#\mathcal A_x^{>R}
 		\big(\beta t-x+\varDelta_{\beta\delta_1}\big)
 		\big]}
 	{\mathbf P_x\big(\tau_0^-\in t+\varDelta_\delta\big)}.
 	\eeqnn
 	By Theorems~\ref{MainThm01} and~\ref{MainThm02} and
 	\eqref{eqn.30}, the first ratio on the right-hand side
 	goes to $1$ as $t\to\infty$.
 	By spatial homogeneity and  \eqref{eqn5.5} with $\delta_1=\delta_2=\delta$, the second ratio vanishes as $t\to\infty$ and then $R\to\infty$.
 	This proves the limit \eqref{eqn.601} with $i=2$
 	and completes the whole proof.
 \qed

 \subsection{Proof of Lemma~\ref{Lemma.MainThm04.transfer}}
 
 In this section, we provide the detailed proof for Lemma~\ref{Lemma.MainThm04.transfer} with the help of the following auxiliary asymptotic result on the occurrence time of the first big jump of $X$ before passing into $(-\infty,0)$. 
 More precisely, for any $x\geq 0$ and $h>0$, we define
 \beqnn
 \mathcal B_x^h:= \big\{ s \in (0,\tau_{-x}^-] : \Delta X_s>h \big\} 
 \quad \mbox{and}\quad 
 \mathcal J_x^h := \inf  \big\{ s \in (0,\tau_{-x}^-] : \Delta X_s>h \big\},
 \eeqnn 
 with $ \inf\emptyset=\infty$ for convention. 
 The one-big-jump principle requires that trajectories containing two positive jumps of macroscopic size contribute only at a smaller order.

    \begin{proposition}\label{Prop.203}
 	For any $t>0$ and $x,y,z>0$, we always have
 	\beqnn
 	\mathbf E_z\bigg[
 	\sum_{0<s<r\leq t\wedge\tau_0^-}
 	\mathbf{1}_{\{\Delta X_s>x,\,\Delta X_r>y\}}
 	\bigg]
 	\leq
 	t\cdot\mathbf E_z[\tau_0^-]\cdot
 	\overline{\nu}(x)\cdot\overline{\nu}(y).
 	\eeqnn
 \end{proposition}
 
 \proof
 By the Markov property, we first have
 \beqnn
 \mathbf E_z\bigg[
 \sum_{0<s<r\leq t\wedge\tau_0^-}
 \mathbf{1}_{\{\Delta X_s>x,\,\Delta X_r>y\}}
 \bigg]
 =
 \mathbf E_z\bigg[
 \sum_{0<s\leq t\wedge\tau_0^-}
 \mathbf{1}_{\{\Delta X_s>x\}}
 \cdot
 \mathbf E_{X_{s-}+\Delta X_s}\Big[
 \sum_{0<r\leq t-s}
 \mathbf{1}_{\{\tau_0^->r,\,\Delta X_r>y\}}
 \Big]
 \bigg].
 \eeqnn
 Using the compensation formula 
 \eqref{eqn.CompenstationX} successively yields that the right-hand side equals
 \beqnn
 \lefteqn{
 	\int_0^t ds
 	\int_0^\infty
 	\mathbf P_z\big(
 	X_s\in du,\tau_0^->s
 	\big)
 	\int_x^\infty
 	\mathbf E_{u+w}\bigg[
 	\sum_{0<r\leq t-s}
 	\mathbf{1}_{\{\tau_0^->r,\,\Delta X_r>y\}}
 	\bigg]
 	\,\nu(dw)
 }
 \cr
 \ar=\ar
 \int_0^t ds
 \int_0^\infty
 \mathbf P_z\big(
 X_s\in du,\tau_0^->s
 \big)
 \int_x^\infty\nu(dw)
 \int_0^{t-s}
 \mathbf P_{u+w}\big(
 \tau_0^->r
 \big)\,dr
 \cdot\overline{\nu}(y)
 \cr\cr
 \ar\leq\ar
 t\cdot
 \int_0^t
 \mathbf P_z\big(
 \tau_0^->s
 \big)\,ds
 \cdot
 \overline{\nu}(x)
 \cdot
 \overline{\nu}(y),
 \eeqnn
 which is smaller than
 $
 t\cdot\mathbf E_z[\tau_0^-]\cdot
 \overline{\nu}(x)\cdot\overline{\nu}(y).
 $
 The proof ends.
 \qed
 
 \begin{proposition}\label{Lemma.MainThm04.bigjump}
 	For any $x,\delta,R>0$ and $\lambda\in(0,\beta)$,  recall $\mathcal T_x$ defined in \eqref{eqnn1.9}, we have, as $t\to\infty$,
 	\beqnn
 	\Big| \mathbf P_x\big( \mathcal J^{\lambda t}_0\leq R \,\big|\, \tau_{0}^-\in t+\varDelta_\delta \big) - \mathbf P(\mathcal T_x\leq R) \Big| \to 0.
 	\eeqnn
 \end{proposition}
 \proof In view of  \eqref{eqnn1.9},  for any $K>0$, the triangle inequality gives
 \beqlb\label{eqn.6011}
 \lefteqn{ \Big|	\mathbf P_x\big(	\mathcal J^{\lambda t}_0\leq R	\,\big|\,	\tau_0^-\in t+\varDelta_\delta	\big)	-	\mathbf P(\mathcal T_x\leq R)	\Big|}\cr
 \ar=\ar\Bigg|\mathbf P_x\big(\mathcal J^{\lambda t}_0\leq R\,\big|\,\tau_0^-\in t+\varDelta_\delta\big)-\frac{1}{\mathbf E_x[\tau_0^-]}\int_0^R ds\int_{\mathbb R_+}\mathbf P_x\big(X_{s-}\in dz,\tau_0^->s\big)\Bigg|\cr
 \ar\leq\ar\mathbf P_x\big(\mathcal J^{\lambda t}_0\leq R,\,X_{\mathcal J^{\lambda t}_0-}>K\,\big|\,\tau_0^-\in t+\varDelta_\delta\big)+\frac{1}{\mathbf E_x[\tau_0^-]}\int_0^R\mathbf P_x\big(X_{s-}>K,\tau_0^->s\big)\,ds\cr
 \ar\ar+\Bigg|\mathbf P_x\big(\mathcal J^{\lambda t}_0\leq R,\,X_{\mathcal J^{\lambda t}_0-}\leq K\,\big|\,\tau_0^-\in t+\varDelta_\delta\big) -\frac{1}{\mathbf E_x[\tau_0^-]}\int_0^R\mathbf P_x\big(X_{s-}\leq K,\tau_0^->s\big)\,ds \Bigg|.
 \eeqlb
 Since $\mathbf E_x[\tau_0^-]<\infty$, the second term on the right-hand side of the last inequality goes to $0$ as $K\to\infty$. 
 For the first one, we first have  
 \beqnn
 \mathbf P_x\big(\mathcal J^{\lambda t}_0\leq R,\,X_{\mathcal J^{\lambda t}_0-}>K,\,\tau_0^-\in t+\varDelta_\delta\big)
 \leq \mathbf E_x\bigg[\sum_{0<s\leq R\wedge\tau_0^-}\mathbf 1_{\{X_{s-}>K,\,	\Delta X_s>\lambda t,\,	\tau_0^-\in t+\varDelta_\delta	\}}\bigg].
 \eeqnn
 By the Markov property and then the compensation formula \eqref{eqn.CompenstationX}, the last expectation equals
 \beqnn
 \int_0^R ds\int_K^\infty\mathbf P_x\big(X_{s-}\in dz,\tau_0^->s\big)\int_{\lambda t}^\infty\mathbf P_{z+w}\big(\tau_0^-\in t-s+\varDelta_\delta\big)\nu(dw),
 \eeqnn
 which further can be bounded by
 \beqnn
 \lefteqn{\sup_{s\in[0,R]}\sup_{z\geq K}\int_{\lambda t}^\infty\mathbf P_{z+w}\big(\tau_0^-\in t-s+\varDelta_\delta\big)\nu(dw)\cdot\int_0^R\mathbf P_x\big(X_{s-}>K,\tau_0^->s\big)ds} \cr
 \ar\leq\ar \sup_{z\geq K}\int_{\lambda t}^\infty\mathbf P_{z+w}\big(\tau_0^-\in t-R+\varDelta_{R+\delta}\big)\nu(dw)\cdot\int_0^\infty\mathbf P_x\big(X_{s-}>K,\tau_0^->s\big)ds.
 \eeqnn
 By Proposition 6.6 in \cite{WuXu2026b}, we have for sufficiently large $t$,
 \beqnn
 \mathbf P_x\Big(\mathcal J^{\lambda t}_0\leq R,\,X_{\mathcal J^{\lambda t}_0-}>K\,\Big|\,\tau_0^-\in t+\varDelta_\delta\Big)
 \ar\leq\ar C\cdot\frac{\nu\big(\beta t+\varDelta_\delta\big)}{	\mathbf P_x\big(	\tau_0^-\in t+\varDelta_\delta	\big)}\cdot\int_0^\infty\mathbf P_x\big(X_{s-}>K,\tau_0^->s\big)ds.
 \eeqnn
 Since $\mathbf E_x[\tau_0^-]<\infty$, the last integral vanishes as $K\to\infty$, which along with Theorem~\ref{MainThm01} yields 
 \beqnn
 \lim_{K\to\infty} \limsup_{t\to\infty}\mathbf P_x\Big(\mathcal J^{\lambda t}_0\leq R,\,X_{\mathcal J^{\lambda t}_0 -}>K\,\Big|\,\tau_0^-\in t+\varDelta_\delta\Big) =0.
 \eeqnn
 It remains to consider the last absolute value in \eqref{eqn.6011}. 
 For $R,K>0$, we first introduce
 \beqnn
 \mathcal B_{\leq R}^{\lambda t}
 := \{s \in (0,\tau_0^-\wedge R]: \Delta X_s>\lambda t  \} 
 \quad\mbox{and}\quad
 \mathcal B_{\leq R}^{\lambda t}(K)
 :=
 \mathcal B_{\leq R}^{\lambda t}
 \cap
 \big\{
 0<s\leq\tau_0^-:
 X_{s-}\leq K
 \big\}.
 \eeqnn
 On the event
 $
 \big\{
 \mathcal J^{\lambda t}_0\leq R,\,
 X_{\mathcal J^{\lambda t}_0-}\leq K,\,
 \tau_0^-\in t+\varDelta_\delta
 \big\},
 $
 the first jump of size greater than $\lambda t$ occurs before time
 $R$ and its pre-jump position is bounded by $K$. Hence this jump
 belongs to $\mathcal B_{\leq R}^{\lambda t}(K)$ and therefore
 \beqnn
 \big\{
 \mathcal J^{\lambda t}_0\leq R,\,
 X_{\mathcal J^{\lambda t}_0-}\leq K,\,
 \tau_0^-\in t+\varDelta_\delta
 \big\}\subset
 \big\{
 \#\mathcal B_{\leq R}^{\lambda t}(K)\geq1,\,
 \tau_0^-\in t+\varDelta_\delta
 \big\}.
 \eeqnn
 Conversely, suppose that
 $\#\mathcal B_{\leq R}^{\lambda t}(K)\geq1$, but
 $X_{\mathcal J^{\lambda t}_0-}>K$.
 Since the existence of a jump in
 $\mathcal B_{\leq R}^{\lambda t}(K)$ already implies that
 $\mathcal J^{\lambda t}_0\leq R$, there must then be another jump
 of size greater than $\lambda t$, occurring after
 $\mathcal J^{\lambda t}_0$ but no later than time $R$, whose
 pre-jump position is at most $K$. Consequently, there are at least
 two jumps of size greater than $\lambda t$ before time $R$. Thus
 \beqnn
 \big\{
 \#\mathcal B_{\leq R}^{\lambda t}(K)\geq1,\,
 \tau_0^-\in t+\varDelta_\delta
 \big\}
 \setminus
 \big\{
 \mathcal J^{\lambda t}_0\leq R,\,
 X_{\mathcal J^{\lambda t}_0-}\leq K,\,
 \tau_0^-\in t+\varDelta_\delta
 \big\}\subset
 \big\{
 \#\mathcal B_{\leq R}^{\lambda t}\geq2,\,
 \tau_0^-\in t+\varDelta_\delta
 \big\}.
 \eeqnn
 The preceding two containment relations imply that
 \beqnn
 \Big|
 \mathbf P_x\big(
 \mathcal J^{\lambda t}_0\leq R,\,
 X_{\mathcal J^{\lambda t}_0-}\leq K,\,
 \tau_0^-\in t+\varDelta_\delta
 \big)
 -
 \mathbf E_x\big[
 \#\mathcal B_{\leq R}^{\lambda t}(K);
 \tau_0^-\in t+\varDelta_\delta
 \big]
 \Big|
 \leq 
 \mathbf E_x\big[
 \#\mathcal B_{\leq R}^{\lambda t}
 \big(
 \#\mathcal B_{\leq R}^{\lambda t}-1
 \big)
 \big].
 \eeqnn
 Proposition~\ref{Prop.203} gives
 \beqnn
 \mathbf E_x\big[
 \#\mathcal B_{\leq R}^{\lambda t}
 \big(
 \#\mathcal B_{\leq R}^{\lambda t}-1
 \big)
 \big]=
 2\mathbf E_x\bigg[
 \sum_{0<s<r\leq R\wedge\tau_0^-}
 \mathbf 1_{\{
 	\Delta X_s>\lambda t,\,
 	\Delta X_r>\lambda t
 	\}}
 \bigg]
 \leq 
 2R\mathbf E_x[\tau_0^-]\,
 \overline\nu(\lambda t)^2,
 \eeqnn
 which is
 $
 o\big(
 \mathbf P_x(
 \tau_0^-\in t+\varDelta_\delta
 )
 \big)
 $
 as $t\to\infty$, by Theorem~\ref{MainThm01} and \eqref{eqn.30}.
 Consequently,
 \beqnn
 \limsup_{t\to\infty}
 \Bigg|
 \mathbf P_x\Big(
 \mathcal J^{\lambda t}_0\leq R,\,
 X_{\mathcal J^{\lambda t}_0-}\leq K
 \,\Big|\,
 \tau_0^-\in t+\varDelta_\delta
 \Big)
 -
 \frac{
 	\mathbf E_x\big[
 	\#\mathcal B_{\leq R}^{\lambda t}(K);
 	\tau_0^-\in t+\varDelta_\delta
 	\big]
 }{
 	\mathbf P_x(
 	\tau_0^-\in t+\varDelta_\delta
 	)
 }
 \Bigg|
 =0.
 \eeqnn
 Hence it suffices to prove
 \beqlb\label{eqn.P.901.third}
 \lim_{t\to\infty}
 \frac{
 	\mathbf E_x\big[
 	\#\mathcal B_{\leq R}^{\lambda t}(K);
 	\tau_0^-\in t+\varDelta_\delta
 	\big]
 }{
 	\mathbf P_x(
 	\tau_0^-\in t+\varDelta_\delta
 	)
 }
 =
 \frac{1}{\mathbf E_x[\tau_0^-]}
 \int_0^R
 \mathbf P_x\big(
 X_{s-}\leq K,\tau_0^->s
 \big)ds.
 \eeqlb
 Using the Markov property and then the compensation formula \eqref{eqn.CompenstationX} again,
 \beqnn
 \mathbf E_x\big[\#\mathcal B_{\leq R}^{\lambda t}(K);\tau_0^-\in t+\varDelta_\delta\big]
 \ar=\ar\int_0^R ds\int_{0}^{K}\mathbf P_x\big(X_{s-}\in dz,\tau_0^->s\big)\int_{\lambda t}^\infty\mathbf P_{z+w}\big(\tau_0^-\in t-s+\varDelta_\delta\big)\nu(dw).
 \eeqnn
 Then \eqref{eqn.P.901.third} follows directly from Proposition 6.6 in \cite{WuXu2026b}  and Theorem~\ref{MainThm01}.
 \qed

 \medskip
 
 \textbf{\textit{Proof of Lemma~\ref{Lemma.MainThm04.transfer}.}}
 For any fixed $R>0$ and all sufficiently large $t$, we have
 \beqlb\label{eqn.6012}
 \lefteqn{
 	\mathbf P\Big(
 	\#\mathcal A_{x}^{\leq R}\big(\beta t-x-\theta_t/2+\varDelta_{\theta_t}
 	\big)=0
 	\,\Big|\,
 	\tau_{-x}^-\in t+\varDelta_\delta
 	\Big)}
 \cr\cr
 \ar\leq\ar
 \mathbf P\Big(
 \#\mathcal A_{x}^{\leq R}\big(
 (0,\infty)
 \big)=0
 \,\Big|\,
 \tau_{-x}^-\in t+\varDelta_\delta
 \Big)
 \cr\cr
 \ar\ar+
 \mathbf P\Big(
 \#\mathcal A_{x}^{\leq R}\big(
 (0,\beta t-x-\theta_t/4]
 \big)
 +
 \#\mathcal A_{x}^{\leq R}\big(
 (\beta t-x+\theta_t/4,\infty)
 \big)
 \geq1
 \,\Big|\,
 \tau_{-x}^-\in t+\varDelta_\delta
 \Big).
 \eeqlb
 Consequently, it suffices to prove that both the two conditional probabilities on the right-hand side vanish as $t\to\infty$ and then $R\to\infty$.  
 
 \textbf{(i)}
 Fix $\lambda\in(0,\beta)$.  
 Let $e_{u^{\lambda t}}$ be the excursion that contains the first big jump larger than $\lambda t$.  
 For all sufficiently large $t$, on the event $\{\mathcal J_{x}^{\lambda t}<\tau_{-x}^-\}$ we have 
 \beqnn
 -\widehat H_{u^{\lambda t}-}+\bar e_{u^{\lambda t}}
 >
 \lambda t-x>0
 \quad\mbox{and}\quad
 \widehat L^{-1}_{{u^{\lambda t}}-}
 \leq\mathcal J_{x}^{\lambda t}.
 \eeqnn
 If $\mathcal J_{x}^{\lambda t}\leq R/4$ and $\tau_{-x}^--\widehat L^{-1}_{u^{\lambda t}}\leq R/2$, we have $ \widehat L^{-1}_{u^{\lambda t}-} +\tau_{-x}^- -\widehat L^{-1}_{u^{\lambda t}}<R$ and hence ${u^{\lambda t}}\in\mathcal A_{x}^{\leq R}((0,\infty))$.
 Thus, if $\#\mathcal A_{x}^{\leq R}((0,\infty))=0$, we see that either $\mathcal J_{x}^{\lambda t}>R/4$, or $\mathcal J_{x}^{\lambda t}\leq R/4$ and $\tau_{-x}^--\widehat L^{-1}_{u^{\lambda t}}>R/2$.
 Consequently, we have 
 \beqnn
 \lefteqn{\mathbf P\Big( \#\mathcal A_{x}^{\leq R}\big( (0,\infty) \big)=0\,\Big|\, \tau_{-x}^-\in t+\varDelta_\delta \Big)}\cr\cr
 \ar\leq\ar
 \mathbf P\Big(
 \mathcal J_{x}^{\lambda t}>R/4
 \,\Big|\,
 \tau_{-x}^-\in t+\varDelta_\delta
 \Big) +
 \mathbf P\Big(
 \mathcal J_{x}^{\lambda t}\leq R/4,\,
 \tau_{-x}^--\widehat L^{-1}_{u^{\lambda t}}>R/2
 \,\Big|\,
 \tau_{-x}^-\in t+\varDelta_\delta
 \Big).
 \eeqnn
 By the spatial homogeneity and Proposition~\ref{Lemma.MainThm04.bigjump}, the first conditional probability on the right-hand side vanishes as  $t\to\infty$ and then $R\to\infty$. 
 For the second one, we decompose it into  two terms:
 \beqnn
 P_1(t,R)\ar:=\ar \mathbf P\Big(\mathcal J_{x}^{\lambda t}\leq R/4,\, R/2<\tau_{-x}^--\widehat L^{-1}_{u^{\lambda t}}\leq t/2\,\Big|\,\tau_{-x}^-\in t+\varDelta_\delta\Big),\cr\cr
 P_2(t,R)\ar:=\ar\mathbf P\Big(\mathcal J_{x}^{\lambda t}\leq R/4,\,\tau_{-x}^--\widehat L^{-1}_{u^{\lambda t}}>t/2\,\Big|\,\tau_{-x}^-\in t+\varDelta_\delta\Big).
 \eeqnn
 We now prove that  both of them go to $0$ as $t\to\infty$ and then $R\to\infty$. 
 On the event $ \{\tau_{-x}^-\in t+\varDelta_\delta\}$, one can easily verify that 
 \beqnn
 \zeta(e_u)  =\tau_{-x}^--\widehat L^{-1}_{u-} - \big( \tau_{-x}^--\widehat L^{-1}_u \big)
 \in t-\widehat L^{-1}_{u-}-k-1+\varDelta_{\delta+1},
 \eeqnn 
 whenever $\widehat L^{-1}_{u-}\leq R/4 $ and $k<\tau_{-x}^--\widehat L^{-1}_u\leq k+1$.
 Consequently,
 \beqnn
 \lefteqn{ \mathbf P\Big(\mathcal J_{x}^{\lambda t}\leq R/4,\, R/2<\tau_{-x}^--\widehat L^{-1}_{u^{\lambda t}}\leq t/2,\, \tau_{-x}^-\in t+\varDelta_\delta \Big)} \cr\cr
 \ar\leq\ar \sum_{k=\lfloor R/2\rfloor}^{\lfloor t/2\rfloor} \mathbf E\Bigg[	\sum_{0\leq u<\widehat L_{\tau_{-x}^- -}}	\mathbf 1_{\{		e_u\neq\Upsilon,\,	\widehat L^{-1}_{u-}\leq R/4	\}}	\cdot \mathbf 1_{\{ 	k<\tau_{-x}^--\widehat L^{-1}_u\leq k+1	\}}	\cdot \mathbf 1_{\{	\tau_{-x}^-\in t+\varDelta_\delta	\}}	\Bigg]	\cr\cr
 \ar\leq\ar \sum_{k=\lfloor R/2\rfloor}^{\lfloor t/2\rfloor}	\mathbf E\Bigg[	\sum_{0\leq u<\widehat L_{\tau_{-x}^- -}}	\mathbf 1_{\{	e_u\neq\Upsilon,\,	\widehat L^{-1}_{u-}\leq R/4\}}\cdot	\mathbf 1_{\{\zeta(e_u)	\in	t-\widehat L^{-1}_{u-}-k-1+\varDelta_{\delta+1}	\}} \cdot \mathbf 1_{\{\tau_{-x}^--\widehat L^{-1}_u>k	\}}	\Bigg].	
 \eeqnn
 By the strong Markov property and the compensation formula \eqref{eqnn5.12}, the last sum can be bounded by
 \beqnn
 \sum_{k=\lfloor R/2\rfloor}^{\lfloor t/2\rfloor}	\mathbf P\big(	\tau_{-x}^->k	\big)	\mathbf E\Bigg[	\int_0^{\widehat L_{\tau_{-x}^- -}}	\mathbf 1_{\{\widehat L^{-1}_{u-}\leq R/4\}}	\underline n\Big(	\zeta\in	t-\widehat L^{-1}_{u-}-k-1+\varDelta_{\delta+1}	\Big)\,du	\Bigg], 
 \eeqnn
 which, by the fact $0\leq\widehat H_u\leq x$ for all $u<\widehat L_{\tau_{-x}^- -}$, further can be bounded by 
 \beqnn
 \mathbf E\big[	\widehat L_{\tau_{-x}^-}	\big]\cdot	\sup_{t/3\leq u\leq t}	\underline n\big(	\zeta\in u+\varDelta_{\delta+1}	\big)\cdot	\sum_{k=\lfloor R/2\rfloor}^{\lfloor t/2\rfloor}	\mathbf P\big(	\tau_{-x}^->k	\big). 
 \eeqnn
 Combining these estimates with \eqref{eqn.DefV} yields
 \beqnn
 P_1(t,R) \ar\leq\ar
 C\cdot \sup_{t/3\leq u\leq t}
 \frac{ \underline n\big(	\zeta\in u+\varDelta_{\delta+1}	\big)}{	\mathbf P_x\big(	\tau_{0}^-\in t+\varDelta_\delta	\big)}\cdot
 \sum_{k=\lfloor R/2\rfloor}^{\infty}
 \mathbf P\big(
 \tau_{-x}^->k
 \big).
 \eeqnn
 The ratio remains bounded as $t\to\infty$ by Theorem~1.3 in \cite{WuXu2026b}, Theorem~\ref{MainThm01} and \eqref{eqn.30}.
 Since $\mathbf E[\tau_{-x}^-]<\infty$, the last sum tends to $0$ as $R\to\infty$. Hence  $P_1(t,R)  \to 0$ as $t\to\infty$ and then $R\to\infty$. 
 For  $P_2(t,R)$,  since the excursion $e_{u^{\lambda t}}$ starts no later than $R/4$ and contains a jump larger than $\lambda t$, we have 
 \beqnn
 \lefteqn{	\mathbf P\Big(	\mathcal J_{x}^{\lambda t}\leq R/4,\,	\tau_{-x}^--\widehat L^{-1}_{u^{\lambda t}}>t/2,\,	\tau_{-x}^-\in t+\varDelta_\delta	\Big)} \cr\cr	
 \ar\leq\ar	\mathbf E\bigg[	\sum_{0\leq u<\widehat L_{\tau_{-x}^- -}}	\mathbf 1_{\{ 	\widehat L^{-1}_{u-}\leq R/4	\}} \cdot	\mathbf 1_{\{	\tau_{-x}^--\widehat L^{-1}_u>t/2	\}}	\sum_{0<s\leq\zeta(e_u)}	\mathbf 1_{\{\Delta e_u(s)>\lambda t\}}	\bigg].	
 \eeqnn
 By the strong Markov property, the compensation formulas \eqref{eqnn5.12} and \eqref{eqn.Compensation}, the last expectation can be bounded by
 \beqnn
 \lefteqn{\mathbf P\big(	\tau_{-x}^->t/2	\big)\cdot	\mathbf E\Bigg[	\sum_{0\leq u<\widehat L_{\tau_{-x}^- -}}	\mathbf 1_{\{	 	\widehat L^{-1}_{u-}\leq R/4		\}}	\sum_{0<s\leq\zeta(e_u)}	\mathbf 1_{\{\Delta e_u(s)>\lambda t\}}	\Bigg]}	\cr\cr
 \ar\leq\ar	\mathbf P_x\big(	\tau_{0}^->t/2	\big) \cdot	\mathbf E\Bigg[	\int_0^{\widehat L_{\tau_{-x}^- -}}	\mathbf 1_{\{\widehat L^{-1}_{u-}\leq R/4\}}	\underline n\Bigg(	\sum_{0<s\leq\zeta}	\mathbf 1_{\{\Delta\epsilon_s>\lambda t\}}	\Bigg)\,du	\Bigg]	
 \leq C\cdot	\mathbf P_x\big(	\tau_{0}^->t/2	\big)\cdot	 	\overline\nu(\lambda t).
 \eeqnn
 Combining these estimates with   \eqref{eqn.501} yields that  for any large $R$, 
 \beqnn
 P_2(t,R)\leq C\cdot	\frac{	\mathbf P_x\big(	\tau_{0}^->t/2	\big)\cdot\overline\nu(\lambda t)	}{	\mathbf P_x\big(	\tau_{0}^-\in t+\varDelta_\delta	\big)	},
 \eeqnn 
 which goes to $0$ as $t\to\infty$ by \eqref{Asytau}, Theorem~\ref{MainThm01} and \eqref{eqn.30}. 
 Here we have proved that the first conditional probability on the right-hand side of \eqref{eqn.6012} vanishes as $t\to\infty$ and then $R\to\infty$.  
 
 \textbf{(ii)} We now start to prove that for any fixed large $R>0$,  the second conditional probability on the right-hand side of \eqref{eqn.6012} goes to $0$ as $t\to\infty$. Firstly, we have 
 \beqlb\label{eqn.6013}
 \lefteqn{	\mathbf P\Big(	\#\mathcal A_{x}^{\leq R}\big(	(0,\beta t-x-\theta_t/4]	\big)	+	\#\mathcal A_{x}^{\leq R}\big(	(\beta t-x+\theta_t/4,\infty)	\big)	\geq1,\,	\tau_{-x}^-\in t+\varDelta_\delta	\Big)} \cr\cr
 \ar\leq\ar  \mathbf E\Big[ \#\mathcal A_{x}^{\leq R}\big((0,\beta t-x-\theta_t/4]\big) + \#\mathcal A_{x}^{\leq R}\big((\beta t-x+\theta_t/4,\infty)\big) ;\tau_{-x}^-\in t+\varDelta_\delta\Big].
 \eeqlb
 On the event $\{\tau_{-x}^-\in t+\varDelta_\delta\}$,  we have for any $u\in\mathcal A_{x}^{\leq R}\big((0,\beta t-x-\theta_t/4]\big)\cup\mathcal A_{x}^{\leq R}\big((\beta t-x+\theta_t/4,\infty)\big)$, 
 \beqnn
 \zeta(e_u) =\tau_{-x}^--\big(\widehat L^{-1}_{u-}+\tau_{-x}^--\widehat L^{-1}_u\big)\in t-R+\varDelta_{R+\delta}.
 \eeqnn
 Moreover, since $0\leq\widehat H_{u-}\leq x$ and $\theta_t\to\infty$, for all sufficiently large $t$, we also have 
 \beqnn
 \big| \bar e_u- \beta(t-R)\big|\ar\geq\ar\big| x-\widehat H_{u-}+\bar e_u -\beta t \big| - \big| x-\widehat H_{u-}-\beta R \big| \geq  \theta_t/4-x-\beta R > \theta_t/8.
 \eeqnn
 By the preceding two estimates and the compensation formula \eqref{eqnn5.12},  the expectation on the right-hand side of \eqref{eqn.6013} can be bounded by 
 \beqnn 
 \lefteqn{\mathbf E\Bigg[\sum_{0\leq u\leq\widehat L_{\tau_{-x}^- -}}\mathbf 1_{\{ 	\zeta(e_u)\in	t-R+\varDelta_{R+\delta}	\}} \cdot \mathbf 1_{\{|\bar e_u-\beta(t-R)|>\theta_t/8	\}}\Bigg] }\cr\cr
 \ar=\ar	\mathbf E\Big[\widehat L_{\tau_{-x}^- -}\Big]\cdot\underline n\Big(|\bar\epsilon-\beta(t-R)|>\theta_t/8,\,\zeta\in t-R+\varDelta_{R+\delta}\Big).
 \eeqnn
 Combining the preceding estimates together with \eqref{eqn.DefV} yields 
 \beqnn
 \ar\ar \mathbf P\Big(
 \#\mathcal A_{x}^{\leq R}\big(
 (0,\beta t-x-\theta_t/4]
 \big)
 +
 \#\mathcal A_{x}^{\leq R}\big(
 (\beta t-x+\theta_t/4,\infty)
 \big)
 \geq1
 \,\Big|\,
 \tau_{-x}^-\in t+\varDelta_\delta
 \Big)\cr\cr
 \ar\ar \leq C\cdot
 \frac{
 	\underline n\big(
 	\zeta\in
 	t-R+\varDelta_{R+\delta}
 	\big)
 }{
 	\mathbf P_x\big(
 	\tau_{0}^-\in t+\varDelta_\delta
 	\big)
 }\cdot
 \underline n\Big(
 |\bar\epsilon-\beta(t-R)|
 >\theta_t/8
 \,\Big|\,
 \zeta\in
 t-R+\varDelta_{R+\delta}
 \Big),
 \eeqnn
 which goes to $0$ as $t\to\infty$, since  the last conditional probability tends to $0$ by Theorem~1.6 in \cite{WuXu2026b} and the ratio can be bounded uniformly for all sufficiently large $t$ by Theorem~1.3 in  \cite{WuXu2026b}, Theorem~\ref{MainThm01} and \eqref{eqn.30}. 
 This completes the whole proof. 
 \qed

  \section{Proofs of Theorem~\ref{MainThm05}}
  \label{Sec.ProofThm05}
  \setcounter{equation}{0}

  In this section, we prove the two conditional scaling limits
  stated in Theorem~\ref{MainThm05} under $\mathbf P_x$.
  As explained in detail in the Introduction, a path with a large
  first-passage time or a large maximum is typically generated by
  a sufficiently large jump occurring at an early stage.
  Moreover, the independent-increment property of the underlying
  L\'evy process suggests decomposing the path at this large jump
  and analyzing the corresponding pre- and post-jump contributions
  separately, i.e., for any fixed $\lambda\in(0,\beta)$, recall  $\mathcal{J}^{\lambda t}_0$  the occurrence time of first big-jump with size larger than $\lambda t$, 
  \beqnn
  X_s
  =
  \big(
  X_s-\Delta X_{\mathcal{J}^{\lambda t}_0}
  \cdot\mathbf 1_{\{s\geq\mathcal{J}^{\lambda t}_0\}}
  \big)
  +
  \Delta X_{\mathcal{J}^{\lambda t}_0}
  \cdot\mathbf 1_{\{s\geq\mathcal{J}^{\lambda t}_0\}},
  \qquad
  s\in[0,\tau_0^-].
  \eeqnn
  On the event $\{\mathcal{J}^{\lambda t}_0\geq\tau_0^-\}$, we set
  $\Delta X_{\mathcal{J}^{\lambda t}_0}=0$ and
  $X_{s\vee\mathcal{J}^{\lambda t}_0}-X_{\mathcal{J}^{\lambda t}_0}=0$
  for the expressions appearing below.
  The following lemmas establish the asymptotic behavior of the
  first occurrence time of the large jump, as well as that of the
  pre-jump and post-jump path segments, respectively.
  
  \begin{lemma}\label{Lemma.901}
  	For any constants $\delta,R>0$ and $\lambda\in(0,\beta)$, we have as
  	$t\to\infty$,
  	\beqnn
  	\Big|
  	\mathbf P_x\big(
  	\mathcal{J}^{\lambda t}_0\leq R
  	\,\big|\,
  	\tau_0^-\in t+\varDelta_\delta
  	\big)
  	-
  	\mathbf P(\mathcal T_x\leq R)
  	\Big|
  	+
  	\Big|
  	\mathbf P_x\big(
  	\mathcal{J}^{\lambda t}_0\leq R
  	\,\big|\,
  	\overline{X}_{\tau_0^-}\in\beta(t+\varDelta_\delta)
  	\big)
  	-
  	\mathbf P(\mathcal T_x\leq R)
  	\Big|
  	\to0.
  	\eeqnn
  \end{lemma}
  
  \begin{lemma}\label{Lemma.902}
  	For any constants $\lambda\in(0,\beta)$ and $\delta>0$, for any
  	$\varepsilon>0$, we have as $t\to\infty$,
  	\beqnn
  	\mathbf P_x\bigg(
  	\Big|
  	\frac{\Delta  X_{\mathcal{J}^{\lambda t}_0}}{t}-\beta
  	\Big|\geq\varepsilon
  	\,\bigg|\,
  	\tau_0^-\in t+\varDelta_\delta
  	\bigg)
  	+
  	\mathbf P_x\bigg(
  	\Big|
  	\frac{\Delta  X_{\mathcal{J}^{\lambda t}_0}}{t}-\beta
  	\Big|\geq\varepsilon
  	\,\bigg|\,
  	\overline{X}_{\tau_0^-}\in\beta(t+\varDelta_\delta)
  	\bigg)
  	\to0.
  	\eeqnn
  \end{lemma}
  
  \begin{lemma}\label{Lemma.903}
  	For any constants $\lambda\in(0,\beta)$ and $\delta,T>0$, for any
  	$\varepsilon>0$, we have as $t\to\infty$,
  	\beqnn
  	\mathbf P_x\bigg(
  	\sup_{s\leq T}
  	\frac{X_{(s\wedge\mathcal{J}^{\lambda t}_0)-}}{t}
  	\geq\varepsilon
  	\,\bigg|\,
  	\tau_0^-\in t+\varDelta_\delta
  	\bigg)
  	+
  	\mathbf P_x\bigg(
  	\sup_{s\leq T\wedge\tau_0^-}
  	\frac{X_{(s\wedge\mathcal{J}^{\lambda t}_0)-}}{t}
  	\geq\varepsilon
  	\,\bigg|\,
  	\overline{X}_{\tau_0^-}\in\beta(t+\varDelta_\delta)
  	\bigg)
  	\to0.
  	\eeqnn
  \end{lemma}

  \begin{lemma}\label{Lemma.904}
  	For any constants $\lambda\in(0,\beta)$ and $\delta,T,x>0$, for any
  	$\varepsilon>0$, we have as $t\to\infty$,
  	\beqnn
  	\mathbf P_x\bigg(
  	\sup_{s\leq T}
  	\Big|
  	\frac{
  		X_{s\vee\mathcal{J}^{\lambda t}_0}
  		-
  		X_{\mathcal{J}^{\lambda t}_0}
  	}{t}
  	\Big|
  	\geq\varepsilon
  	\,\bigg|\,
  	\tau_0^-\in t+\varDelta_\delta
  	\bigg) \to 0
  	\eeqnn
  	and 
  	\beqnn
  	\mathbf P_x\bigg(		\sup_{s\leq T\wedge \tau_0^-}
  	\Big|
  	\frac{
  		X_{s\vee\mathcal{J}^{\lambda t}_0}
  		-
  		X_{\mathcal{J}^{\lambda t}_0}
  	}{t}
  	\Big|
  	\geq\varepsilon
  	\,\bigg|\,
  	\overline{X}_{\tau_0^-}\in\beta(t+\varDelta_\delta)
  	\bigg) \to0.
  	\eeqnn
  \end{lemma}

  \textbf{\textit{Proof of Theorem~\ref{MainThm05}(1).}} Here we just give a detailed proof for the uniform weak convergence under $\mathbf P_x\big(\cdot\,|\, \tau_0^-\in t+\varDelta_\delta\big)$.
  The case of $\mathbf P_x\big(\cdot\,|\, \overline X_{\tau_0^-}\in\beta(t+\varDelta_\delta)\big)$ can be proved similarly.
  
  \textbf{Step 1.} We first prove the weak convergence for any fixed $\delta>0$.
  By Lemma~\ref{Lemma.901}, we have $\mathbf P_x\big(\mathcal{J}^{\lambda t}_0\geq\tau_0^-\,\big|\,\tau_0^-\in t+\varDelta_\delta\big)\to0$ as $t\to\infty$.
  For any fixed $\lambda\in(0,\beta)$, we can decompose $X$ as follows:
  \beqnn
  X_s
  =
  \big( X_s-\Delta X_{\mathcal{J}^{\lambda t}_0}\cdot\mathbf 1_{\{s\geq\mathcal{J}^{\lambda t}_0\}} \big)
  +
  \Delta X_{\mathcal{J}^{\lambda t}_0}\cdot\mathbf 1_{\{s\geq\mathcal{J}^{\lambda t}_0\}},
  \quad s\in[0,\tau_0^-].
  \eeqnn
  By using Slutsky's theorem along with Lemma~\ref{Lemma.901} and \ref{Lemma.902}, we have $\big(\mathcal{J}^{\lambda t}_0,\Delta X_{\mathcal{J}^{\lambda t}_0}/t\big)$ under $\mathbf P_x\big(\cdot\,|\,\tau_0^-\in t+\varDelta_\delta\big)$ converges to $(\mathcal T_x,\beta)$ as $t\to\infty$, which induces that $\{ t^{-1}\cdot \Delta  X_{\mathcal{J}^{\lambda t}_0}\cdot\mathbf 1_{\{s\geq\mathcal{J}^{\lambda t}_0\}}:s\geq0\}$ converges weakly to $\{\beta\cdot\mathbf 1_{\{s\geq\mathcal T_x\}}:s\geq0\}$ in $D(\mathbb R_+;\mathbb R)$; see Example~1.19 in \cite[p.329]{JacodShiryaev2003}.
  By Theorem 16.4 in \cite[p.~170]{Billingsley1999}, it remains to prove that for any $T,\varepsilon>0$, as $t\to\infty$,
  \beqlb\label{eqn.906}
  \mathbf P_x\Big(
  \sup_{s\leq T}
  \big|X_s-\Delta X_{\mathcal{J}^{\lambda t}_0}\cdot\mathbf 1_{\{s\geq\mathcal{J}^{\lambda t}_0\}}\big|
  \geq\varepsilon t,\,
  \mathcal{J}^{\lambda t}_0<\tau_0^-
  \,\Big|\,
  \tau_0^-\in t+\varDelta_\delta
  \Big)
  \to0.
  \eeqlb
  On the event $\{\mathcal{J}^{\lambda t}_0<\tau_0^-\}$, we have for any $s\leq T$,
  \beqnn
  X_s-\Delta X_{\mathcal{J}^{\lambda t}_0}\cdot\mathbf 1_{\{s\geq\mathcal{J}^{\lambda t}_0\}}
  =
  \big( X_{s\wedge\mathcal{J}^{\lambda t}_0}-\Delta X_{\mathcal{J}^{\lambda t}_0}\cdot\mathbf 1_{\{s\geq\mathcal{J}^{\lambda t}_0\}} \big)
  +
  \big( X_{s\vee\mathcal{J}^{\lambda t}_0}-X_{\mathcal{J}^{\lambda t}_0} \big),
  \eeqnn
  which follows that
  \beqnn
  \sup_{s\leq T}
  \big|X_s-\Delta X_{\mathcal{J}^{\lambda t}_0}\cdot\mathbf 1_{\{s\geq\mathcal{J}^{\lambda t}_0\}}\big|
  \leq
  \sup_{s\leq T}X_{(s\wedge\mathcal{J}^{\lambda t}_0)-}
  +
  \sup_{s\leq T}\big|X_{s\vee\mathcal{J}^{\lambda t}_0}-X_{\mathcal{J}^{\lambda t}_0}\big|.
  \eeqnn
  This induces that the left-hand side of \eqref{eqn.906} can be bounded by
  \beqnn
  \mathbf P_x\Big(
  \sup_{s\leq T}X_{(s\wedge\mathcal{J}^{\lambda t}_0)-}
  \geq\frac{\varepsilon t}{2}
  \,\Big|\,
  \tau_0^-\in t+\varDelta_\delta
  \Big)
  +
  \mathbf P_x\Big(
  \sup_{s\leq T}\big|X_{s\vee\mathcal{J}^{\lambda t}_0}-X_{\mathcal{J}^{\lambda t}_0}\big|
  \geq\frac{\varepsilon t}{2}
  \,\Big|\,
  \tau_0^-\in t+\varDelta_\delta
  \Big),
  \eeqnn
  which tends to $0$ as $t\to\infty$ by Lemma~\ref{Lemma.903} and \ref{Lemma.904}.
  
  \textbf{Step 2.} We extend the convergence uniformly in
  $\delta\in[\theta_0,\theta_t]$. It suffices to prove that for any $T>0$ and any bounded, non-negative, uniformly continuous functional $F$ on $D\big([0,T];\mathbb R\big)$, as $t\to\infty$,
  \beqlb\label{eqn.9009}
  \sup_{\delta\in[\theta_0,\theta_t]}
  \Big|
  \mathbf E_x\big[ F(X_s/t:s\in[0,T])\,\big|\,\tau_0^-\in t+\varDelta_\delta \big]
  -
  \mathbf E\big[ F(\beta\cdot\mathbf 1_{\{s\geq\mathcal T_x\}}:s\in[0,T])\big]
  \Big|
  \to0.
  \eeqlb
  For any large integer $n\geq1$, we have 
  \beqnn
  \mathbf E_x\big[ F(X_s/t:s\in[0,T]);\,\tau_0^-\in t+\varDelta_\delta \big] 
  \ar \leq\ar \sum_{k=0}^{[n\delta]} \mathbf E_x\Big[ F(X_s/t:s\in[0,T]),\, \tau_0^-\in t+\frac{k}{n}+\varDelta_{1/n} \Big],
  \eeqnn
  which further can be decomposed into the following two terms:
  \beqnn
  I_n(t,\delta) \ar:=\ar  \sum_{k=0}^{[n\delta]} \mathbf E_x\Big[ F\Big(\frac{X_s}{t+k/n}:s\in[0,T]\Big); \, \tau_0^-\in t+\frac{k}{n}+\varDelta_{1/n} \Big],
  \cr\cr
  \varepsilon_n(t,\delta) \ar:=\ar  \sum_{k=0}^{[n\delta]} \mathbf E_x\Big[ F\Big(\frac{X_s}{t}:s\in[0,T]\Big) - F\Big(\frac{X_s}{t+k/n}:s\in[0,T]\Big);\, \tau_0^-\in t+\frac{k}{n}+\varDelta_{1/n} \Big].
  \eeqnn
  Recall that the Skorokhod topology is weaker than the uniform topology.
  The boundedness and uniform continuity of $F$ imply that there exists
  a constant $C>0$, independent of $\epsilon$ and $t,k,n$ such that
  \beqnn
  \bigg|
  F\Big(\frac{X_s}{t}:s\in[0,T]\Big)
  -
  F\Big(\frac{X_s}{t+k/n}:s\in[0,T]\Big)
  \bigg|
  \ar\leq\ar
  C\cdot\mathbf 1_{\{\overline X_{\tau_0^-}\geq2\beta t\}}
  +
  C\cdot\frac{k/n}{t^2}\cdot\overline X_{\tau_0^-}\cdot\mathbf 1_{\{\overline X_{\tau_0^-}\leq2\beta t\}}
  \cr\cr
  \ar\leq\ar
  C\cdot\Big(
  \mathbf 1_{\{\overline X_{\tau_0^-}\geq2\beta t\}}
  +
  \frac{k/n}{t}
  \Big).
  \eeqnn
  Applying this and the law of total probability to $|\varepsilon_n(t,\delta)|$ yields
  \beqlb\label{eqn.9012}
  \sup_{\delta\in[\theta_0,\theta_t]}\frac{|\varepsilon_n(t,\delta)|}{\mathbf P_x\big(\tau_0^-\in t+\varDelta_\delta\big)}
  \ar\leq\ar
  C\cdot\sup_{\delta\in[\theta_0,\theta_t]}
  \frac{
  	\mathbf P_x\big(
  	\overline X_{\tau_0^-}\geq2\beta t,\,
  	\tau_0^-\in t+\varDelta_{([n\delta]+1)/n}
  	\big)
  }{
  	\mathbf P_x\big(\tau_0^-\in t+\varDelta_\delta\big)
  }
  \cr\cr
  \ar\ar+
  \frac{C\cdot\theta_t}{t}\cdot 
  \sup_{\delta\in[\theta_0,\theta_t]}
  \frac{
  	\mathbf P_x\big(\tau_0^-\in t+\varDelta_{([n\delta]+1)/n}\big)
  }{
  	\mathbf P_x\big(\tau_0^-\in t+\varDelta_\delta\big)
  },
  \eeqlb
  which goes to $0$ as $t\to\infty$ by using Theorem~\ref{MainThm04} and~\ref{MainThm01} to the two terms on the right-hand side respectively. 
  For $I_n(t,\delta)$, our result established in the first step states that there exists a positive function $\rho_n$ on $\mathbb R_+$
  with $\rho_n(t)\downarrow0$ as $t\to\infty$ such that for all $t>0$ and $k\geq0$,
  \beqnn
  \mathbf E_x\Big[
  F\Big(\frac{X_s}{t+k/n}:s\in[0,T]\Big)
  \,\Big|\,
  \tau_0^-\in t+\frac{k}{n}+\varDelta_{1/n}
  \Big]
  \leq
  \mathbf E\big[ F(\beta\cdot\mathbf 1_{\{s\geq\mathcal T_x\}}:s\in[0,T])\big]
  +
  \rho_n(t).
  \eeqnn
  Equivalently,
  \beqnn
  \lefteqn{
  	\mathbf E_x\Big[
  	F\Big(\frac{X_s}{t+k/n}:s\in[0,T]\Big),\,
  	\tau_0^-\in t+\frac{k}{n}+\varDelta_{1/n}
  	\Big]
  }
  \cr\cr
  \ar\leq\ar
  \Big(
  \mathbf E\big[ F(\beta\cdot\mathbf 1_{\{s\geq\mathcal T_x\}}:s\in[0,T])\big]
  +
  \rho_n(t)
  \Big)
  \cdot
  \mathbf P_x\Big(
  \tau_0^-\in t+\frac{k}{n}+\varDelta_{1/n}
  \Big).
  \eeqnn
  Substituting this back into the summation in $I_n(t,\delta)$ yields
  \beqnn
  \frac{I_n(t,\delta)}{\mathbf P_x\big(\tau_0^-\in t+\varDelta_\delta\big)}
  \ar\leq\ar
  \sum_{k=0}^{[n\delta]}
  \Big(
  \mathbf E\big[ F(\beta\cdot\mathbf 1_{\{s\geq\mathcal T_x\}}:s\in[0,T])\big]
  +
  \rho_n(t)
  \Big)
  \cdot
  \frac{
  	\mathbf P_x\big(\tau_0^-\in t+k/n+\varDelta_{1/n}\big)
  }{
  	\mathbf P_x\big(\tau_0^-\in t+\varDelta_\delta\big)
  }
  \cr\cr
  \ar=\ar
  \Big(
  \mathbf E\big[ F(\beta\cdot\mathbf 1_{\{s\geq\mathcal T_x\}}:s\in[0,T])\big]
  +
  \rho_n(t)
  \Big)
  \cdot
  \frac{
  	\mathbf P_x\big(\tau_0^-\in t+\varDelta_{([n\delta]+1)/n}\big)
  }{
  	\mathbf P_x\big(\tau_0^-\in t+\varDelta_\delta\big)
  }.
  \eeqnn
  This induces that
  \beqnn
  \lefteqn{\limsup_{t\to\infty}
  	\sup_{\delta\in[\theta_0,\theta_t]}
  	\frac{I_n(t,\delta)}{\mathbf P_x\big(\tau_0^-\in t+\varDelta_\delta\big)}}\cr\cr
  \ar\leq\ar
  \mathbf E\big[
  F(\beta\cdot\mathbf 1_{\{s\geq\mathcal T_x\}}:s\in[0,T])
  \big]
  \cdot
  \limsup_{t\to\infty}
  \sup_{\delta\in[\theta_0,\theta_t]}
  \frac{
  	\mathbf P_x\big(\tau_0^-\in t+\varDelta_{([n\delta]+1)/n}\big)
  }{
  	\mathbf P_x\big(\tau_0^-\in t+\varDelta_\delta\big)
  }
  \cr\cr
  \ar\leq\ar
  \mathbf E\big[
  F(\beta\cdot\mathbf 1_{\{s\geq\mathcal T_x\}}:s\in[0,T])
  \big]
  \cdot
  \Big(1+\frac{1/n}{\theta_0}\Big).
  \eeqnn 
  Combining this together with \eqref{eqn.9012} and then letting $n\to\infty$, we obtain
  \beqnn
  \limsup_{t\to\infty}\sup_{\delta\in[\theta_0,\theta_t]}
  \mathbf E_x\big[ F(X_s/t:s\in[0,T])\,\big|\,\tau_0^-\in t+\varDelta_\delta \big]
  \leq
  \mathbf E\big[ F(\beta\cdot\mathbf 1_{\{s\geq\mathcal T_x\}}:s\in[0,T])\big].
  \eeqnn
  On the other hand, using the sum from $k=0$ to $[n\delta]-1$, an analogous argument gives
  \beqnn
  \liminf_{t\to\infty}\inf_{\delta\in[\theta_0,\theta_t]}
  \mathbf E_x\big[ F(X_s/t:s\in[0,T])\,\big|\,\tau_0^-\in t+\varDelta_\delta \big]
  \ar\geq\ar
  \mathbf E\big[ F(\beta\cdot\mathbf 1_{\{s\geq\mathcal T_x\}}:s\in[0,T])\big].
  \eeqnn
  Combining them together yields \eqref{eqn.9009} immediately.
  \qed

  \medskip \textbf{\textit{Proof of Theorem~\ref{MainThm05}(2).}}
  Similarly to the proof of claim (1), here we still just prove that
  the weak convergence in claim (2) under
  $\mathbf P_x\big(\cdot\,|\,\tau_0^-\in t+\varDelta_\delta\big)$
  holds for any fixed $\delta$.
  It can be extended uniformly in $\delta\in[\theta_0,\theta_t]$ as before.
  The uniform weak convergence under
  $\mathbf P_x\big(\cdot\,|\,\overline X_{\tau_0^-}
  \in\beta(t+\varDelta_\delta)\big)$
  can be proved in the same way. 
  
  We start to prove the weak convergence under
  $\mathbf P_x\big(\cdot\,|\,\tau_0^-\in t+\varDelta_\delta\big)$.
  By Lemma~\ref{Lemma.901}, it suffices to prove that for any
  $\varepsilon\in(0,1)$, $0<s_1<s_2<\infty$ and
  $0<\eta\leq\varepsilon(s_1\wedge1)/(4\beta+1)<1$,
  as $t\to\infty$,
  \beqlb\label{eqn.9013}
  \mathbf P_x\Big(
  \sup_{s_1\leq s\leq s_2}
  \Big|
  \frac{X_{ts}}{t}-\beta(1-s)
  \Big|
  \geq\varepsilon,\,
  \mathcal{J}^{\lambda t}_0\leq\eta t
  \,\Big|\,
  \tau_0^-\in t+\varDelta_\delta
  \Big)
  \to0.
  \eeqlb
  On the event
  $\{\mathcal{J}^{\lambda t}_0\leq\eta t<\tau_0^-\}$,
  we observe that for any $s\in[s_1,s_2]$,
  \beqnn
  \frac{X_{ts}}{t}
  \ar=\ar
  \Big(
  \frac{X_{\mathcal{J}^{\lambda t}_0}}{t}-\beta
  \Big)
  +
  \Big[
  \frac{X_{ts}-X_{\mathcal{J}^{\lambda t}_0}}{t}
  +
  \beta\Big(
  s-\frac{\mathcal{J}^{\lambda t}_0}{t}
  \Big)
  \Big]
  +
  \beta(1-s)
  +
  \beta\cdot\frac{\mathcal{J}^{\lambda t}_0}{t},
  \eeqnn
  which along with the triangle inequality yields
  \beqnn
  \Big|
  \frac{X_{ts}}{t}-\beta(1-s)
  \Big|
  \ar\leq\ar
  \Big|
  \frac{X_{\mathcal{J}^{\lambda t}_0}}{t}-\beta
  \Big|
  +
  \Big|
  \frac{X_{ts}-X_{\mathcal{J}^{\lambda t}_0}}{t}
  +
  \beta\cdot\frac{ts-\mathcal{J}^{\lambda t}_0}{t}
  \Big|
  +
  \beta\cdot\frac{\mathcal{J}^{\lambda t}_0}{t}.
  \eeqnn
  Consequently, the left-hand side of \eqref{eqn.9013} can be bounded
  by the sum of the following three terms:
  \beqnn
  A_1(t)
  \ar:=\ar
  \mathbf P_x\Big(
  \Big|
  \frac{X_{\mathcal{J}^{\lambda t}_0}}{t}-\beta
  \Big|
  \geq\frac{\varepsilon}{4},\,
  \mathcal{J}^{\lambda t}_0\leq\eta t
  \,\Big|\,
  \tau_0^-\in t+\varDelta_\delta
  \Big),
  \cr\cr
  A_2(t)
  \ar:=\ar
  \mathbf P_x\Big(
  \sup_{s_1\leq s\leq s_2}
  \Big|
  \frac{X_{ts}-X_{\mathcal{J}^{\lambda t}_0}}{t}
  +
  \beta\cdot\frac{ts-\mathcal{J}^{\lambda t}_0}{t}
  \Big|
  \geq\frac{\varepsilon}{4},\,
  \mathcal{J}^{\lambda t}_0\leq\eta t
  \,\Big|\,
  \tau_0^-\in t+\varDelta_\delta
  \Big),
  \cr\cr
  A_3(t)
  \ar:=\ar
  \mathbf P_x\Big(
  \beta\cdot\frac{\mathcal{J}^{\lambda t}_0}{t}
  \geq\frac{\varepsilon}{4},\,
  \mathcal{J}^{\lambda t}_0\leq\eta t
  \,\Big|\,
  \tau_0^-\in t+\varDelta_\delta
  \Big).
  \eeqnn
  It is obvious that $A_3(t)=0$ since
  $\beta\eta\leq\beta\varepsilon(s_1\wedge1)/(4\beta+1)
  <\varepsilon/4$.
  For $A_1(t)$, by the triangle inequality we have
  \beqnn
  A_1(t)
  \ar\leq\ar
  \mathbf P_x\Big(
  \Big|
  \frac{\Delta X_{\mathcal{J}^{\lambda t}_0}}{t}-\beta
  \Big|
  \geq\frac{\varepsilon}{8},\,
  \mathcal{J}^{\lambda t}_0\leq\eta t
  \,\Big|\,
  \tau_0^-\in t+\varDelta_\delta
  \Big)\cr
  \ar\ar +
  \mathbf P_x\Big(
  \frac{X_{\mathcal{J}^{\lambda t}_0-}}{t}
  \geq\frac{\varepsilon}{8},\,
  \mathcal{J}^{\lambda t}_0\leq\eta t
  \,\Big|\,
  \tau_0^-\in t+\varDelta_\delta
  \Big).
  \eeqnn
  The first term goes to $0$ as $t\to\infty$ by
  Lemma~\ref{Lemma.902}.
  For the second term, we choose $R>0$ and obtain that
  for all $t>R/\eta$,
  \beqnn
  \lefteqn{
  	\mathbf P_x\Big(
  	\frac{X_{\mathcal{J}^{\lambda t}_0-}}{t}
  	\geq\frac{\varepsilon}{8},\,
  	\mathcal{J}^{\lambda t}_0\leq\eta t
  	\,\Big|\,
  	\tau_0^-\in t+\varDelta_\delta
  	\Big)}
  \cr\cr
  \ar\leq\ar
  \mathbf P_x\Big(
  \frac{X_{\mathcal{J}^{\lambda t}_0-}}{t}
  \geq\frac{\varepsilon}{8},\,
  \mathcal{J}^{\lambda t}_0\leq R
  \,\Big|\,
  \tau_0^-\in t+\varDelta_\delta
  \Big)
  +
  \mathbf P_x\big(
  R<\mathcal{J}^{\lambda t}_0\leq\eta t
  \,\big|\,
  \tau_0^-\in t+\varDelta_\delta
  \big)
  \cr\cr
  \ar\leq\ar
  \mathbf P_x\Big(
  \sup_{s\leq R}
  \frac{X_{(s\wedge\mathcal{J}^{\lambda t}_0)-}}{t}
  \geq\frac{\varepsilon}{8}
  \,\Big|\,
  \tau_0^-\in t+\varDelta_\delta
  \Big)
  +
  \mathbf P_x\big(
  \mathcal{J}^{\lambda t}_0>R
  \,\big|\,
  \tau_0^-\in t+\varDelta_\delta
  \big),
  \eeqnn
  which tends to $0$ as $t\to\infty$ and then $R\to\infty$
  by Lemmas~\ref{Lemma.901} and~\ref{Lemma.903}.
  Consequently, we have $A_1(t)\to0$ as $t\to\infty$. It remains to prove $A_2(t)\to0$ as $t\to\infty$.
  On the event $\{\mathcal{J}^{\lambda t}_0\leq\eta t\}$,
  we have $\mathcal{J}^{\lambda t}_0\leq ts$
  for all $s\in[s_1,s_2]$, and hence
  \beqlb\label{eqn.9014}
  A_2(t)
  \ar\leq\ar
  \mathbf P_x\Big(
  \sup_{0\leq u\leq s_2t}
  \Big|
  \frac{
  	X_{\mathcal{J}^{\lambda t}_0+u}
  	-
  	X_{\mathcal{J}^{\lambda t}_0}
  }{t}
  +
  \beta\cdot\frac{u}{t}
  \Big|
  \geq\frac{\varepsilon}{4},\,
  \mathcal{J}^{\lambda t}_0\leq\eta t
  \,\Big|\,
  \tau_0^-\in t+\varDelta_\delta
  \Big).
  \eeqlb
  Obviously, the last expression can be bounded by
  \beqnn
  \frac{1}{
  	\mathbf P_x\big(
  	\tau_0^-\in t+\varDelta_\delta
  	\big)
  }
  \cdot
  \mathbf E_x\Big[
  \sum_{0<r\leq\eta t\wedge\tau_0^-}
  \mathbf 1_{\{\Delta X_r>\lambda t\}}
  \cdot
  \mathbf 1_{\{
  	\sup_{0\leq u\leq s_2t}
  	|X_{r+u}-X_r+\beta u|
  	\geq\varepsilon t/4,\,
  	\tau_0^-\in t+\varDelta_\delta\}}
  \Big].
  \eeqnn
  By using the Markov property and then the compensation formula
  \eqref{eqn.CompenstationX}, the expectation in the preceding
  expression can further be bounded by
  \beqnn
  \lefteqn{	\int_0^{\eta t}dr	\int_{0}^{\infty}	\mathbf P_x\big(	X_{r-}\in dz,\tau_0^->r	\big)	\int_{\lambda t}^{\infty}	\nu(dw)\,	\mathbf P_{z+w}\Big(	\sup_{0\leq u\leq s_2t}	|X_u-z-w+\beta u|	\geq\frac{\varepsilon t}{4},\,	\tau_0^-\in t-r+\varDelta_\delta	\Big)}	\cr\cr
  \ar=\ar
  \int_0^{\eta t}dr
  \int_{0}^{\infty}
  \mathbf P_x\big(
  X_{r-}\in dz,\tau_0^->r
  \big)
  \int_{\lambda t}^{\infty}
  \nu(dw)\,
  \mathbf P\Big(
  \sup_{0\leq u\leq s_2t}
  |X_u+\beta u|
  \geq\frac{\varepsilon t}{4},\,
  \tau_{-(z+w)}^-\in t-r+\varDelta_\delta
  \Big),
  \eeqnn
  which, by Fubini's theorem, can be bounded by 
  \beqnn 
  \int_0^{\eta t}dr
  \int_{0}^{\infty}
  \mathbf P_x\big(
  X_{r-}\in dz,\tau_0^->r
  \big)
  \mathbf E\Big[
  \nu\Big(
  \big((
  -\underline X_{t-r}-z)\vee \lambda t,\,
  -\underline X_{t-r+\delta}-z
  \big];
  \sup_{0\leq u\leq s_2t}
  |X_u+\beta u|
  \geq\frac{\varepsilon t}{4}
  \Big].
  \eeqnn
  By Corollary 2.2 in \cite{WuXu2026b} and the monotonicity of $\overline\nu$,
  there exists a constant $C>0$ such that for all sufficiently
  large $t$, all $z\geq0$ and $0\leq r\leq\eta t$,
  \beqnn
  \nu\Big(
  \big((
  -\underline X_{t-r}-z)\vee \lambda t,\,
  -\underline X_{t-r+\delta}-z
  \big]
  \leq
  C\cdot\frac{\overline\nu(\lambda t)}{\lambda t}
  \cdot
  \big(
  \underline X_{t-r}-\underline X_{t-r+\delta}+2
  \big).
  \eeqnn
  Combining these two estimates together and then taking them back
  into \eqref{eqn.9014}, we obtain that $A_2(t)$ is bounded by
  \beqnn
  \frac{
  	C\cdot\overline\nu(\lambda t)/\lambda t
  }{
  	\mathbf P_x\big(
  	\tau_0^-\in t+\varDelta_\delta
  	\big)
  }
  \cdot
  \int_0^{\eta t}
  \mathbf P_x(\tau_0^->r)\,dr\cdot
  \sup_{0\leq r\leq\eta t}
  \mathbf E\Big[
  \underline X_{t-r}-\underline X_{t-r+\delta}+2;
  \sup_{0\leq u\leq s_2t}
  |X_u+\beta u|
  \geq\frac{\varepsilon t}{4}
  \Big].
  \eeqnn
  By Lemma~2.7 in \cite{WuXu2026b} and Proposition 4.2 in \cite{CuiWuXu2025},
  the last supremum is uniformly bounded in $t\geq0$ and also
  converges to $0$ as $t\to\infty$.
  By using Theorem~\ref{MainThm01}, \eqref{eqn.30}
  and then the dominated convergence theorem along with
  $\mathbf E_x[\tau_0^-]<\infty$, we obtain that
  $A_2(t)\to0$ as $t\to\infty$.
  Combining all preceding estimates together yields the limit
  \eqref{eqn.9013}.
  This completes the whole proof.
  \qed
  
  \subsection{Proofs of Lemma~\ref{Lemma.901},~\ref{Lemma.902},~\ref{Lemma.903} and~\ref{Lemma.904}}
  
  Let $\vartheta\in\boldsymbol{\varTheta}$ and, for $h,\delta,K>0$, define
  \beqnn
  \mathcal A_{\delta}(h,K)
  :=
  \Big\{
  0<s<\tau_0^-:
  \Delta X_s>\vartheta_h,\,
  X_{s-}\leq K,\,
  \sup_{s\leq r\leq\tau_0^-}X_r\in h+\varDelta_\delta,\,
  \sup_{s\leq r\leq\tau_0^-}X_r-X_s\leq K
  \Big\}.
  \eeqnn
  For $R>0$, write
  \beqnn
  \mathcal A_{\delta}^{\leq R}(h,K)
  :=
  \mathcal A_{\delta}(h,K)\cap(0,R]
  \quad\mbox{and}\quad
  \mathcal A_{\delta}^{>R}(h,K)
  :=
  \mathcal A_{\delta}(h,K)\cap(R,\infty).
  \eeqnn
  For any $s\in\mathcal A_{\delta}(h,K)$, we have
  \beqlb\label{eqn.P.height.jump-window}
  \Delta X_s\in(h-2K,h+\delta].
  \eeqlb
  For any interval $I\subset(0,\infty)$, the Markov property  and
  the compensation formula \eqref{eqn.CompenstationX} give
  \beqnn
  \mathbf E_x\big[
  \#\big(\mathcal A_{\delta}(h,K)\cap I\big)
  \big]
  \ar=\ar
  \int_I ds
  \int_{0}^{K}
  \mathbf P_x\big(
  X_{s-}\in dz,\tau_0^->s
  \big)
  G_\delta(h,K,z),
  \eeqnn
  where
  \beqlb\label{eqn.50123}
  G_\delta(h,K,z)
  :=
  \int_{\vartheta_h}^{\infty}
  \mathbf P_{z+w}\Big(
  \overline X_{\tau_0^-}\in h+\varDelta_\delta,\,
  \overline X_{\tau_0^-}-z-w\leq K
  \Big)\,\nu(dw).
  \eeqlb
  Since $\mathbf E_x[\tau_0^-]<\infty$,
  Proposition~4.10 in \cite{WuXu2026b} and dominated convergence
  yield, as $h\to\infty$,
  \beqlb\label{eqn.P.height.count}
  \frac{
  	\mathbf E_x\big[
  	\#\big(\mathcal A_{\delta}(h,K)\cap I\big)
  	\big]
  }{
  	\nu(h+\varDelta_\delta)
  }
  \to
  \mathbf P(\overline X_\infty\leq K)\cdot
  \int_I
  \mathbf P_x\big(
  X_{s-}\leq K,\tau_0^->s
  \big)\,ds.
  \eeqlb

  \noindent\textbf{\textit{Proof of Lemma~\ref{Lemma.901}.}}
  The convergence of the first absolute value  follows directly from Proposition~\ref{Lemma.MainThm04.bigjump}. 
  We now consider the second absolute value.
  By \eqref{eqnn1.9} and the triangle inequality,
  \beqnn
  \lefteqn{
  	\Big|
  	\mathbf P_x\big(
  	\mathcal{J}^{\lambda t}_0\leq R
  	\,\big|\,
  	\overline X_{\tau_0^-}\in\beta(t+\varDelta_\delta)
  	\big)
  	-
  	\mathbf P(\mathcal T_x\leq R)
  	\Big|
  }
  \cr\cr
  \ar\leq\ar
  \Big|
  \mathbf P_x\big(
  \mathcal{J}^{\lambda t}_0\leq R
  \,\big|\,
  \overline X_{\tau_0^-}\in\beta(t+\varDelta_\delta)
  \big)
  -
  \frac{
  	\mathbf P_x\big(
  	\mathcal{J}^{\lambda t}_0\leq R,\,
  	\#\mathcal A_{\beta\delta}^{\leq R}(\beta t,K)=1
  	\big)
  }{
  	\mathbf P_x\big(
  	\overline X_{\tau_0^-}\in\beta(t+\varDelta_\delta)
  	\big)
  }
  \Big|
  \cr\cr
  \ar\ar+
  \Big|
  \frac{
  	\mathbf P_x\big(
  	\mathcal{J}^{\lambda t}_0\leq R,\,
  	\#\mathcal A_{\beta\delta}^{\leq R}(\beta t,K)=1
  	\big)
  }{
  	\mathbf P_x\big(
  	\overline X_{\tau_0^-}\in\beta(t+\varDelta_\delta)
  	\big)
  }
  -
  \frac{
  	\mathbf E_x[R\wedge\tau_0^-]
  }{
  	\mathbf E_x[\tau_0^-]
  }
  \Big|.
  \eeqnn
  Thus, it suffices to prove that
  \beqlb\label{eqn.P.height.1}
  \lim_{K\to\infty}\limsup_{t\to\infty}
  \Big|
  \mathbf P_x\big(
  \mathcal{J}^{\lambda t}_0\leq R
  \,\big|\,
  \overline X_{\tau_0^-}\in\beta(t+\varDelta_\delta)
  \big)
  -
  \frac{
  	\mathbf P_x\big(
  	\mathcal{J}^{\lambda t}_0\leq R,\,
  	\#\mathcal A_{\beta\delta}^{\leq R}(\beta t,K)=1
  	\big)
  }{
  	\mathbf P_x\big(
  	\overline X_{\tau_0^-}\in\beta(t+\varDelta_\delta)
  	\big)
  }
  \Big|
  =0
  \eeqlb
  and
  \beqlb\label{eqn8.10}
  \lim_{K\to\infty}\limsup_{t\to\infty}
  \Big|
  \frac{
  	\mathbf P_x\big(
  	\mathcal{J}^{\lambda t}_0\leq R,\,
  	\#\mathcal A_{\beta\delta}^{\leq R}(\beta t,K)=1
  	\big)
  }{
  	\mathbf P_x\big(
  	\overline X_{\tau_0^-}\in\beta(t+\varDelta_\delta)
  	\big)
  }
  -
  \frac{
  	\mathbf E_x[R\wedge\tau_0^-]
  }{
  	\mathbf E_x[\tau_0^-]
  }
  \Big|
  =0.
  \eeqlb
  We prove the two limits separately.

  \textbf{(i)} 
  We first prove \eqref{eqn.P.height.1}.
  Since
  $
  \#\mathcal A_{\beta\delta}(\beta t,K)
  =
  \#\mathcal A_{\beta\delta}^{\leq R}(\beta t,K)
  +
  \#\mathcal A_{\beta\delta}^{>R}(\beta t,K),
  $
  we have
  \beqnn
  \Big|
  \mathbf 1_{\{
  	\overline X_{\tau_0^-}\in\beta(t+\varDelta_\delta)
  	\}}
  -
  \mathbf 1_{\{
  	\#\mathcal A_{\beta\delta}^{\leq R}(\beta t,K)=1
  	\}}
  \Big|
  \leq 
  \Big|
  \mathbf 1_{\{
  	\overline X_{\tau_0^-}\in\beta(t+\varDelta_\delta)
  	\}}
  -
  \mathbf 1_{\{
  	\#\mathcal A_{\beta\delta}(\beta t,K)=1
  	\}}
  \Big|
  +
  \mathbf 1_{\{
  	\#\mathcal A_{\beta\delta}^{>R}(\beta t,K)\geq1
  	\}}.
  \eeqnn
  Multiplying both sides by
  $\mathbf 1_{\{\mathcal{J}^{\lambda t}_0\leq R\}}$
  and taking expectations, we obtain that the absolute value in
  \eqref{eqn.P.height.1} is bounded by
  \beqlb\label{eqqn8.11}
  \frac{
  	\mathbf P_x\Big(
  	\big\{
  	\overline X_{\tau_0^-}\in\beta(t+\varDelta_\delta)
  	\big\}
  	\mathbin{\triangle}
  	\big\{
  	\#\mathcal A_{\beta\delta}(\beta t,K)=1
  	\big\}
  	\Big)
  }{
  	\mathbf P_x\big(
  	\overline X_{\tau_0^-}\in\beta(t+\varDelta_\delta)
  	\big)
  }
  +
  \frac{
  	\mathbf P_x\big(
  	\mathcal{J}^{\lambda t}_0\leq R,\,
  	\#\mathcal A_{\beta\delta}^{>R}(\beta t,K)\geq1
  	\big)
  }{
  	\mathbf P_x\big(
  	\overline X_{\tau_0^-}\in\beta(t+\varDelta_\delta)
  	\big)
  }.
  \eeqlb
  It remains to show that each term vanishes by letting first $t\to\infty$ and then $K\to\infty$.

  \noindent$\bullet$  For the first term in \eqref{eqqn8.11}, note that 
  \beqlb\label{eqnn8.11}
  \ar\ar \mathbf P_x\big(
  \big\{
  \overline X_{\tau_0^-}\in\beta(t+\varDelta_\delta)
  \big\}
  \mathbin{\triangle}
  \big\{
  \#\mathcal A_{\beta\delta}(\beta t,K)=1
  \big\}
  \big)\cr\cr
  \ar\ar\qquad\leq \Big|
  \mathbf P_x\big(
  \overline X_{\tau_0^-}\in\beta(t+\varDelta_\delta)
  \big)
  -
  \mathbf P_x\big(
  \#\mathcal A_{\beta\delta}(\beta t,K)=1
  \big)
  \Big|\cr\cr
  \ar\ar\qquad\qquad+2\mathbf P_x\big(
  \#\mathcal A_{\beta\delta}(\beta t,K)=1,\,
  \overline X_{\tau_0^-}\notin\beta(t+\varDelta_\delta)
  \big).
  \eeqlb
  For any $S>0$, applying the inequality
  $
  m+n-\mathbf 1_{\{m+n=1\}}
  \leq m(m-1)+2n
  $
  for $m,n\in\mathbb Z_+$ with
  $
  m=\#\mathcal A_{\beta\delta}^{\leq S}(\beta t,K)
  $,
  $
  n=\#\mathcal A_{\beta\delta}^{>S}(\beta t,K)
  $
  and then using \eqref{eqn.P.height.jump-window}, we obtain for all sufficiently large $t$,
  \beqnn
  \lefteqn{
  	0\leq
  	\mathbf E_x\big[
  	\#\mathcal A_{\beta\delta}(\beta t,K)
  	\big]
  	-
  	\mathbf P_x\big(
  	\#\mathcal A_{\beta\delta}(\beta t,K)=1
  	\big)
  }
  \cr\cr
  \ar\leq\ar
  \mathbf E_x\Big[
  \#\mathcal A_{\beta\delta}^{\leq S}(\beta t,K)
  \Big(
  \#\mathcal A_{\beta\delta}^{\leq S}(\beta t,K)-1
  \Big)
  \Big]
  +
  2\mathbf E_x\big[
  \#\mathcal A_{\beta\delta}^{>S}(\beta t,K)
  \big]
  \cr\cr
  \ar\leq\ar
  \mathbf E_x\Big[
  \#\mathcal B_{\leq S}^{\lambda t}
  \Big(
  \#\mathcal B_{\leq S}^{\lambda t}-1
  \Big)
  \Big]
  +
  2\mathbf E_x\big[
  \#\mathcal A_{\beta\delta}^{>S}(\beta t,K)
  \big].
  \eeqnn
  By Proposition~\ref{Prop.203}, we have 
  \beqnn
  \mathbf E_x\Big[
  \#\mathcal B_{\leq S}^{\lambda t}
  \Big(
  \#\mathcal B_{\leq S}^{\lambda t}-1
  \Big)
  \Big]
  =
  2\mathbf E_x\bigg[
  \sum_{0<s<r\leq S\wedge\tau_0^-}
  \mathbf 1_{\{
  	\Delta X_s>\lambda t,\,
  	\Delta X_r>\lambda t
  	\}}
  \bigg]
  \leq 
  2S\mathbf E_x[\tau_0^-]\,
  \big|\overline\nu(\lambda t)\big|^2,
  \eeqnn
  and hence by \eqref{eqn.30},
  \beqlb\label{eqnnn8.11}
  \limsup_{t\to\infty}
  \frac{
  	\mathbf E_x\Big[
  	\#\mathcal B_{\leq S}^{\lambda t}
  	\Big(
  	\#\mathcal B_{\leq S}^{\lambda t}-1
  	\Big)
  	\Big]
  }{
  	\nu(\beta t+\varDelta_{\beta\delta})
  }
  \ar\leq\ar
  2S\mathbf E_x[\tau_0^-]
  \limsup_{t\to\infty}
  \frac{
  	\big|\overline\nu(\lambda t)\big|^2
  }{
  	\nu(\beta t+\varDelta_{\beta\delta})
  }
  =
  0.
  \eeqlb
  Moreover, by \eqref{eqn.P.height.count} with $I=(S,\infty)$,
  \beqnn
  \lim_{t\to\infty}
  \frac{
  	\mathbf E_x\big[
  	\#\mathcal A_{\beta\delta}^{>S}(\beta t,K)
  	\big]
  }{
  	\nu(\beta t+\varDelta_{\beta\delta})
  }
  \ar=\ar
  \mathbf P(\overline X_\infty\leq K)
  \int_S^\infty
  \mathbf P_x\big(
  X_{s-}\leq K,\tau_0^->s
  \big)ds\leq \int_S^\infty
  \mathbf P_x\big(\tau_0^->s
  \big)ds,
  \eeqnn
  which goes to $0$ as $S\to\infty$,
  since $\mathbf E_x[\tau_0^-]<\infty$.
  Therefore,  as $t\to\infty$,  we have 
  \beqlb\label{eqn.P.height.multiple}
  \mathbf E_x\big[
  \#\mathcal A_{\beta\delta}(\beta t,K)
  \big]
  -
  \mathbf P_x\big(
  \#\mathcal A_{\beta\delta}(\beta t,K)=1
  \big)
  =
  o\big(
  \nu(\beta t+\varDelta_{\beta\delta})
  \big).
  \eeqlb
  By \eqref{eqn.P.height.multiple},  \eqref{eqn.P.height.count} with $I=(0,\infty)$ and Theorem~\ref{MainThm02},
  for any fixed $K>0$,
  \beqnn
  \lim_{t\to\infty}
  \frac{
  	\mathbf P_x\big(
  	\#\mathcal A_{\beta\delta}(\beta t,K)=1
  	\big)
  }{
  	\mathbf P_x\big(
  	\overline X_{\tau_0^-}\in\beta(t+\varDelta_\delta)
  	\big)
  }
  =
  \lim_{t\to\infty}
  \frac{	\mathbf E_x\big[
  	\#\mathcal A_{\beta\delta}(\beta t,K)
  	\big]
  }{	\mathbf E_x[\tau_0^-]  	\nu(\beta t+\varDelta_{\beta\delta})
  }
  =
  \frac{
  	\mathbf P(\overline X_\infty\leq K)
  }{
  	\mathbf E_x[\tau_0^-]
  }
  \int_0^\infty
  \mathbf P_x\big(
  X_{s-}\leq K,\tau_0^->s
  \big)ds,
  \eeqnn
  which goes to $1$ as $K\to\infty$, by  using monotone convergence. 
  Consequently,
  \beqlb\label{eqn8.11}
  \lim_{K\to\infty}\limsup_{t\to\infty}
  \frac{
  	\Big|
  	\mathbf P_x\big(
  	\overline X_{\tau_0^-}\in\beta(t+\varDelta_\delta)
  	\big)
  	-
  	\mathbf P_x\big(
  	\#\mathcal A_{\beta\delta}(\beta t,K)=1
  	\big)
  	\Big|
  }{
  	\mathbf P_x\big(
  	\overline X_{\tau_0^-}\in\beta(t+\varDelta_\delta)
  	\big)
  }
  =0.
  \eeqlb
  If $\#\mathcal A_{\beta\delta}(\beta t,K)=1$
  but $\overline X_{\tau_0^-}\notin\beta(t+\varDelta_\delta)$,
  then for the unique $s\in\mathcal A_{\beta\delta}(\beta t,K)$,
  $
  \sup_{0\leq r<s}X_r>\beta t+\beta\delta.
  $
  Consequently, by the compensation formula \eqref{eqn.CompenstationX} and the definition of
  $G_{\beta\delta}$ in \eqref{eqn.50123}, we obtain
  \beqlb\label{eqnn8.15}
  \lefteqn{
  	\frac{
  		1
  	}{
  		\nu(\beta t+\varDelta_{\beta\delta})
  	}\cdot
  	\mathbf P_x\big(
  	\#\mathcal A_{\beta\delta}(\beta t,K)=1,\,
  	\overline X_{\tau_0^-}
  	\notin\beta(t+\varDelta_\delta)
  	\big)
  }
  \cr\cr
  \ar\leq\ar
  \frac{1}{
  	\nu(\beta t+\varDelta_{\beta\delta})
  }\cdot 
  \mathbf E_x\bigg[
  \sum_{0<s<\tau_0^-}
  \mathbf 1_{\{
  	\sup_{0\leq r<s}X_r>\beta t+\beta\delta,\,
  	s\in\mathcal A_{\beta\delta}(\beta t,K)
  	\}}
  \bigg]
  \cr\cr
  \ar=\ar
  \int_0^\infty ds
  \int_{0}^{K}
  \mathbf P_x\big(
  \sup_{0\leq r<s}X_r>\beta t+\beta\delta,\,
  X_{s-}\in dz,\tau_0^->s
  \big)
  \frac{
  	G_{\beta\delta}(\beta t,K,z)
  }{
  	\nu(\beta t+\varDelta_{\beta\delta})
  }.
  \eeqlb
  By Proposition~4.10 in \cite{WuXu2026b}, we have for any fixed
  $K>0$ and all sufficiently large $t$, 
  \beqnn
  \sup_{z\in[0,K]}	\frac{G_{\beta\delta}(\beta t,K,z)}{\nu(\beta t+\varDelta_{\beta\delta})}
  \leq 1.
  \eeqnn
  Thus, the right-hand side of \eqref{eqnn8.15} is bounded by 
  \beqnn
  \int_0^\infty
  \mathbf P_x\big(
  \sup_{0\leq r<s}X_r>\beta t+\beta\delta,\,
  \tau_0^->s
  \big)ds
  \leq 
  \mathbf E_x\big[
  \tau_0^-;
  \overline X_{\tau_0^-}>\beta t+\beta\delta
  \big],
  \eeqnn
  which goes to $0$,  as  $t\to\infty$. Putting these estimates together along with  Theorem~\ref{MainThm02}, we obtain for any fixed $K>0$,
  \beqlb\label{eqn8.12}
  \lim_{t\to\infty}
  \frac{  	\mathbf P_x\big(  	\#\mathcal A_{\beta\delta}(\beta t,K)=1,\,  	\overline X_{\tau_0^-}	\notin\beta(t+\varDelta_\delta)  	\big)  }{  	\mathbf P_x\big(  	\overline X_{\tau_0^-}\in\beta(t+\varDelta_\delta)  	\big)  }  =0 .
  \eeqlb
  Combining \eqref{eqn8.11} and \eqref{eqn8.12} with \eqref{eqnn8.11} yields
  \beqlb\label{eqn8.18}
  \lim_{K\to\infty}\limsup_{t\to\infty}
  \frac{
  	\mathbf P_x\big(
  	\big\{
  	\overline X_{\tau_0^-}\in\beta(t+\varDelta_\delta)
  	\big\}
  	\mathbin{\triangle}
  	\big\{
  	\#\mathcal A_{\beta\delta}(\beta t,K)=1
  	\big\}
  	\big)
  }{
  	\mathbf P_x\big(
  	\overline X_{\tau_0^-}\in\beta(t+\varDelta_\delta)
  	\big)
  }
  =0.
  \eeqlb
  
  \noindent$\bullet$
  For the second term in \eqref{eqqn8.11}, for any $S>0$ and  all sufficiently large $t$, it holds that
  \beqnn
  \mathbf P_x\big(
  \mathcal{J}^{\lambda t}_0\leq R,\,
  \#\mathcal A_{\beta\delta}^{>R}(\beta t,K)\geq1
  \big)
  \ar\leq\ar
  \mathbf P_x\big(
  \mathcal{J}^{\lambda t}_0\leq R,\,
  \#\mathcal A_{\beta\delta}^{>R+S}(\beta t,K)\geq1
  \big)
  \cr\cr
  \ar\ar+
  \mathbf P_x\big(
  \mathcal{J}^{\lambda t}_0\leq R,\,
  \#\mathcal A_{\beta\delta}^{>R}(\beta t,K)
  -
  \#\mathcal A_{\beta\delta}^{>R+S}(\beta t,K)
  \geq1
  \big)
  \cr\cr
  \ar\leq\ar\mathbf E_x\big[
  \#\mathcal A_{\beta\delta}^{>R+S}(\beta t,K)
  \big]+
  \mathbf P_x\big(
  \#\mathcal B_{\leq R+S}^{\lambda t}\geq2
  \big),
  \eeqnn
  where the last inequality follows from Markov's inequality and
  \eqref{eqn.P.height.jump-window}, since
  $\mathcal{J}^{\lambda t}_0\leq R$ together with
  $\mathcal A_{\beta\delta}(\beta t,K)\cap(R,R+S]\neq\emptyset$
  implies that at least two jumps larger than $\lambda t$
  occur in $(0,R+S]$. By  \eqref{eqn.P.height.count} with
  $I=(R+S,\infty)$ and Theorem~\ref{MainThm02},  we have 
  \beqnn
  \limsup_{t\to\infty}
  \frac{
  	\mathbf E_x\big[
  	\#\mathcal A_{\beta\delta}^{>R+S}(\beta t,K)
  	\big]
  }{
  	\mathbf P_x\big(
  	\overline X_{\tau_0^-}\in\beta(t+\varDelta_\delta)
  	\big)
  }
  \ar=\ar
  \frac{
  	\mathbf P(\overline X_\infty\leq K)
  }{
  	\mathbf E_x[\tau_0^-]
  }
  \int_{R+S}^\infty
  \mathbf P_x\big(
  X_{s-}\leq K,\tau_0^->s
  \big)\,ds\cr
  \ar\leq\ar 
  C\cdot 
  \int_{R+S}^\infty
  \mathbf P_x(\tau_0^->s)\,ds,
  \eeqnn
  which goes to $0$ as $S\to\infty$.
  Moreover, applying \eqref{eqnnn8.11} with $S$ replaced by $R+S$
  and using Theorem~\ref{MainThm02}, we obtain
  \beqlb\label{eqn8.16}
  \limsup_{t\to\infty}
  \frac{
  	\mathbf P_x\big(
  	\#\mathcal B_{\leq R+S}^{\lambda t}\geq2
  	\big)
  }{
  	\mathbf P_x\big(
  	\overline X_{\tau_0^-}\in\beta(t+\varDelta_\delta)
  	\big)
  }
  \ar\leq\ar
  \limsup_{t\to\infty}
  \frac{
  	\mathbf E_x\Big[
  	\#\mathcal B_{\leq R+S}^{\lambda t}
  	\Big(
  	\#\mathcal B_{\leq R+S}^{\lambda t}-1
  	\Big)
  	\Big]
  }{
  	\mathbf P_x\big(
  	\overline X_{\tau_0^-}\in\beta(t+\varDelta_\delta)
  	\big)
  }
  =0.
  \eeqlb
  Combining these estimates and letting first $t\to\infty$
  and then $S\to\infty$, we obtain, for any fixed $K>0$,
  \beqlb\label{eqnn8.20}
  \limsup_{t\to\infty}
  \frac{
  	\mathbf P_x\big(
  	\mathcal{J}^{\lambda t}_0\leq R,\,
  	\#\mathcal A_{\beta\delta}^{>R}(\beta t,K)\geq1
  	\big)
  }{
  	\mathbf P_x\big(
  	\overline X_{\tau_0^-}\in\beta(t+\varDelta_\delta)
  	\big)
  }
  =0.
  \eeqlb

  \textbf{(ii)} 
  We now prove \eqref{eqn8.10}. By \eqref{eqn.P.height.jump-window}, for all
  sufficiently large $t$, we have 
  \beqnn
  \mathbf P_x\big(
  \mathcal{J}^{\lambda t}_0\leq R,\,
  \#\mathcal A_{\beta\delta}^{\leq R}(\beta t,K)=1
  \big)
  \ar=\ar
  \mathbf P_x\big(
  \#\mathcal A_{\beta\delta}^{\leq R}(\beta t,K)=1
  \big).
  \eeqnn
  Since
  $
  \#\mathcal A_{\beta\delta}^{\leq R}(\beta t,K)
  \leq
  \#\mathcal A_{\beta\delta}(\beta t,K)
  $
  and $m-\mathbf 1_{\{m=1\}}$ is non-decreasing on $\mathbb Z_+$, we have
  \beqnn
  \Big|\mathbf E_x\big[
  \#\mathcal A_{\beta\delta}^{\leq R}(\beta t,K)
  \big]
  -
  \mathbf P_x\big(
  \#\mathcal A_{\beta\delta}^{\leq R}(\beta t,K)=1
  \big)\Big|
  \leq 
  \mathbf E_x\big[
  \#\mathcal A_{\beta\delta}(\beta t,K)
  \big]
  -
  \mathbf P_x\big(
  \#\mathcal A_{\beta\delta}(\beta t,K)=1
  \big).
  \eeqnn
  Hence, by \eqref{eqn.P.height.multiple}, as $t\to\infty$,
  \beqnn
  \mathbf E_x\big[
  \#\mathcal A_{\beta\delta}^{\leq R}(\beta t,K)
  \big]
  -
  \mathbf P_x\big(
  \#\mathcal A_{\beta\delta}^{\leq R}(\beta t,K)=1
  \big)
  =
  o\big(
  \nu(\beta t+\varDelta_{\beta\delta})
  \big).
  \eeqnn
  Combining this with \eqref{eqn.P.height.count} for $I=(0,R]$ and
  Theorem~\ref{MainThm02}, for any fixed $K>0$,
  \beqnn
  \lim_{t\to\infty}
  \frac{
  	\mathbf P_x\big(
  	\#\mathcal A_{\beta\delta}^{\leq R}(\beta t,K)=1
  	\big)
  }{
  	\mathbf P_x\big(
  	\overline X_{\tau_0^-}\in\beta(t+\varDelta_\delta)
  	\big)
  }
  \ar=\ar
  \lim_{t\to\infty}
  \frac{	\mathbf E_x\big[
  	\#\mathcal A^{\leq R}_{\beta\delta}(\beta t,K)
  	\big]
  }{	\mathbf E_x[\tau_0^-]  	\nu(\beta t+\varDelta_{\beta\delta})
  }\cr
  \ar=\ar
  \frac{
  	\mathbf P(\overline X_\infty\leq K)
  }{
  	\mathbf E_x[\tau_0^-]
  }
  \int_0^R
  \mathbf P_x\big(
  X_{s-}\leq K,\tau_0^->s
  \big)ds.
  \eeqnn
  Letting $K\to\infty$ and using monotone convergence, we obtain
  \beqnn
  \lim_{K\to\infty}\lim_{t\to\infty}
  \frac{
  	\mathbf P_x\big(
  	\mathcal{J}^{\lambda t}_0\leq R,\,
  	\#\mathcal A_{\beta\delta}^{\leq R}(\beta t,K)=1
  	\big)
  }{
  	\mathbf P_x\big(
  	\overline X_{\tau_0^-}\in\beta(t+\varDelta_\delta)
  	\big)
  }
  \ar=\ar
  \frac{1}{
  	\mathbf E_x[\tau_0^-]
  }
  \int_0^R
  \mathbf P_x(\tau_0^->s)ds
  =
  \frac{
  	\mathbf E_x[R\wedge\tau_0^-]
  }{
  	\mathbf E_x[\tau_0^-]
  }.
  \eeqnn
  Thus \eqref{eqn8.10} follows.
  \qed

  \noindent\textbf{\textit{Proof of Lemma~\ref{Lemma.902}.}}
  For any constants $R,K>0$, the law of total probability gives
  \beqlb\label{eqnn8.17}
  \lefteqn{
  	\mathbf P_x\bigg(
  	\Big|\frac{\Delta  X_{\mathcal{J}^{\lambda t}_0}}{t}-\beta\Big|
  	\geq\varepsilon
  	\,\bigg|\,
  	\tau_0^-\in t+\varDelta_\delta
  	\bigg)
  }
  \cr\cr
  \ar\leq\ar
  \mathbf P_x\bigg(
  \Big|\frac{\Delta  X_{\mathcal{J}^{\lambda t}_0}}{t}-\beta\Big|
  \geq\varepsilon,\,
  \mathcal{J}^{\lambda t}_0>R
  \,\bigg|\,
  \tau_0^-\in t+\varDelta_\delta
  \bigg)
  \cr\cr
  \ar\ar+
  \mathbf P_x\bigg(
  \Big|\frac{\Delta  X_{\mathcal{J}^{\lambda t}_0}}{t}-\beta\Big|
  \geq\varepsilon,\,
  \mathcal{J}^{\lambda t}_0\leq R,\,
  \sup_{0\leq r<\mathcal{J}^{\lambda t}_0}X_r>K
  \,\bigg|\,
  \tau_0^-\in t+\varDelta_\delta
  \bigg)
  \cr\cr
  \ar\ar+
  \mathbf P_x\bigg(
  \Big|\frac{\Delta  X_{\mathcal{J}^{\lambda t}_0}}{t}-\beta\Big|
  \geq\varepsilon,\,
  \mathcal{J}^{\lambda t}_0\leq R,\,
  \sup_{0\leq r<\mathcal{J}^{\lambda t}_0}X_r\leq K
  \,\bigg|\,
  \tau_0^-\in t+\varDelta_\delta
  \bigg).
  \eeqlb
  The first term can be bounded by
  $
  \mathbf P_x(
  \mathcal{J}^{\lambda t}_0>R
  \mid
  \tau_0^-\in t+\varDelta_\delta)
  $,
  which  tends to $0$ as $t\to\infty$
  and then $R\to\infty$, by Lemma~\ref{Lemma.901}.  
  For the second term, it is easy to see that  
  \beqlb\label{eqn8.20}
  \lefteqn{
  	\mathbf P_x\big(
  	\mathcal{J}^{\lambda t}_0\leq R,\,
  	\sup_{0\leq r<\mathcal{J}^{\lambda t}_0}X_r>K
  	\,\big|\,
  	\tau_0^-\in t+\varDelta_\delta
  	\big)}
  \cr\cr
  \ar\leq\ar
  \frac{1}{
  	\mathbf P_x\big(
  	\tau_0^-\in t+\varDelta_\delta
  	\big)}
  \mathbf E_x\bigg[
  \sum_{0<s\leq R\wedge\tau_0^-}
  \mathbf 1_{\{
  	\sup_{0\leq r<s}X_r>K,\,
  	\Delta X_s>\lambda t,\,
  	\tau_0^-\in t+\varDelta_\delta
  	\}}
  \bigg].
  \eeqlb
  Applying  the  strong Markov property at the
  jump time and then  compensation formula
  \eqref{eqn.CompenstationX}, the expectation in \eqref{eqn8.20} equals  
  \beqnn
  \lefteqn{\mathbf E_x\bigg[
  	\sum_{0<s\leq R\wedge\tau_0^-}
  	\mathbf 1_{\{
  		\sup_{0\leq r<s}X_r>K,\,
  		\Delta X_s>\lambda t
  		\}}
  	\mathbf P_{X_{s-}+\Delta X_s}\big(
  	\tau_0^-\in t-s+\varDelta_\delta
  	\big)
  	\bigg]}
  \cr\cr
  \ar=\ar
  \int_0^R ds
  \int_{0}^{\infty}
  \mathbf P_x\big(
  \sup_{0\leq r<s}X_r>K,\,
  X_{s-}\in dz,\,
  \tau_0^->s
  \big)
  \int_{\lambda t}^{\infty}
  \mathbf P_{z+w}\big(
  \tau_0^-\in t-s+\varDelta_\delta
  \big)\nu(dw),
  \eeqnn
  which, by Proposition~6.6 in \cite{WuXu2026b},  can be bounded by 
  \beqnn
  C \cdot \nu\big(\beta t+\varDelta_\delta \big)\cdot  \int_0^R
  \mathbf P_x\big(
  \sup_{0\leq r<s}X_r>K,\,
  \tau_0^->s
  \big)ds\leq C\cdot   \nu\big(\beta t+\varDelta_\delta \big)\cdot R\cdot \mathbf P_x\big(
  \overline X_R>K
  \big).
  \eeqnn
  Combining the preceding estimates gives 
  \beqlb\label{eqn.P.902.pre-truncation}
  \mathbf P_x\big(
  \mathcal{J}^{\lambda t}_0\leq R,\,
  \sup_{0\leq r<\mathcal{J}^{\lambda t}_0}X_r>K
  \,\big|\,
  \tau_0^-\in t+\varDelta_\delta
  \big)
  \leq 
  C\cdot
  \frac{
  	\nu\big(\beta t+\varDelta_\delta\big)
  }{
  	\mathbf P_x\big(
  	\tau_0^-\in t+\varDelta_\delta
  	\big)
  }\cdot 
  R\cdot \mathbf P_x\big(
  \overline X_R>K
  \big).
  \eeqlb
  By Theorem~\ref{MainThm01}, the ratio on the right-hand side is
  bounded for all sufficiently large $t$. Hence, letting first
  $t\to\infty$ and then $K\to\infty$, the second term in \eqref{eqnn8.17} tends to $0$. Finally, consider the third term. By the definition of
  $\mathcal{J}^{\lambda t}_0$ and the same argument as in the proof of
  \eqref{eqn.P.902.pre-truncation}, we obtain
  \beqnn
  \lefteqn{
  	\mathbf P_x\bigg(
  	\Big|\frac{\Delta  X_{\mathcal{J}^{\lambda t}_0}}{t}-\beta\Big|
  	\geq\varepsilon,\,
  	\mathcal{J}^{\lambda t}_0\leq R,\,
  	\sup_{0\leq r<\mathcal{J}^{\lambda t}_0}X_r\leq K
  	\,\bigg|\,
  	\tau_0^-\in t+\varDelta_\delta
  	\bigg)}
  \cr\cr
  \ar\leq\ar
  \frac{1}{
  	\mathbf P_x\big(
  	\tau_0^-\in t+\varDelta_\delta
  	\big)}
  \int_0^R ds
  \int_{0}^{K}
  \mathbf P_x\big(
  X_{s-}\in dz,\tau_0^->s
  \big)
  \int_{\lambda t}^{\infty}
  \mathbf 1_{\{|w/t-\beta|\geq\varepsilon\}}
  \mathbf P_{z+w}\big(
  \tau_0^-\in t-s+\varDelta_\delta
  \big)\nu(dw)
  \cr\cr
  \ar\leq\ar 
  \frac{
  	C	\nu(\beta t+\varDelta_\delta)
  }{
  	\mathbf P_x\big(
  	\tau_0^-\in t+\varDelta_\delta
  	\big)
  }\cdot 
  \sup_{s\leq R}	\sup_{0\leq z\leq K}
  \int_{\lambda t}^{\infty}
  \mathbf 1_{\{|w/t-\beta|\geq\varepsilon\}}
  \mathbf P_{z+w}\big(
  \tau_0^-\in t-s+\varDelta_\delta
  \big)\frac{\nu(dw)}{\nu(\beta t+\varDelta_\delta)}.
  \eeqnn
  It remains to show that the second factor tends to $0$. This will imply
  that the third term in \eqref{eqnn8.17} vanishes as $t\to\infty$. By spatial homogeneity
  and Fubini's theorem, it follows that 
  \beqlb\label{eqn8.22}
  \lefteqn{
  	\int_{\lambda t}^{\infty}
  	\mathbf 1_{\{|w/t-\beta|\geq\varepsilon\}}
  	\mathbf P_{z+w}\big(
  	\tau_0^-\in t-s+\varDelta_\delta
  	\big)\nu(dw)}
  \cr
  \ar=\ar
  \mathbf E\Big[
  \nu\Big(
  \big[
  -\underline X_{t-s}-z,\,
  -\underline X_{t-s+\delta}-z
  \big)
  \cap(\lambda t,\infty)
  \cap
  \big\{
  w:|w/t-\beta|\geq\varepsilon
  \big\}
  \Big)
  \Big].
  \eeqlb
  For $s\leq R$ and $0\leq z\leq K$, if
  $
  w\in
  \big[
  -\underline X_{t-s}-z,\,
  -\underline X_{t-s+\delta}-z
  \big)$ and 
  $\big|w/t-\beta\big|\geq\varepsilon,
  $
  then for all sufficiently large $t$,
  \beqnn
  \Big|
  \frac{-\underline X_{t-s}}{t}-\beta
  \Big|
  >
  \frac{\varepsilon}{4}
  \quad\mbox{or}\quad
  \underline X_{t-s}-\underline X_{t-s+\delta}
  >
  \frac{\varepsilon t}{4}.
  \eeqnn
  Moreover, the interval
  $
  \big[
  -\underline X_{t-s}-z,\,
  -\underline X_{t-s+\delta}-z
  \big)
  \cap(\lambda t,\infty)
  $
  has left endpoint at least $\lambda t$ and length at most $\underline X_{t-s}-\underline X_{t-s+\delta}$. 
  Therefore, by Corollary~2.2 in \cite{WuXu2026b}, the expectation in \eqref{eqn8.22} can be bounded, uniformly for all sufficiently large $t$, by 
  \beqnn
  \ar\ar C\cdot\frac{\overline\nu(\lambda t)}{\lambda t}\cdot  \mathbf E\big[ \underline X_{t-s}-\underline X_{t-s+\delta}+1; \big|-\underline X_{t-s}/t-\beta \big|>\varepsilon/4 \big]\cr
  \ar\ar+  C\cdot\frac{\overline\nu(\lambda t)}{\lambda t}\cdot \mathbf E\big[ \underline X_{t-s}-\underline X_{t-s+\delta}+1; \underline X_{t-s}-\underline X_{t-s+\delta} >\varepsilon t/4 \big] .
  \eeqnn
  Consequently, by \eqref{eqn.30}, we have 
  \beqnn
  \lefteqn{ \sup_{s\leq R}	\sup_{0\leq z\leq K}\int_{\lambda t}^{\infty}\mathbf 1_{\{|w/t-\beta|\geq\varepsilon\}}\mathbf P_{z+w}\big(\tau_0^-\in t-s+\varDelta_\delta\big)\frac{\nu(dw)}{\nu(\beta t+\varDelta_\delta)}}\cr\cr
  \ar\leq\ar C\cdot\sup_{s\leq R}\mathbf E\Big[\underline X_{t-s}-\underline X_{t-s+\delta}+1;\Big|\frac{-\underline X_{t-s}}{t}-\beta\Big|>\frac{\varepsilon}{4}\Big]\cr\cr
  \ar\ar+ C\cdot \sup_{s\leq R}\mathbf E\Big[\underline X_{t-s}-\underline X_{t-s+\delta}+1;\underline X_{t-s}-\underline X_{t-s+\delta}>\frac{\varepsilon t}{4}\Big].
  \eeqnn
  By Lemma~2.6 and 2.7 in \cite{WuXu2026b}, both terms on the
  right-hand side go to $0$, as $t\to\infty$.
  We now consider the second term in Lemma~\ref{Lemma.902}. For any $R>0$,
  the law of total probability also gives
  \beqlb\label{eqn8.25}
  \lefteqn{
  	\mathbf P_x\bigg(
  	\Big|\frac{\Delta  X_{\mathcal{J}^{\lambda t}_0}}{t}-\beta\Big|
  	\geq\varepsilon
  	\,\bigg|\,
  	\overline X_{\tau_0^-}\in\beta(t+\varDelta_\delta)
  	\bigg)
  }
  \cr\cr
  \ar=\ar
  \mathbf P_x\bigg(
  \Big|\frac{\Delta  X_{\mathcal{J}^{\lambda t}_0}}{t}-\beta\Big|
  \geq\varepsilon,\,
  \mathcal{J}^{\lambda t}_0>R
  \,\bigg|\,
  \overline X_{\tau_0^-}\in\beta(t+\varDelta_\delta)
  \bigg)
  \cr\cr
  \ar\ar+
  \mathbf P_x\bigg(
  \Big|\frac{\Delta  X_{\mathcal{J}^{\lambda t}_0}}{t}-\beta\Big|
  \geq\varepsilon,\,
  \mathcal{J}^{\lambda t}_0\leq R
  \,\bigg|\,
  \overline X_{\tau_0^-}\in\beta(t+\varDelta_\delta)
  \bigg).
  \eeqlb
  The first term on the right-hand side can be bounded by
  $
  \mathbf P_x(
  \mathcal{J}^{\lambda t}_0>R
  \mid
  \overline X_{\tau_0^-}\in\beta(t+\varDelta_\delta))
  $,
  which goes to $0$ as $t\to\infty$
  and then $R\to\infty$ by Lemma~\ref{Lemma.901}. 
  For the second term, the same estimate as the absolute value in
  \eqref{eqn.P.height.1}  shows that
  \beqnn
  \lefteqn{
  	\bigg|
  	\mathbf P_x\Big(
  	\Big|\frac{\Delta  X_{\mathcal{J}^{\lambda t}_0}}{t}-\beta\Big|
  	\geq\varepsilon,\,
  	\mathcal{J}^{\lambda t}_0\leq R
  	\,\big|\,
  	\overline X_{\tau_0^-}\in\beta(t+\varDelta_\delta)
  	\Big)}
  \cr\cr
  \ar\ar-
  \frac{
  	\mathbf P_x\big(
  	\big| \Delta  X_{\mathcal{J}^{\lambda t}_0}/t-\beta\big|
  	\geq\varepsilon,\,
  	\mathcal{J}^{\lambda t}_0\leq R,\,
  	\#\mathcal A_{\beta\delta}^{\leq R}(\beta t,K)=1
  	\big)
  }{
  	\mathbf P_x\big(
  	\overline X_{\tau_0^-}\in\beta(t+\varDelta_\delta)
  	\big)
  }\bigg|
  \cr\cr
  \ar\leq\ar
  \frac{
  	\mathbf P_x\big(
  	\big\{
  	\overline X_{\tau_0^-}\in\beta(t+\varDelta_\delta)
  	\big\}
  	\mathbin{\triangle}
  	\big\{
  	\#\mathcal A_{\beta\delta}(\beta t,K)=1
  	\big\}
  	\big)
  }{
  	\mathbf P_x\big(
  	\overline X_{\tau_0^-}\in\beta(t+\varDelta_\delta)
  	\big)
  }+
  \frac{
  	\mathbf P_x\big(
  	\mathcal{J}^{\lambda t}_0\leq R,\,
  	\#\mathcal A_{\beta\delta}^{>R}(\beta t,K)\geq1
  	\big)
  }{
  	\mathbf P_x\big(
  	\overline X_{\tau_0^-}\in\beta(t+\varDelta_\delta)
  	\big)
  }.
  \eeqnn
  Here $A\triangle B$ denotes the symmetric difference of the two sets $A$ and $B$. 
  By \eqref{eqn8.18} and \eqref{eqnn8.20}, the right-hand side tends to
  $0$ by letting first $t\to\infty$ and then $K\to\infty$.
  Hence, it remains to prove that
  \beqlb\label{eeqn.P.902.height}
  \lim_{K\to\infty}\limsup_{t\to\infty}
  \frac{
  	\mathbf P_x\big(
  	\big|\Delta X_{\mathcal{J}^{\lambda t}_0}/t-\beta\big|
  	\geq\varepsilon,\,
  	\mathcal{J}^{\lambda t}_0\leq R,\,
  	\#\mathcal A_{\beta\delta}^{\leq R}(\beta t,K)=1
  	\big)
  }{
  	\mathbf P_x\big(
  	\overline X_{\tau_0^-}\in\beta(t+\varDelta_\delta)
  	\big)
  }
  =0.
  \eeqlb
  Splitting according to whether
  $\mathcal{J}^{\lambda t}_0\in
  \mathcal A_{\beta\delta}^{\leq R}(\beta t,K)$, the fraction on the
  left-hand side is bounded by
  \beqlb\label{eqn8.17}
  \frac{
  	\mathbf P_x\big(
  	\big|\Delta  X_{\mathcal{J}^{\lambda t}_0}/t-\beta\big|
  	\geq\varepsilon,\,
  	\mathcal{J}^{\lambda t}_0
  	\in
  	\mathcal A_{\beta\delta}^{\leq R}(\beta t,K)
  	\big)
  }{
  	\mathbf P_x\big(
  	\overline X_{\tau_0^-}\in\beta(t+\varDelta_\delta)
  	\big)
  }+
  \frac{
  	\mathbf P_x\big(
  	\mathcal{J}^{\lambda t}_0
  	\notin
  	\mathcal A_{\beta\delta}^{\leq R}(\beta t,K),\,
  	\#\mathcal A_{\beta\delta}^{\leq R}(\beta t,K)=1
  	\big)
  }{
  	\mathbf P_x\big(
  	\overline X_{\tau_0^-}\in\beta(t+\varDelta_\delta)
  	\big)
  }.
  \eeqlb
  For all
  sufficiently large $t$, on the event
  $
  \mathcal{J}^{\lambda t}_0
  \in
  \mathcal A_{\beta\delta}^{\leq R}(\beta t,K)
  $,
  we have
  $
  \Delta X_{\mathcal{J}^{\lambda t}_0}-\beta t
  \in
  (-2K,\beta\delta]
  $
  and hence the first term in \eqref{eqn8.17} equals $0$. Moreover, on the event
  $
  \big\{
  \mathcal{J}^{\lambda t}_0
  \notin
  \mathcal A_{\beta\delta}^{\leq R}(\beta t,K),\,
  \#\mathcal A_{\beta\delta}^{\leq R}(\beta t,K)=1
  \big\},
  $ we have that at least two jumps larger than $\lambda t$ occur in the time interval $[0,R]$. This induces that
  \beqnn
  \frac{
  	\mathbf P_x\big(
  	\mathcal{J}^{\lambda t}_0
  	\notin
  	\mathcal A_{\beta\delta}^{\leq R}(\beta t,K),\,
  	\#\mathcal A_{\beta\delta}^{\leq R}(\beta t,K)=1
  	\big)
  }{
  	\mathbf P_x\big(
  	\overline X_{\tau_0^-}\in\beta(t+\varDelta_\delta)
  	\big)
  }
  \leq 
  \frac{
  	\mathbf P_x\big(
  	\#\mathcal B_{\leq R}^{\lambda t}\geq2
  	\big)
  }{
  	\mathbf P_x\big(
  	\overline X_{\tau_0^-}\in\beta(t+\varDelta_\delta)
  	\big)
  },
  \eeqnn
  which goes to $0$ as $t\to\infty$ by \eqref{eqn8.16}. These two estimates together yield the limit \eqref{eeqn.P.902.height}. This completes the proof.
  \qed

  \noindent\textbf{\textit{Proof of Lemma~\ref{Lemma.903}.}}
  For any constant $R>0$, the law of total probability shows that the first
  term in Lemma~\ref{Lemma.903} can be bounded by
  \beqnn
  \mathbf P_x\Big(
  \sup_{s\leq T}
  X_{(s\wedge\mathcal{J}^{\lambda t}_0)-}
  \geq\varepsilon t,\,
  \mathcal{J}^{\lambda t}_0\leq R
  \,\Big|\,
  \tau_0^-\in t+\varDelta_\delta
  \Big)
  +
  \mathbf P_x\big(
  \mathcal{J}^{\lambda t}_0>R
  \,\big|\,
  \tau_0^-\in t+\varDelta_\delta
  \big).
  \eeqnn
  For any $K>0$ and all sufficiently large $t$ such that
  $\varepsilon t>K$, we have
  \beqnn
  \Big\{
  \sup_{s\leq T}
  X_{(s\wedge\mathcal{J}^{\lambda t}_0)-}
  \geq\varepsilon t,\,
  \mathcal{J}^{\lambda t}_0\leq R
  \Big\}
  \subset
  \Big\{
  \mathcal{J}^{\lambda t}_0\leq R,\,
  \sup_{0\leq r<\mathcal{J}^{\lambda t}_0}X_r>K
  \Big\}.
  \eeqnn
  Hence, by \eqref{eqn.P.902.pre-truncation}, it follows that
  \beqnn
  \limsup_{t\to\infty}
  \mathbf P_x\Big(
  \sup_{s\leq T}
  X_{(s\wedge\mathcal{J}^{\lambda t}_0)-}
  \geq\varepsilon t,\,
  \mathcal{J}^{\lambda t}_0\leq R
  \,\Big|\,
  \tau_0^-\in t+\varDelta_\delta
  \Big)
  \leq
  C\cdot R \cdot \mathbf P_x\big(\overline X_R>K\big).
  \eeqnn
  Letting $K\to\infty$, the right-hand side vanishes. 
  The second one goes to $0$ as $t\to\infty$ and then $R\to \infty$ by Lemma~\ref{Lemma.901}. Consequently, the first term in
  Lemma~\ref{Lemma.903} converges to $0$ as $t\to\infty$.

  We now consider the second term in Lemma~\ref{Lemma.903}.
  For any $R>0$, the law of total probability also shows that it can
  be bounded by
  \beqnn
  \mathbf P_x\Big(
  \sup_{s\leq T\wedge\tau_0^-}
  X_{(s\wedge\mathcal{J}^{\lambda t}_0)-}
  \geq\varepsilon t,\,
  \mathcal{J}^{\lambda t}_0\leq R
  \,\Big|\,
  \overline X_{\tau_0^-}\in\beta(t+\varDelta_\delta)
  \Big)
  +
  \mathbf P_x\Big(
  \mathcal{J}^{\lambda t}_0>R
  \,\Big|\,
  \overline X_{\tau_0^-}\in\beta(t+\varDelta_\delta)
  \Big).
  \eeqnn
  The second term goes to $0$ as $t\to\infty$ and then $R\to\infty$
  by Lemma~\ref{Lemma.901}. It remains to prove that the first term also tends to $0$ as
  $t\to\infty$. By an argument analogous to the second term in
  \eqref{eqn8.25},
  it suffices to prove that, for any fixed $R,K>0$,
  \beqlb\label{eqn8.29}
  \lim_{t\to\infty}\frac{
  	\mathbf P_x\Big(
  	\sup_{s\leq T\wedge\tau_0^-}
  	X_{(s\wedge\mathcal{J}^{\lambda t}_0)-}
  	\geq\varepsilon t,\,
  	\mathcal{J}^{\lambda t}_0\in
  	\mathcal A_{\beta\delta}^{\leq R}(\beta t,K)
  	\Big)
  }{
  	\mathbf P_x\big(
  	\overline X_{\tau_0^-}\in\beta(t+\varDelta_\delta)
  	\big)
  }
  =0.
  \eeqlb
  By the definition of
  $\mathcal A_{\beta\delta}^{\leq R}(\beta t,K)$, the strong Markov
  property and the compensation formula \eqref{eqn.CompenstationX},
  the last fraction is bounded by
  \beqnn
  \frac{
  	\nu\big(\beta (t+\varDelta_{\delta})\big)
  }{
  	\mathbf P_x\big(
  	\overline X_{\tau_0^-}\in\beta(t+\varDelta_\delta)
  	\big)
  }\cdot 
  \int_0^R ds\int_{0}^{K}
  \mathbf P_x\Big(
  \sup_{0\leq r\leq T\wedge s}X_{r-}
  \geq\varepsilon t,\,
  X_{s-}\in dz,\,
  \tau_0^->s
  \Big)
  \frac{
  	G_{\beta\delta}(\beta t,K,z)
  }{
  	\nu\big(\beta (t+\varDelta_{\delta})\big)
  }.
  \eeqnn
  By Theorem~\ref{MainThm02} and Proposition~4.10 in
  \cite{WuXu2026b}, for all sufficiently large $t$, this is bounded by
  \beqnn
  C\cdot
  \int_0^R
  \mathbf P_x\Big(
  \sup_{0\leq r\leq T\wedge s}X_{r-}
  \geq\varepsilon t,\,
  \tau_0^->s
  \Big)ds
  \ar\leq\ar
  C\cdot R \cdot \mathbf P_x\big(
  \overline X_R\geq\varepsilon t
  \big),
  \eeqnn
  which goes to $0$, as $t\to\infty$.
  Hence \eqref{eqn8.29} follows. This completes the proof.
  \qed

  \noindent\textbf{\textit{Proof of Lemma~\ref{Lemma.904}.}}
  On the event
  $
  \big\{
  \sup_{s\leq T}
  |X_{s\vee\mathcal{J}^{\lambda t}_0}
  -X_{\mathcal{J}^{\lambda t}_0}|
  \geq\varepsilon t
  \big\},
  $
  we have $\mathcal{J}^{\lambda t}_0\leq T$. Hence, for $t>T$,
  \beqnn
  \mathbf P_x\Big(
  \sup_{s\leq T}
  |X_{s\vee\mathcal{J}^{\lambda t}_0}
  -X_{\mathcal{J}^{\lambda t}_0}|
  \geq\varepsilon t,\,
  \tau_0^-\in t+\varDelta_\delta
  \Big)\leq 
  \mathbf E_x\Big[
  \sum_{0<s\leq T}
  \mathbf 1_{\{
  	\Delta X_s>\lambda t,\,
  	\sup_{s\leq r\leq T}|X_r-X_s|
  	\geq\varepsilon t,\,
  	\tau_0^-\in t+\varDelta_\delta
  	\}}
  \Big].
  \eeqnn
  By the strong Markov property and the compensation formula
  \eqref{eqn.CompenstationX}, the right-hand side equals
  \beqnn
  \int_0^T ds\int_{0}^{\infty}
  \mathbf P_x\big(
  X_{s-}\in dz,\tau_0^->s
  \big)
  \int_{\lambda t}^{\infty}
  \mathbf P_{z+w}\Big(
  \sup_{0\leq r\leq T-s}
  |X_r-z-w|
  \geq\varepsilon t,\,
  \tau_0^-\in t-s+\varDelta_\delta
  \Big)\nu(dw),
  \eeqnn
  which can be bounded by
  \beqnn
  \mathbf E_x[\tau_0^-]\cdot
  \sup_{s\in[0,T]}\sup_{z\geq0}
  \int_{\lambda t}^{\infty}
  \mathbf P_{z+w}\Big(
  \sup_{0\leq r\leq T-s}
  |X_r-z-w|
  \geq\varepsilon t,\,
  \tau_0^-\in t-s+\varDelta_\delta
  \Big)\nu(dw).
  \eeqnn
  By the spatial homogeneity of $X$ and Fubini's theorem,
  \beqnn
  \lefteqn{
  	\sup_{s\in[0,T]}\sup_{z\geq0}
  	\int_{\lambda t}^{\infty}
  	\mathbf P_{z+w}\Big(
  	\sup_{0\leq r\leq T-s}
  	|X_r-z-w|
  	\geq\varepsilon t,\,
  	\tau_0^-\in t-s+\varDelta_\delta
  	\Big)\nu(dw)}
  \cr\cr
  \ar=\ar
  \sup_{s\in[0,T]}\sup_{z\geq0}
  \int_{\lambda t}^{\infty}
  \mathbf P\Big(
  \sup_{0\leq r\leq T-s}|X_r|
  \geq\varepsilon t,\,
  \tau_{-z-w}^-\in t-s+\varDelta_\delta
  \Big)\nu(dw)
  \cr\cr
  \ar\leq\ar
  \sup_{s\in[0,T]}\sup_{z\geq0}
  \mathbf E\Big[
  \nu\Big(
  \big(
  (-\underline X_{t-s}-z)\vee \lambda t,\,
  -\underline X_{t-s+\delta}-z
  \big]
  \Big);
  \sup_{0\leq r\leq T}|X_r|
  \geq\varepsilon t
  \Big].
  \eeqnn
  By Corollary~2.2 in \cite{WuXu2026b} and \eqref{eqn.30},
  the last expression is bounded, for all sufficiently large $t$, by
  \beqnn
  \lefteqn{
  	C\cdot
  	\sup_{s\in[0,T]}
  	\mathbf E\Big[
  	\underline X_{t-s}-\underline X_{t-s+\delta}+1;
  	\sup_{0\leq r\leq T}|X_r|
  	\geq\varepsilon t
  	\Big]\cdot \nu\big(\beta (t+\varDelta_\delta)\big)}
  \cr\cr
  \ar\leq\ar
  C\cdot
  \mathbf E\Big[
  \underline X_{t-T}-\underline X_{t+\delta}+1;
  \sup_{0\leq r\leq T}|X_r|
  \geq\varepsilon t
  \Big]\cdot \nu\big(\beta( t+\varDelta_\delta)\big).
  \eeqnn
  Combining these estimates and using
  Theorem~\ref{MainThm01},
  we obtain, for all sufficiently large $t$,
  \beqnn
  \mathbf P_x\Big(
  \sup_{s\leq T}
  t^{-1}|X_{s\vee\mathcal{J}^{\lambda t}_0}
  -X_{\mathcal{J}^{\lambda t}_0}|
  \geq\varepsilon
  \,\Big|\,
  \tau_0^-\in t+\varDelta_\delta
  \Big)\leq 
  C\cdot\mathbf E\Big[
  \underline X_{t-T}-\underline X_{t+\delta}+1;
  \sup_{0\leq r\leq T}|X_r|
  \geq\varepsilon t
  \Big],
  \eeqnn
  which goes to $0$ as $t\to\infty$, since
  $\sup_{0\leq r\leq T}|X_r|<\infty$ a.s. and the family
  $\{\underline X_{t-T}-\underline X_{t+\delta}\}_{t\geq T}$
  is uniformly integrable; see Lemma~2.7 in \cite{WuXu2026b}. 
  This proves the convergence of the first term in
  Lemma~\ref{Lemma.904} to $0$. 
  
  We now consider the second term in Lemma~\ref{Lemma.904}.
  Similarly as in the proof of Lemma~\ref{Lemma.903}, it suffices to prove
  that, for any fixed $K>0$, as $t\to\infty$,
  \beqlb\label{eqn.402}
  \frac{
  	\mathbf P_x\Big(
  	\sup_{s\leq T\wedge\tau_0^-}
  	|X_{s\vee\mathcal{J}^{\lambda t}_0}
  	-X_{\mathcal{J}^{\lambda t}_0}|
  	\geq\varepsilon t,\,
  	\#\mathcal A_{\beta\delta}^{\leq T}(\beta t,K)=1,\,
  	\mathcal{J}^{\lambda t}_0\in
  	\mathcal A_{\beta\delta}^{\leq T}(\beta t,K)
  	\Big)
  }{
  	\mathbf P_x\big(
  	\overline X_{\tau_0^-}\in\beta(t+\varDelta_\delta)
  	\big)
  }
  \to0.
  \eeqlb
  In view of \eqref{eqn.P.height.jump-window}, on the event
  $\{
  \#\mathcal A_{\beta\delta}^{\leq T}(\beta t,K)=1,\,
  \mathcal{J}^{\lambda t}_0\in
  \mathcal A_{\beta\delta}^{\leq T}(\beta t,K)
  \}$,
  we have $\Delta X_{\mathcal{J}^{\lambda t}_0} \in(\beta t-2K,\beta t+\beta\delta]$ and the numerator of the preceding fraction can be bounded by
  \beqnn
  \mathbf P_x\bigg( \sup_{s\leq T\wedge\tau_0^-} |X_{s\vee\mathcal{J}^{\lambda t}_0}-X_{\mathcal{J}^{\lambda t}_0}|\geq\varepsilon t,\,\Delta X_{\mathcal{J}^{\lambda t}_0}\in(\beta t-2K,\beta t+\beta\delta],\,\mathcal{J}^{\lambda t}_0\leq T\bigg)
  \eeqnn
  which further can be bounded by 
  \beqnn
  \lefteqn{\mathbf E_x\Big[
  	\sum_{0<s\leq T\wedge\tau_0^-}
  	\mathbf 1_{\{
  		\Delta X_s\in(\beta t-2K,\beta t+\beta\delta],\,
  		\sup_{s\leq r\leq T}|X_r-X_s|\geq\varepsilon t
  		\}}
  	\Big]}\cr\cr
  \ar\leq\ar 	\nu\big((\beta t-2K,\beta t+\beta\delta]\big) \cdot \int_0^T\mathbf P_x(\tau_0^->s)\,ds
  \cdot
  \mathbf P\Big(
  \sup_{0\leq r\leq T}|X_r|
  \geq\varepsilon t
  \Big)
  \eeqnn
  In conclusion, the left-hand side of  \eqref{eqn.402} can be bounded by  
  \beqnn
  \frac{	\nu\big((\beta t-2K,\beta t+\beta\delta]\big)}{	\mathbf P_x\big(	\overline X_{\tau_0^-}\in\beta(t+\varDelta_\delta)	\big)}\cdot\int_0^T\mathbf P_x(\tau_0^->s)\,ds\cdot\mathbf P\Big( \sup_{0\leq r\leq T}|X_r|\geq\varepsilon t\Big),
  \eeqnn
  which goes to $0$ as $t\to\infty$ by
  Theorem~\ref{MainThm02}, \eqref{eqn.30} and $\sup_{0\leq r\leq T}|X_r|<\infty$ a.s.
  \qed

    \section{Proof of Theorem \ref{MainThm06}} 
   \label{Sec.ProofThm06}
   \setcounter{equation}{0} 
   
   Fix $x>0$. For notational convenience in the proof of Theorem~\ref{MainThm06}, we introduce the following family of finite measures on $\mathbb R_+$:
   \beqlb\label{eqn.8011}
   \mu_{t,\delta}^{x}(dz) :=  \int_{t}^{t+\delta} \frac{  \mathbf P_x\big(X_s\in dz,\tau_0^->s\big)}{ \widehat V(x)\nu\big( \beta t+\varDelta_\delta\big) } \, ds,\quad t,\delta>0 .
   \eeqlb 
   The vague convergence of these measures as $t\to\infty$ is established in the following proposition.
   
   \begin{proposition}\label{Prop.801}
   	For any constants $\delta>0$, $0\leq a<b<\infty$ and any continuous function $f$ on $\mathbb{R}_+$, we have as $t\to\infty$,
   	\beqnn
   	\bigg|\int_a^b f(z)\,\mu_{t,\delta}^{x}(dz)
   	-\int_a^b f(z)V(z)\,dz\bigg|\to0.
   	\eeqnn
   \end{proposition}
   \proof
   By \eqref{eqn.8011} and the change of variable $r=s-t$, we have
   \beqnn
   \lefteqn{\bigg|\int_a^b f(z)\,\mu_{t,\delta}^{x}(dz)
   	-\int_a^b f(z)V(z)\,dz\bigg|}
   \cr\cr
   \ar=\ar
   \bigg|\int_0^\delta dr\int_a^b f(z)
   \frac{\mathbf P_x\big(X_{t+r}\in dz,\tau_0^->t+r\big)}
   {\widehat V(x)\nu\big(\beta t+\varDelta_\delta\big)}
   -\int_a^b f(z)V(z)\,dz\bigg|\cr\cr
   \ar\leq\ar
   \frac{1}{\delta}\int_0^\delta
   \bigg|\int_a^b f(z)
   \frac{\delta\cdot\mathbf P_x\big(X_{t+r}\in dz,\tau_0^->t+r\big)}
   {\widehat V(x)\nu\big(\beta t+\varDelta_\delta\big)}
   -\int_a^b f(z)V(z)\,dz\bigg|\,dr\cr\cr
   \ar\leq\ar
   \sup_{r\in[0,\delta]}
   \bigg|\int_a^b f(z)
   \frac{\delta\cdot\mathbf P_x\big(X_{t+r}\in dz,\tau_0^->t+r\big)}
   {\widehat V(x)\nu\big(\beta t+\varDelta_\delta\big)}
   -\int_a^b f(z)V(z)\,dz\bigg|.
   \eeqnn
   By the triangle inequality, the last supremum can be bounded by the sum of the following two terms:
   \beqnn
   \sup_{r\in[0,\delta]}
   \bigg|
   \frac{\alpha\cdot\overline\nu\big(\beta(t+r)\big)\cdot\delta}
   {\beta(t+r)\cdot\nu\big(\beta t+\varDelta_\delta\big)}
   -1\bigg|
   \cdot\int_a^b |f(z)|V(z)\,dz
   \eeqnn
   and
   \beqnn
   \ar\ar
   \sup_{r\in[0,\delta]}
   \frac{\alpha\cdot\overline\nu\big(\beta(t+r)\big)\cdot\delta}
   {\beta(t+r)\cdot\nu\big(\beta t+\varDelta_\delta\big)}
   \bigg|
   \int_a^b f(z)\,\frac{\mathbf P_x\big(X_{t+r}\in dz,\tau_0^->t+r\big)}
   {\alpha\cdot\widehat V(x)\cdot \overline\nu\big(\beta(t+r)\big)/\big(\beta(t+r)\big)}
   -\int_a^b f(z)V(z)\,dz
   \bigg|.
   \eeqnn
   By \eqref{eqn.30}, the first term goes to $0$ as $t\to\infty$.
   Moreover, by Corollary~1.6 in \cite{CuiWuXu2026b},
   \beqnn
   \lim_{t\to\infty}\sup_{r\geq0}
   \bigg|
   \int_a^b f(z)\,\frac{\mathbf P_x\big(X_{t+r}\in dz,\tau_0^->t+r\big)}
   {\alpha\cdot\widehat V(x)\cdot \overline\nu\big(\beta(t+r)\big)/\big(\beta(t+r)\big)}
   -\int_a^b f(z)V(z)\,dz
   \bigg|=0,
   \eeqnn
   and then the second term also vanishes as $t\to\infty$.
   This completes the proof.
   \qed

   \textbf{\textit{Proof of Theorem \ref{MainThm06}.}}
   By the identity \eqref{eqn1.9}, we first have for all $\delta\in [\delta_0,\infty]$, 
   \beqnn
   \mathbf E_x\big(f(0,0)\,\big|\,\tau_0^-\in t+\varDelta_\delta\big)
   =
   f(0,0)
   =
   \frac{1}{\mathbf{E}\big[\widehat{H}_1\big]}
   \Big(
   \underline{\tt a}\cdot f(0,0)
   +
   \int_0^\infty du\int_0^\infty
   f(0,0)V(u)\,\nu(-u-dv)
   \Big).
   \eeqnn
   Consequently, it suffices to prove that the uniform convergence \eqref{eqn.MainThm06} holds for any non-negative, bounded and
   continuous function $f$ on $[0,\infty)^2$ with
   $\|f\|_\infty\leq1$ and $f(0,0)=0$, with $f(u,v)=0$ if $u<0$ or $v<0$;
   that is, as $t\to\infty$,
   \beqlb\label{eqn.801}
   \bigg|
   \mathbf E_x\Big(
   f\big(X_{\tau_0^--},-X_{\tau_0^-}\big)
   \,\Big|\,\tau_0^-\in t+\varDelta_\delta
   \Big)
   -
   \frac{1}{\mathbf{E}[\widehat{H}_1]}
   \int_0^\infty du\int_0^\infty
   f(u,v)V(u)\,\nu(-u-dv)
   \bigg|
   \to0.
   \eeqlb
   The proof consists of two steps.
   We first establish that the limit \eqref{eqn.801} holds for any
   fixed $\delta>0$ and then extend this convergence uniformly to all
   $\delta\in[\delta_0,\infty]$ by the method of interval partition.
   
   \textbf{Step 1.}
   We now prove that the limit \eqref{eqn.801} holds for any fixed $\delta>0$.
   By the triangle inequality, the left-hand side of \eqref{eqn.801}
   can be bounded by the sum of the following two terms:
   \beqnn
   \int_0^\infty du\int_0^\infty
   f(u,v)V(u)\,\nu(-u-dv)
   \cdot
   \bigg|
   \frac{\widehat V(x)\nu\big(\beta t+\varDelta_\delta\big)}
   {\mathbf P_x\big(\tau_0^-\in t+\varDelta_\delta\big)}
   -
   \frac{1}{\mathbf{E}\big[\widehat{H}_1\big]}
   \bigg|
   \eeqnn
   and
   \beqnn
   \frac{\widehat V(x)\nu\big(\beta t+\varDelta_\delta\big)}
   {\mathbf P_x\big(\tau_0^-\in t+\varDelta_\delta\big)}
   \cdot
   \bigg|
   \frac{
   	\mathbf E_x\big(
   	f(X_{\tau_0^--},-X_{\tau_0^-}),
   	\tau_0^-\in t+\varDelta_\delta
   	\big)}
   {\widehat V(x)\nu\big(\beta t+\varDelta_\delta\big)}
   -
   \int_0^\infty du\int_0^\infty
   f(u,v)V(u)\,\nu(-u-dv)
   \bigg|.
   \eeqnn
   By the fact that $\|f\|_\infty\leq1$ and \eqref{eqn1.9},
   the first term can be bounded by
   \beqnn
   \mathbf{E}\big[\widehat{H}_1\big]
   \cdot
   \bigg|
   \frac{\widehat V(x)\nu\big(\beta t+\varDelta_\delta\big)}
   {\mathbf P_x\big(\tau_0^-\in t+\varDelta_\delta\big)}
   -
   \frac{1}{\mathbf{E}\big[\widehat{H}_1\big]}
   \bigg|,
   \eeqnn
   which goes to $0$ as $t\to\infty$; see
   Theorem~\ref{MainThm01}, \eqref{eqn.EH1} and \eqref{Asytau}.
   To obtain the limit \eqref{eqn.801}, it remains to prove that
   as $t\to\infty$,
   \beqlb\label{eqn.802}
   \bigg|
   \frac{
   	\mathbf E_x\big(
   	f(X_{\tau_0^--},-X_{\tau_0^-}),
   	\tau_0^-\in t+\varDelta_\delta
   	\big)}
   {\widehat V(x)\nu\big(\beta t+\varDelta_\delta\big)}
   -
   \int_0^\infty du\int_0^\infty
   f(u,v)V(u)\,\nu(-u-dv)
   \bigg|
   \to0.
   \eeqlb
   For $t\geq1$, by using the Markov property at time $1$ we get
   \beqnn
   \lefteqn{\mathbf E_x\big(
   	f(X_{\tau_0^--},-X_{\tau_0^-}),
   	\tau_0^-\in t+\varDelta_\delta
   	\big) }\cr\cr
   \ar=\ar
   \int_{0}^{\infty}
   \mathbf{E}_z\Big[
   f\big(X_{\tau_0^--},-X_{\tau_0^-}\big),
   \tau_0^-\in t-1+\varDelta_\delta
   \Big]\,\mathbf P_x\big(X_1\in dz,\tau_0^->1\big).
   \eeqnn
   For any $z>0$, it is easy to identify that
   \beqnn
   f\big(X_{\tau_0^--},-X_{\tau_0^-}\big)
   \cdot
   \mathbf{1}_{\{\tau_0^-\in t-1+\varDelta_\delta\}}
   \ar=\ar
   \sum_{s\in t-1+\varDelta_\delta}
   f\big(X_{s-},-X_s\big)
   \cdot
   \mathbf{1}_{\{\underline X_{s-}>0,X_s<0\}}
   \cr\cr
   \ar=\ar
   \sum_{s\in t-1+\varDelta_\delta}
   f\big(X_{s-},-X_{s-}-\Delta X_s\big)
   \cdot
   \mathbf{1}_{\{
   	\underline X_{s-}>0,\Delta X_s<-X_{s-}
   	\}},
   \eeqnn
   which along with the compensation formula
   \eqref{eqn.CompenstationX} yields that
   \beqnn
   \lefteqn{\mathbf{E}_z\Big[
   	f\big(X_{\tau_0^--},-X_{\tau_0^-}\big),
   	\tau_0^-\in t-1+\varDelta_\delta
   	\Big]}\cr\cr
   \ar=\ar
   \mathbf{E}_z\bigg[
   \sum_{s\in t-1+\varDelta_\delta}
   f\big(X_{s-},-X_{s-}-\Delta X_s\big)
   \cdot
   \mathbf{1}_{\{
   	\underline X_{s-}>0,\Delta X_s<-X_{s-}
   	\}}
   \bigg]
   \cr\cr
   \ar=\ar
   \int_{t-1}^{t-1+\delta}
   \mathbf{E}_z\Big[
   \int_{X_{s-}}^\infty
   f(X_{s-},-X_{s-}+w)\,\nu(-dw),
   \underline X_{s-}>0
   \Big]\,ds.
   \eeqnn
   Recall that $f(0,0)=0$ and $f(u,v)=0$ if $u<0$ or $v<0$.
   For convenience, we define
   \beqlb\label{eqn.8001}
   g(u):=
   \int_{(u,\infty)} f(u,-u+w)\,\nu(-dw),
   \quad u\geq0.
   \eeqlb
   Using Fubini's theorem and then the Markov property,
   \beqnn
   \mathbf E_x\big[
   f(X_{\tau_0^--},-X_{\tau_0^-}),
   \tau_0^-\in t+\varDelta_\delta
   \big] 
   \ar=\ar
   \int_{0}^{\infty}\mathbf P_x\big(X_1\in dz,\tau_0^->1\big)
   \int_{t-1}^{t-1+\delta}
   \mathbf{E}_z\big[
   g(X_{s-}),\tau_0^->s
   \big]\,ds
   \cr\cr
   \ar=\ar
   \int_{t-1}^{t-1+\delta}
   \mathbf E_x\Big(
   \mathbf{E}_{X_1}\big[
   g(X_{s-}),\tau_0^->s
   \big],\tau_0^->1
   \Big)\,ds
   \cr\cr
   \ar=\ar
   \int_{t-1}^{t-1+\delta}
   \mathbf E_x\big(
   g(X_{(s+1)-}),\tau_0^->s+1
   \big)\,ds
   \cr\cr
   \ar=\ar
   \int_{t-1}^{t-1+\delta}ds
   \int_0^\infty g(u)\,\mathbf P_x\big(X_{s+1}\in du,\tau_0^->s+1\big) \cr\cr
   \ar=\ar 
   \int_0^\infty g(u)
   \int_t^{t+\delta}\mathbf P_x\big(X_s\in du,\tau_0^->s\big)\,ds
   \eeqnn
   and
   \beqnn
   \int_0^\infty V(u)\cdot g(u)\,du
   =
   \int_0^\infty du\int_0^\infty
   f(u,v)V(u)\,\nu(-u-dv).
   \eeqnn
   Recall the measure $\mu_{t,\delta}^{x}(dz)$ defined in \eqref{eqn.8011}.
   The left-hand side of \eqref{eqn.802} equals 
   \beqnn
   \bigg|
   \int_0^\infty g(u)\,\mu_{t,\delta}^{x}(du)
   -
   \int_0^\infty g(u)V(u)\,du
   \bigg|.
   \eeqnn
   For any $\varepsilon\in(0,1)$ and integer $M\geq1$,
   it can be bounded by the sum of the following three terms:
   \beqnn
   I_1(t,\varepsilon)
   \ar:=\ar
   \bigg|
   \int_0^\varepsilon g(u)\,\mu_{t,\delta}^{x}(du)
   -
   \int_0^\varepsilon g(u)V(u)\,du
   \bigg|,
   \cr\cr
   I_2(t,\varepsilon,M)
   \ar:=\ar
   \bigg|
   \int_\varepsilon^M g(u)\,\mu_{t,\delta}^{x}(du)
   -
   \int_\varepsilon^M g(u)V(u)\,du
   \bigg|,
   \cr\cr
   I_3(t,M)
   \ar:=\ar
   \bigg|
   \int_M^\infty g(u)\,\mu_{t,\delta}^{x}(du)
   -
   \int_M^\infty g(u)V(u)\,du
   \bigg|.
   \eeqnn
   Firstly, note that $g$ is bounded on $[\varepsilon,M]$,
   Proposition~\ref{Prop.801} yields that 
   $I_2(t,\varepsilon,M)\to0$ as $t\to\infty$.
   To obtain the limit \eqref{eqn.802} and then the limit
   \eqref{eqn.MainThm06} for any fixed $\delta>0$, it remains to prove
   \beqlb\label{eqn.804}
   \limsup_{\varepsilon\to0+}
   \limsup_{t\to\infty}I_1(t,\varepsilon)=0
   \quad\mbox{and}\quad
   \limsup_{M\to\infty}
   \limsup_{t\to\infty}I_3(t,M)=0.
   \eeqlb
   
   \textbf{(i)}
   Firstly, by the triangle inequality we have
   \beqlb\label{eqn.80111}
   I_1(t,\varepsilon)
   \leq
   \int_0^\varepsilon |g(u)|\,\mu_{t,\delta}^{x}(du)
   +
   \int_0^\varepsilon |g(u)|V(u)\,du.
   \eeqlb
   Note that $g(u)\leq\nu\big((-\infty,-u)\big)$,
   by \eqref{eqn1.9} the second term on the right-hand side goes to
   $0$ as $\varepsilon\to0+$.
   For the first term, in view of \eqref{eqn.8001} we decompose $g$
   into the following two parts: for any $\varepsilon'\in(0,1)$,
   \beqnn
   g_1(\varepsilon',u):=
   \int_u^{u+\varepsilon'} f(u,-u+w)\,\nu(-dw)
   \quad\mbox{and}\quad
   g_2(\varepsilon',u):=
   \int_{u+\varepsilon'}^\infty f(u,-u+w)\,\nu(-dw).
   \eeqnn
   This along with the triangle inequality yields that the first term
   on the right-hand side of \eqref{eqn.80111} can be bounded by
   $I_{1,1}(t,\varepsilon,\varepsilon')
   +I_{1,2}(t,\varepsilon,\varepsilon')$ with
   \beqnn
   I_{1,i}(t,\varepsilon,\varepsilon')
   \ar:=\ar
   \int_0^\varepsilon
   g_i(\varepsilon',u)\,\mu_{t,\delta}^{x}(du),
   \quad i=1,2.
   \eeqnn
   For $i=1$, since
   $g_1(\varepsilon',u)
   \leq\sup_{z\in[0,\varepsilon']}|f(u,z)|
   \cdot\nu\big((-\infty,-u)\big)$, we have
   \beqnn
   I_{1,1}(t,\varepsilon,\varepsilon')
   \ar\leq\ar
   \sup_{v\in[0,\varepsilon]}
   \sup_{z\in[0,\varepsilon']}|f(v,z)|
   \cdot
   \int_0^\varepsilon
   \nu\big((-\infty,-u)\big)\,\mu_{t,\delta}^{x}(du).
   \eeqnn
   It follows from  \eqref{eqn.8011} and the compensation formula
   \eqref{eqn.CompenstationX} applied to terminal jumps that
   \beqnn
   \int_0^\varepsilon
   \nu\big((-\infty,-u)\big)\,\mu_{t,\delta}^{x}(du)
   \ar\leq\ar
   \int_0^\infty
   \nu\big((-\infty,-u)\big)\,\mu_{t,\delta}^{x}(du)
   \cr\cr
   \ar=\ar
   \frac{1}{\widehat V(x)\nu\big(\beta t+\varDelta_\delta\big)}
   \int_t^{t+\delta}ds\int_0^\infty
   \nu\big((-\infty,-u)\big)\,\mathbf P_x\big(X_s\in du,\tau_0^->s\big)
   \cr\cr
   \ar=\ar
   \frac{1}{\widehat V(x)\nu\big(\beta t+\varDelta_\delta\big)}
   \mathbf E_x\bigg[
   \sum_{0<s\leq\tau_0^-}
   \mathbf{1}_{\{
   	s\in t+\varDelta_\delta,\,
   	X_{s-}>0,\,
   	\Delta X_s<-X_{s-}
   	\}}
   \bigg]
   \cr\cr
   \ar\leq\ar
   \frac{\mathbf P_x\big(\tau_0^-\in t+\varDelta_\delta\big)}
   {\widehat V(x)\nu\big(\beta t+\varDelta_\delta\big)}.
   \eeqnn
   Hence, by  Theorem~\ref{MainThm01}, \eqref{Asytau} and \eqref{eqn.30},
   we have as $t\to\infty$,
   \beqnn
   I_{1,1}(t,\varepsilon,\varepsilon')
   \ar\leq\ar
   \sup_{v\in[0,\varepsilon]}
   \sup_{z\in[0,\varepsilon']}|f(v,z)|
   \cdot
   \frac{\mathbf P_x\big(\tau_0^-\in t+\varDelta_\delta\big)}
   {\widehat V(x)\nu\big(\beta t+\varDelta_\delta\big)}
   \to
   \beta\underline n(\zeta)\cdot
   \sup_{v\in[0,\varepsilon]}
   \sup_{z\in[0,\varepsilon']}|f(v,z)|,
   \eeqnn
   which goes to $0$ as $\varepsilon\to0$ and then
   $\varepsilon'\to0$ because of the continuity of $f$ and $f(0,0)=0$.
   For $i=2$, we have
   $g_2(\varepsilon',u)
   \leq\nu\big((-\infty,-u-\varepsilon']\big)
   \leq\nu\big((-\infty,-\varepsilon']\big)$
   for all $u\geq0$ and then
   \beqnn
   I_{1,2}(t,\varepsilon,\varepsilon')
   \ar\leq\ar
   \nu\big((-\infty,-\varepsilon']\big)
   \cdot
   \mu_{t,\delta}^{x}\big((0,\varepsilon]\big).
   \eeqnn
   By Proposition~\ref{Prop.801} with $a=0$ and $b=\varepsilon$,
   we have
   \beqnn
   \limsup_{t\to\infty}I_{1,2}(t,\varepsilon,\varepsilon')
   \ar\leq\ar
   \nu\big((-\infty,-\varepsilon']\big)
   \cdot
   \int_0^\varepsilon V(z)\,dz,
   \eeqnn
   which goes to $0$ as $\varepsilon\to0$.
   These two estimates together yield the first limit in
   \eqref{eqn.804} by letting first $t\to\infty$,
   then $\varepsilon\to0+$ and finally $\varepsilon'\to0+$.
   
   \textbf{(ii)} 
   Since  $\|f\|_\infty\leq 1$ and $V(\infty)<\infty$, we also have 
   \beqlb\label{eqn.803}
   \ar\ar I_3(t,M)
   \leq   \int_{t}^{t+\delta} ds \int_M^\infty \nu\big((-\infty,-u)\big) \frac{\mathbf P_x\big(X_s\in du,\tau_0^->s\big) }{\widehat V(x)\nu\big( \beta t+\varDelta_\delta\big)} + \int_M^\infty\nu\big((-\infty,-u)\big)\,V(u)\,du  \cr
   \ar\ar \leq   \sum_{k=M}^{\infty} \nu\big((-\infty,-k)\big) \int_{t}^{t+\delta}
   \frac{\mathbf P_x\big(X_s\in k+\varDelta_1,\tau_0^->s\big)}{\widehat V(x)\nu\big( \beta t+\varDelta_\delta\big)}\, ds + V(\infty)\cdot\int_M^\infty  \nu\big((-\infty,-u)\big) \,du .
   \eeqlb
   By Corollary 1.7 in \cite{CuiWuXu2026b} and $V(\infty)<\infty$, there exists $C>0$ such that, for any sufficiently large $t$, uniformly in $s\in[t,t+\delta]$ and $k\geq M$,
   \beqnn
   \mathbf P_x\big(X_s\in k+\varDelta_1,\tau_0^->s\big)
   \ar\leq\ar
   C\widehat V(x)\cdot\frac{\overline\nu(\beta s+k)}{\beta s+k}
   \int_k^{k+1}V(z)\,dz
   \leq
   C\widehat V(x)\cdot\frac{\overline\nu(\beta s)}{\beta s}.
   \eeqnn
   By using the monotonicity of $\overline{\nu}$ and then  \eqref{eqn.30} with   fixed $\delta>0$, we obtain for all sufficiently large $t$, 
   \beqnn
   \lefteqn{\sum_{k=M}^{\infty}\nu\big((-\infty,-k)\big) \int_{t}^{t+\delta}
   	\frac{\mathbf P_x\big(X_s\in k+\varDelta_1,\tau_0^->s\big)}{\widehat V(x)\nu\big( \beta t+\varDelta_\delta\big)}\, ds}\cr\cr
   \ar\leq\ar \frac{C}{\nu\big( \beta t+\varDelta_\delta\big)} \sum_{k=M}^{\infty} \nu\big((-\infty,-k)\big) \int_{t}^{t+\delta} \frac{\overline\nu(\beta s)}{\beta s}\, ds \cr
   \ar\leq\ar C\delta\cdot \frac{\overline\nu(\beta t)/(\beta t)}{\nu\big( \beta t+\varDelta_\delta\big)} \cdot \sum_{k=M}^{\infty} \nu\big((-\infty,-k)\big)  
   \leq C   \int_{M-1}^{\infty} \nu\big((-\infty,-z)\big)\, dz. 
   \eeqnn 
   Plugging this back into the right-hand side of the second inequality in \eqref{eqn.803}, 
   \beqnn
   \limsup_{t\to\infty} I_3(t,M) \ar\leq \ar  C\cdot  	\int_{M-1}^\infty  \nu\big((-\infty,-u)\big)\,du \to 0,
   \eeqnn 
   as $M\to\infty$, since  $ \int_1^\infty  \nu\big((-\infty,-v)\big)\, dv =    \int_{-\infty}^{-1} (|v|-1) \,\nu(dv) <\infty$.

   \textbf{Step 2.} We now extend \eqref{eqn.MainThm06} to hold uniformly for $\delta\in[\delta_0,\infty]$.
   For any large integer $n>1/\delta_0$ and $\delta\in[\delta_0,\infty]$, since $f\geq0$, we have
   \beqnn
   \lefteqn{\sum_{k=0}^{[n\delta]-1}
   	\mathbf E_x\big(f(X_{\tau_0^--},-X_{\tau_0^-}),
   	\tau_0^-\in t+k/n+\varDelta_{1/n}\big)}\cr\cr
   \ar\leq\ar
   \mathbf E_x\big(f(X_{\tau_0^--},-X_{\tau_0^-}),
   \tau_0^-\in t+\varDelta_\delta\big)
   \leq
   \sum_{k=0}^{[n\delta]}
   \mathbf E_x\big(f(X_{\tau_0^--},-X_{\tau_0^-}),
   \tau_0^-\in t+k/n+\varDelta_{1/n}\big).
   \eeqnn
   For convenience, we write $A$ for the limit on the right-hand side of \eqref{eqn.MainThm06} with $f(0,0)=0$.
   Applying the limit \eqref{eqn.MainThm06} with $\delta=1/n$, we obtain a positive function $\rho_n$ on $\mathbb R_+$
   that decreases to $0$ at infinity and satisfies, for all sufficiently large $t$ and all $k\geq0$,
   \beqnn
   \big(A-\rho_n(t)\big)\cdot
   \mathbf P_x\big(\tau_0^-\in t+k/n+\varDelta_{1/n}\big)
   \ar\leq\ar
   \mathbf E_x\big(f(X_{\tau_0^--},-X_{\tau_0^-}),
   \tau_0^-\in t+k/n+\varDelta_{1/n}\big)
   \cr\cr
   \ar\leq\ar
   \big(A+\rho_n(t)\big)\cdot
   \mathbf P_x\big(\tau_0^-\in t+k/n+\varDelta_{1/n}\big).
   \eeqnn
   Combining these two estimates together yields
   \beqnn
   \mathbf E_x\Big(f(X_{\tau_0^--},-X_{\tau_0^-})
   \,\Big|\,\tau_0^-\in t+\varDelta_\delta\Big)
   \ar\leq\ar
   \frac{A+\rho_n(t)}
   {\mathbf P_x\big(\tau_0^-\in t+\varDelta_\delta\big)}
   \sum_{k=0}^{[n\delta]}
   \mathbf P_x\big(\tau_0^-\in t+k/n+\varDelta_{1/n}\big)
   \cr\cr
   \ar\leq\ar
   \big(A+\rho_n(t)\big)\cdot
   \frac{\mathbf P_x\big(\tau_0^-\in t+\varDelta_{\delta+1/n}\big)}
   {\mathbf P_x\big(\tau_0^-\in t+\varDelta_\delta\big)}.
   \eeqnn
   In view of Theorem~\ref{MainThm01}, the monotonicity of $\nu(x+\varDelta_\delta)$ in $\delta$ and \eqref{eqn.30}, this yields
   \beqnn
   \lefteqn{\limsup_{t\to\infty}
   	\sup_{\delta\in[\delta_0,\infty]}
   	\mathbf E_x\Big(f(X_{\tau_0^--},-X_{\tau_0^-})
   	\,\Big|\,\tau_0^-\in t+\varDelta_\delta\Big)}
   \cr\cr
   \ar\leq\ar
   A\cdot\limsup_{t\to\infty}
   \sup_{\delta\in[\delta_0,\infty]}
   \frac{\mathbf P_x\big(\tau_0^-\in t+\varDelta_{\delta+1/n}\big)}
   {\mathbf P_x\big(\tau_0^-\in t+\varDelta_\delta\big)}
   \cr\cr
   \ar=\ar
   A\cdot\limsup_{t\to\infty}
   \sup_{\delta\in[\delta_0,\infty]}
   \frac{\nu\big(\beta(t+\varDelta_{\delta+1/n})\big)}
   {\nu\big(\beta(t+\varDelta_\delta)\big)}
   \cr\cr
   \ar\leq\ar
   A\cdot\Big(
   1+\limsup_{t\to\infty}
   \sup_{\delta\in[\delta_0,\infty]}
   \frac{\nu\big(\beta(t+\delta)+\varDelta_{\beta/n}\big)}
   {\nu\big(\beta t+\varDelta_{\beta\delta_0}\big)}
   \Big)
   \leq
   A\cdot\Big(1+\frac{1}{n\delta_0}\Big).
   \eeqnn
   On the other hand, an analogous argument gives
   \beqnn
   \lefteqn{\liminf_{t\to\infty}
   	\inf_{\delta\in[\delta_0,\infty]}
   	\mathbf E_x\Big(f(X_{\tau_0^--},-X_{\tau_0^-})
   	\,\Big|\,\tau_0^-\in t+\varDelta_\delta\Big)}\cr\cr
   \ar\geq\ar 
   A\cdot\liminf_{t\to\infty}
   \inf_{\delta\in[\delta_0,\infty]}
   \frac{\nu\big(\beta(t+\varDelta_{\delta-1/n})\big)}
   {\nu\big(\beta(t+\varDelta_\delta)\big)}
   \geq
   A\cdot\Big(1-\frac{1}{n\delta_0}\Big).
   \eeqnn
   Combining these two bounds and then letting $n\to\infty$ yields the uniform convergence \eqref{eqn.MainThm06}. 
   \qed

 \bibliographystyle{plain}
 \bibliography{Reference}
 
 \end{document}